\documentclass[10pt,twoside]{article}
\usepackage{amssymb,amsmath}
\usepackage{amsfonts}

\usepackage[utf8]{inputenc}

\usepackage{makeidx}

\usepackage{epsfig,psfrag,wrapfig}

\numberwithin{equation}{section}

\begin{document}

\title{The theory of Wiener--It\^o integrals in vector valued
Gaussian stationary random fields. \\ Part I}

\author{P\'eter Major \\
Alfr\'ed R\'enyi Institute of Mathematics \\ 
Budapest, P.O.B. 127 H--1364, Hungary, e-mail:major@renyi.hu \\
{\it Dedicated to the memory of Roland Lvovich Dobrushin  } \\
{\it whose ideas appear in this paper} 
\footnote {The author got support from the Hungarian 
Foundation NKFI--EPR No. K-125569.} }


\maketitle
\noindent
{\it Abstract:}\/ The subject of this work is the multivariate 
generalization of the theory of multiple Wiener--It\^o integrals. 
In the scalar valued case this theory was described by the author 
in 2014. The proofs of the present paper apply the 
technique of that work, but in the proof of some results new 
ideas were needed. The motivation for this study was a result 
in the paper ``Limit theorems for nonlinear functionals of a
stationary Gaussian sequence of vectors'' (1994) by Arcones, 
which contained the multivariate generalization of a non-central 
limit theorem for non-linear functionals of Gaussian stationary 
random fields presented in a paper by R. L. Dobrushin and the 
author. However, the formulation of Arcones' result was incorrect. 
To present it in a correct form the multivariate version of the 
theory explained in my work of 2014 has to be worked  out, 
because the notions introduced in this theory are needed in its 
formulation. This is done in the present paper. In its continuation 
it will be explained how to work out a method  with the help of
the results in this work that enables us to prove non-Gaussian 
limit theorems for non-linear functionals of vector valued 
Gaussian stationary random fields. The right version of Arcones' 
result presented also in the introduction of this work will be 
formulated and proved with its help in a future paper of mine. 

\section{Introduction. An overview of the results.}

Let $X(p)=(X_1(p),\dots,X_d(p))$, $p\in{\mathbb Z}^\nu$, where 
${\mathbb Z}^\nu$ denotes the lattice points with integer 
coordinates in the $\nu$-dimensional Euclidean space 
${\mathbb R}^\nu$, be a $d$-dimensional real valued Gaussian 
stationary random field with expectation $EX(p)=0$, 
$p\in{\mathbb Z}^\nu$. We define the notion of Gaussian 
property of a random field in the usual way, i.e., we 
demand that all finite sets $(X(p_1),\dots,X(p_k))$, 
$p_j\in{\mathbb Z}^\nu$, $1\le j\le k$, be a Gaussian 
random vector, and we call a random field $X(p)$, 
$p\in {\mathbb Z}^\nu$, stationary if for all 
$m\in{\mathbb Z}^\nu$ the random field $X^{(m)}(p)=X(p+m)$, 
$p\in{\mathbb Z}^\nu$, has the same finite dimensional 
distributions as the original random field $X(p)$, 
$p\in{\mathbb Z}^\nu$. In most works only the case $\nu=1$ 
is considered, but since we can prove our results without 
any difficulty for stationary random fields with arbitrary 
parameter $\nu\ge1$ we consider such more general models.

Our goal is to work out a good calculus which provides such a 
representation of the non-linear functionals of our vector
valued Gaussian stationary random field which helps us 
in the study of limit theorems for such functionals. 
To understand what kind of limit theorems we have in mind
take the following example.

Let us have a function $H(x_1,\dots,x_d)$ of $d$ variables,
and define with the help of a $d$-dimensional vector valued 
Gaussian stationary random field 
$$
X(p)=(X_1(p),\dots,X_d(p)), \quad p\in{\mathbb Z}^\nu,  
$$
and this function the
random variables $Y(p)=H(X_1(p),\dots,X_d(p))$ for all 
$p\in{\mathbb Z}^\nu$. Let us introduce for all $N=1,2,\dots$
the normalized sum
\begin{equation}
S_N=A_N^{-1}\sum_{p\in B_N} Y(p) \label{1.1} 
\end{equation}
with an appropriate norming constant $A_N>0$, where
\begin{equation}
B_N=\{p=(p_1,\dots,p_\nu)\colon\; 0\le p_k<N \textrm { for all } 
1\le k\le\nu\}. \label{1.2}
\end{equation}
We are interested in a limit theorem for these normalized
sums $S_N$ with an appropriate norming constant $A_N$ as $N\to\infty$.
In particular, we want to know when we get a classical central 
limit theorem with the natural normalization $A_N=N^{\nu/2}$ and 
when appear new kind of limit theorems. These questions were 
studied in the special scalar valued case $d=1$ in 
papers~\cite{2} and~\cite{5}. Arcones investigated the 
multivariate generalization of the results in these papers.

He proved the multivariate version of the result in paper~\cite{2}
which states that if the covariance function of the underlying
Gaussian field tends to zero sufficiently fast at infinity, and the 
function $H(x_1,\dots,x_d)$ has some nice properties, then the central
limit theorem holds with the classical normalization. (He considered
only the case $\nu=1$, but this restriction has no great importance.)
In Theorem~6 of his paper he also formulated a result about 
a non-central limit theorem under appropriate conditions. But there 
are some serious problems with that result. Arcones wanted to prove 
a multivariate generalization of the result in paper~\cite{5}, but 
to do this he should have solved some problems whose discussion he 
omitted.

The Gaussian limit theorem can be proved in the multivariate case
by means of a natural generalization of the method in paper~\cite{2},
or one can apply some more powerful new method, (see for example~\cite{13}),
but in the proof of the multivariate generalization of the non-central
limit theorem~6 in paper~\cite{1} some new problems appear whose 
solution demands hard work.   

The first problem is related to the formulation of the result. In 
paper~\cite{5} the limit distribution is presented by means of a 
multiple Wiener-It\^o integral with respect to the random spectral
measure of a one-dimensional stationary (generalized) Gaussian 
random field. This random integral was introduced in the paper of 
Dobrushin~\cite{4}, and it is explained in more detail in my
Lecture Note~\cite{9}. But this notion was worked out in
Dobrushin's paper only for scalar valued random fields, and the 
limit distribution in Theorem~6 of Arcones' paper is presented with 
the help of Wiener--It\^o integrals with respect to random spectral 
measures corresponding to vector valued stationary Gaussian random 
fields. Such integrals were not defined before, and their definition 
is far from trivial. The goal of the present paper is to fill this
gap. Here the multivariate random spectral measures will be
introduced together with the multiple Wiener--it\^o integrals
with respect to them, and their most important properties will be
proved. This is needed for the right formulation and proof 
of Arcones' result. I shall formulate the right version of this 
result in the introduction of this paper, but its proof will be 
given only in paper~\cite{12} with the help of the results in 
this work and its continuation~\cite{11}.

To understand what kind of problems we meet in this paper 
let us first consider briefly how the theory of  Wiener--It\^o 
integrals was worked out for scalar valued random fields 
by It\^o in~\cite{8} and Dobrushin in~\cite{4}.

It\^o considered a Gaussian random field in~\cite{8} whose 
elements could be expressed as random integrals with respect 
to a Gaussian orthogonal random measure. He also defined multiple 
random integrals (called later Wiener--It\^o integrals in the 
literature) with respect to this orthogonal random measure, 
and expressed all square integrable random variables 
measurable with respect to the $\sigma$-algebra generated by 
the elements of the Gaussian orthogonal random measure as a
sum of such multiple integrals. The introduction of this integral 
turned out to be useful, because it helped in the study of 
non-linear functionals of the Gaussian random field defined 
by means of this integral. In particular, It\^o found a very 
useful relation, called It\^o's formula in the literature, between 
the multiple random integrals he defined and Hermite polynomials. 

Later Dobrushin worked out a version of this theory in~\cite{4},  
where he studied non-linear functionals of a stationary 
Gaussian random field. In such a random field a spectral and a 
random spectral measure can be defined in such a way that the 
elements of the stationary Gaussian random field can be 
expressed in a special form of (one-fold) random integrals
with respect to the random spectral measure. These random
integrals can be considered as the Fourier transforms of the 
random spectral measure. Dobrushin defined also multiple random 
integrals with respect to this random spectral measure, and 
studied their properties. He proved that these random 
integrals defined with respect to the random spectral measure 
have similar properties as the multiple integrals 
introduced by It\^o. In particular, he proved It\^o's 
formula for this new type of random integrals. This enabled
him to express all square integrable random variables 
measurable with respect to the $\sigma$-algebra generated by 
the elements of the original stationary Gaussian random field 
as a sum of multiple random integrals with respect to 
the random spectral measure. He also found a simple and 
useful formula for the calculation of the shift transforms 
of a random variable which is presented as a sum of multiple
random integrals. With the help of these results the normalized 
random sums $S_N$ defined in~(\ref{1.1}) can be expressed in a
simple and useful form if the underlying stationary Gaussian
random field is scalar valued (i.e., $d=1$). This representation 
of the normalized random sums~$S_N$ made possible to prove the 
limit theorems in~\cite{5}. 

We want to prove the generalization of the results in~\cite{5}
for non-linear functionals of vector valued stationary Gaussian
random fields. The first step of this program is to work out
the multivariate version of Dobrushin's theory, and this is
the subject of the present paper.
 
First we have to define the spectral and random spectral 
measure of vector valued stationary Gaussian random fields, 
and this is the subject of Sections~2 and~3. To do this the 
multivariate version of some classical results has to
be proved. In the scalar valued case a spectral measure
can be defined whose Fourier transform is the correlation
function of the stationary random field we are working with.
In the case of a vector valued stationary random field of 
dimension~$d$ the correlation function is a $d\times d$
dimensional matrix valued function. It can be shown that 
there exists a $d\times d$ dimensional matrix valued measure 
on the $d$ dimensional torus $[-\pi,\pi)^d$ for which each 
coordinate of the matrix valued correlation function is 
the Fourier transforms of the corresponding coordinate of 
this matrix valued measure. This measure is called the 
spectral measure of the random field. In the scalar valued 
case, i.e., if $d=1$ the spectral measure is a positive 
measure, while in the vector valued case it is a positive 
semidefinite matrix valued measure. A more detailed description
of these results together with their proofs is given in
Section~2.

In Section~3 the so-called random spectral
measure corresponding to a vector valued stationary
Gaussian random field is defined. It is a vector valued
random measure with the same dimension~$d$ as the
underlying vector valued stationary Gaussian random
field. Its distribution is determined by the spectral
measure of the underlying random field. A random 
integral can be defined with respect to the coordinates
of the random spectral measure, and each coordinate of 
the elements of the underlying vector valued Gaussian 
random field can be expressed by means of an appropriate 
random integral with respect to the corresponding 
coordinate of the random spectral measure. Because of 
the form of this integral this result can be interpreted 
so that the underlying stationary Gaussian random field 
is the Fourier transform of the random spectral measure 
corresponding to it. The construction of the random 
spectral measure and the description of its most 
important properties is given in Section~3.
  
Moreover, we need later the notion of spectral measures 
and random spectral measures corresponding to stationary 
generalized random fields, and they are introduced in 
Section~4. In the main text of this paper a more detailed, 
precise definition of these notions will be given. We 
have to define these objects, because we can formulate 
the limit in the limit theorems we are interested in 
in this paper by means of multiple random integrals with 
respect to the random spectral measures corresponding to 
stationary generalized random fields.

Then I define the multiple Wiener--It\^o integrals with 
respect to the coordinates of a vector valued random
spectral measure in Section~5, and I also prove there their 
most important properties. In Section~6 I prove an important
result, called the diagram formula which enables us to express 
the product of two multiple Wiener--It\^o integrals as the 
sum of appropriately defined multiple Wiener--It\^o 
integrals. The present paper contains these results.

In the continuation of this paper, in work~\cite{11}
I work out the basic tools needed in the proof of
such non-central limit theorems as the multivariate
generalization of the limit theorem in~\cite{5}. First
I prove, with the help of the above mentioned diagram
formula, an important result about the relation between
multiple Wiener--It\^o integrals and Wick polynomials of
Gaussian vectors. Wick polynomials are the several 
dimensional generalizations of Hermite polynomials, and
the result mentioned before is the natural multivariate 
generalization of It\^o's formula. Besides, \cite{11} 
contains a formula that enables us to express the shift 
transforms of a random variable given in the form of a 
sum of multiple random variables in a useful form. These 
results enable us to rewrite the normalized random sums 
$S_N$ defined in~(\ref{1.1}) in a form which helps in 
the study of limit theorems. They enabled me to formulate 
and  prove in~\cite{12} the right version of Theorem~6 
in Arcones' paper~\cite{1}.

Next I briefly describe the right version of Arcones'
non-central limit theorem. In its formulation we consider 
$d$-dimensional stationary Gaussian random fields 
$$
X(p)=(X_1(p),\dots,X_d(p)),\quad EX_j(p)=0 \textrm{ for all } 
1\le j\le\nu \textrm{ and } p\in{\mathbb Z}^\nu,
$$
whose covariance function $r_{j,j'}(p)=EX_j(0)X_{j'}(p)$, 
$1\le j,j'\le d$, $p\in\mathbb Z^\nu$,  is such a matrix 
valued function whose coordinates decrease asymptotically 
polynomially at infinity with some power $0<\alpha<\nu$. 
More generally, this behaviour may be slightly modified 
by multiplication with a slowly varying function. More 
explicitly, we demand that
\begin{equation}
\lim_{T\to\infty}\sup_{p\colon\;p\in{\mathbb Z}^\nu,\,|p|\ge T}
\frac{\left|r_{j,j'}(p)-a_{j,j'}(\frac p{|p|})|p|^{-\alpha}L(|p|)\right|}
{|p|^{-\alpha}L(|p|)}=0 \label{1.3}
\end{equation}
for all $1\le j,j'\le d$, where $0<\alpha<\nu$, $L(t)$, $t\ge1$, 
is a real valued function, slowly varying at infinity, 
bounded in all finite intervals, and $a_{j,j'}(t)$ is a 
real valued continuous function on the unit sphere 
${\cal S}_{\nu-1}=\{x\colon\;x\in {\mathbb R}^\nu,\;|x|=1\}$,
and the identity $a_{j',j}(x)=a_{j,j'}(-x)$ holds 
for all $x\in{\cal S}_{\nu-1}$ and $1\le j,j'\le d$.

For the sake of simpler discussion we also demand that
\begin{equation}
EX_j^2(0)=1 \textrm{ for all } 1\le j\le d, \textrm{ and } 
EX_j(0)X_{j'}(0)=0
\textrm{ if } j\neq j', \;\;1\le j,j'\le d. \label{1.4}
\end{equation}
This is not an essential restriction, as it is explained 
in~\cite{12}.

We want to describe the limit behaviour of some non-linear 
functionals of such a random field. To do this first we 
describe the asymptotic behaviour of its spectral measure. 
To formulate such a result let us introduce the following 
notation.

\medskip
Given a vector valued stationary random field 
$X(p)=(X_1(p),\dots,X_d(p))$, $p\in{\mathbb Z}^\nu$, 
with expectation zero and covariance function 
$r_{j,j'}(p)=EX_j(0)X_{j'}(p)$, $1\le j,j'\le d$, $p\in \mathbb Z^\nu$
that satisfies relation~(\ref{1.3}), let us consider its matrix valued 
spectral measure $G=(G_{j,j'})$,
$1\le j,j'\le d$, on the torus $[-\pi,\pi)^\nu$. Take its rescaled 
version $G^{(N)}=(G^{(N)}_{j,j'}$, $1\le j,j'\le d$,
\begin{equation}
G^{(N)}_{j,j'}(A)=\frac{N^\alpha}{L(N)}G_{j,j'}\left(\frac AN\right),
\quad A\in{\cal B}^\nu, \quad N=1,2,\dots,\;\;1\le j,j'\le d, \label{1.5}
\end{equation}
concentrated on $[-N\pi,N\pi)^\nu$ for all $N=1,2,\dots$, 
where ${\cal B}^\nu$ denotes the $\sigma$-algebra of the 
Borel measurable sets on ${\mathbb R}^\nu$. In the next result
we give the limit of the matrix valued measures $G^{(N)}$, as 
$N\to\infty$. Since the coordinates of the matrices $G^{(N)}$ are 
non-probability measures and their limits are non-finite measures, 
we have to introduce the right form of convergence which will
be applied in the limit theorem we shall describe. In 
paper~\cite{12} the so-called vague convergence of complex 
measures are defined, (more precisely its definition is recalled). 
In this definition also the notion of complex measures with 
locally finite measures appear whose definition is explained 
in Section~4 of this paper. This notion was introduced, because
they are needed in the study of spectral measures of stationary
generalized fields, and we want to work with such objects. In the 
presentation of the limit theorem I want to discuss we  need the 
result of Proposition~1.1 of \cite{12} whose formulation applies 
the above notions. This Proposition~1.1 agrees with the following 
result.

\medskip\noindent
{\bf Proposition~1.1.} {\it Let $G=(G_{j,j'})$ be the matrix
valued spectral measure of a $d$-dimensional vector valued 
stationary random field whose covariance function 
$r_{j,j'}(p)$  satisfies relation~(\ref{1.3}) with some parameter
$0<\alpha<\nu$. Then for all pairs $1\le j,j'\le d$
the sequence of complex measures~$G^{(N)}_{j,j'}$ 
defined in~(\ref{1.5}) with the help of the complex measure 
$G_{j,j'}$ tends vaguely to a complex measure~$G^{(0)}_{j,j'}$ 
on ${\mathbb R}^\nu$ with locally finite total variation. These 
complex measures~$G^{(0)}_{j,j'}$, $1\le j,j'\le d$, 
have the homogeneity property
\begin{equation}
G^{(0)}_{j,j'}(A)=t^{-\alpha}G^{(0)}_{j,j'}(tA) 
\quad \textrm{for all bounded } A\in{\cal B}^\nu, 
\;1\le j,j'\le d,\textrm{ and } t>0.  \label{1.6}
\end{equation}
The complex measure $G^{(0)}_{j,j'}$ with locally 
finite variation is determined by the number $0<\alpha<\nu$
and the function $a_{j,j'}(\cdot)$ on the unit sphere~$S_{\nu-1}$ 
introduced in formula~(\ref{1.3}).

There exists  a vector valued Gaussian stationary generalized
random field on ${\mathbb R}^\nu$ with that matrix valued spectral 
measure  $(G^{(0)}_{j,j'})$, $1\le j,j'\le d$, whose coordinates 
are the above defined complex measures $G^{(0)}_{j,j'}$, 
$1\le j,j'\le d$.} 

\medskip
In the non-central limit theorem I shall describe the limit of 
random variables $S_N$  defined by formulas~(\ref{1.1}) 
and~(\ref{1.2}) with the help of 
a vector valued stationary Gaussian random field whose correlation 
function satisfies relations~(\ref{1.3}) and~(\ref{1.4}) and an
appropriate norming constant $A_N$. To give a 
complete definition of these random variables we must tell what
kind of functions $H(x_1,\dots,x_d)$ we apply in their
definition. I shall choose functions of the following form in this
definition. $H(x_1,\dots,x_d)$ depends on a previously fixed
constant $k$, and it has the form
\begin{equation}
H(x_1,\dots,x_d)=\sum_{\substack{(k_1,\dots, k_d),\; k_j\ge 0,\;1\le j\le d,\\
k_1+\cdots+k_d=k}}
c_{k_1,\dots,k_d}H_{k_1}(x_1)\cdots H_{k_d}(x_d) \label{1.7} 
\end{equation}
with some coefficients $c_{k_1,\dots,k_d}$, where $H_k(\cdot)$ 
denotes the $k$-th Hermite polynomial with leading coefficient~1.

The limit distribution of the above introduced random variable 
$S_N$ is de\-sribed in Theorem 1.2A of~\cite{12}. This theorem 
is written down in the following Theorem~1.2. The limit in this 
result is presented by means of a multiple Wiener--It\^o 
integral with respect to the random spectral measure 
corresponding to the matrix valued spectral measure 
$(G^{(0)}_{j,j'})$, $1\le j,j'\le d$, which appeared in 
Proposition~1.1. Let me remark that because of the 
homogeneity property~(\ref{1.6}) of this measure 
$G^{(0)}_{j,j}(\mathbb R^\nu)=\infty$ for any $1\le j\le d$. 
Hence this matrix valued spectral measure can be defined 
only as the spectral measure of a generalized and not as 
the spectral measure of an ordinary vector valued 
stationary random field. 

\medskip\noindent
{\bf Theorem 1.2.} {\it Fix some integer $k\ge1$,
and let $X(p)=(X_1(p),\dots,X_d(p))$, 
$p\in{\mathbb Z}^\nu$, be a vector valued Gaussian 
stationary random field whose covariance function 
$r_{j,j'}(p)=EX_j(0)X_{j'}(p)$, $1\le j,j'\le d$, 
$p\in{\mathbb Z}^\nu$, satisfies relation~(\ref{1.3}) with 
some $0<\alpha<\frac\nu k$ and relation~(\ref{1.4}). Let 
$H(x_1,\dots,x_d)$ be a function of the form given 
in~(\ref{1.7}) with the parameter $k$ we have fixed in the
formulation of this result. Define the random variables 
$Y(p)=H(X_1(p),\dots,X_d(p))$ for all 
$p\in{\mathbb Z}^\nu$ together with their normalized 
partial sums 
$$
S_N=\frac1{N^{\nu-k\alpha/2}L(N)^{k/2}}\sum_{p\in B_N}Y(p),
$$ 
where the set $B_N$ was defined in~(\ref{1.2}). These random 
variables $S_N$, $N=1,2,\dots$, satisfy the following 
limit theorem.

Let $Z_{G^{(0)}}=(Z_{G^{(0)},1},\dots,Z_{G^{(0)},d})$ be a vector 
valued random spectral measure which corresponds to the 
matrix valued spectral measure $(G^{(0)}_{j,j'})$, 
$1\le j,j'\le d$, defined in  Proposition~1.1 with
the help of the matrix valued spectral measure $G=(G_{j,j'})$,
corresponding the covariance function $r_{j,j'}(p)$ we
are working with. Then the sum of multiple Wiener--It\^o 
integrals
\begin{eqnarray}
S_0&=&\sum_{\substack{(k_1,\dots, k_d),\; k_j\ge 0,\;1\le j\le d, \\
k_1+\cdots+k_d=k}} c_{k_1,\dots,k_d}
\int \prod_{l=1}^\nu\frac{e^{i(x_1^{(l)}+\cdots+x_k^{(l)})}-1}
{i(x^{(l)}_1+\cdots+x_k^{(l)})} \label{1.8} \\
&&\qquad\qquad\qquad Z_{G^{(0)},{j(1|k_1,\dots,k_d)}}(\,dx_1)\dots 
Z_{G^{(0)},{j(k|k_1,\dots,k_d)}}(\,dx_k)
\nonumber
\end{eqnarray}
exists. (These Wiener--It\^o integrals are defined in 
Section~5 of this paper.) Here we use the notation  
$x_p=(x_p^{(1)},\dots,x_p^{(\nu)})$, $p=1,\dots,k$, 
and define the indices $j(s|k_1,\dots,k_d)$,
$1\le s\le k$, as  $j(s|k_1,\dots,k_d)=r$ if 
$\sum_{u=1}^{s-1} k_u< r\le \sum_{u=1}^s k_u$, 
$1\le s\le k$. (For $s=1$ we apply the notation 
$\sum_{u=1}^{0} k_u=0$ in the definition of 
$j(1|k_1,\dots,k_d)$.) The normalized sums 
$S_N$ converge in distribution to the random 
variable $S_0$ defined in~(\ref{1.8}) as $N\to\infty$.}

\medskip
The indexation of the terms $Z_{G^{(0)},{j(s|k_1,\dots,k_d)}}(\,dx_s)$ 
in formula (\ref{1.8}) can be explained in a simpler way. In 
the first $k_1$ arguments $x_1,\dots,x_{k_1}$ we write 
$Z_{G^{(0)},1}(\,dx_s)$, $1\le s\le k_1$, in the next $k_2$ terms 
we write $Z_{G^{(0)},2}(\,dx_s)$, $k_1+1\le s\le k_1+k_2$, and 
so on. In the last $k_d$ terms we write $Z_{G^{(0)},d}(\,dx_s)$,
$k_1+\cdots+k_{d-1}+1\le s\le k$.

\medskip
Actually a more general limit theorem is also proved 
in~\cite{12}, but its proof is based on the result of 
Theorem~1.2. It is worth comparing Theorem~1.2 
with its scalar valued version (i.e., with the result in the 
case $d=1$ proved in~\cite{5}).

In paper~\cite{5} a result similar to Theorem~1.2 is proved 
in the scalar valued case. In that result $CH_k(x)$, $C\neq0$, 
i.e., the $k$-th Hermite polynomial multiplied with a non-zero 
coefficient~$C$ plays the same role as the function $H(\cdot)$
defined in~(\ref{1.7}) in Theorem~1.2, and the condition 
$k\alpha<\nu$ has to be imposed. The limit is given by 
formula~(\ref{1.8}) in the case $d=1$ with $H(x)=CH_k(x)$. 
Let me remark that the Wick polynomials, i.e., the multivariate 
generalizations of Hermite polynomials appeared in Theorem~1.2 
in a hidden way. (See for example Section~2 of \cite{9} for the 
definition of Wick polynomials.) Indeed, the random variables 
$Y(p)=H(X_1(p),\dots,X_d(p))$, $p\in\mathbb Z^\nu$, defined 
with the help of the function $H(\cdot)$ introduced in 
formula~(\ref{1.7})  are Wick polynomials of order $k$ because 
of the relation~(\ref{1.4}). (See Corollary~2.3 in~\cite{9}.) 
This indicates that the role of Hermite polynomials in results 
about scalar valued stationary Gaussian random fields is taken 
by Wick polynomials in the their vector valued counterparts. 
The next results also show such a correspondence. 

The limit theorem in \cite{5} remains valid if we replace the
function $CH_k(x)$ in it with such a function $H(x)$ whose
expansion with respect to the Hermite polynomials contains 
only terms  $H_{k'}(x)$ of order $k'\ge k$, and the term 
$H_k(x)$ has a non-zero coefficient. The limit is the same as 
in the case when we take only the first term 
$\textrm{const.}H_k(x)$ in the expansion of the function $H(x)$.
Similarly, Theorem~1.2 formulated above in the multivariate 
case remains valid if such a random random variable 
$H(X_1(0),\dots,X_d(0))$ is taken whose expansion  with respect 
to Wick polynomials starts with a non-zero Wick polynomial of 
order~$k$, and $k\alpha<\nu$. The limit does not change if we
take only the term of order $k$ of $H(X_1(0),\dots,X_d(0))$ in
this expansion.

Let me finally remark that the Theorem holds only under the 
condition $k\alpha<\nu$. In the case $k\alpha>\nu$ the central 
limit theorem holds for $S_N$ with the usual norming constant 
$A_N=N^{\nu/2}$. This follows from a slight generalization of the 
(correct) results in Arcones' paper~\cite{1}. In the boundary 
case $k\alpha=\nu$ the central limit theorem holds again for 
$S_N$, but in this case the norming constant may have the form 
$A_N=N^{\nu} L'(N)$ with a slowly varying function $L'(N)$ tending
to infinity as $N\to\infty$. Let me also remark that the
definition of the limit distribution in Theorem~1.2 given in 
formula~(\ref{1.8}), is meaningful only for $k\alpha<\nu$. This 
formula contains a multiple Wiener--It\^o integral, and we have 
to check whether this Wiener--It\^o integral is meaningful. It 
is explained at the beginning of Section~5 that the multiple 
Wiener--It\^o integrals  are defined only with such kernel 
functions that satisfy an integrability condition. (This 
condition is formulated in property~(b) in the definition of 
a class of functions ${\cal K}_{n,j_1,\dots,j_n}$.) It can be 
seen that the Wiener--It\^o integral appearing in 
formula~(\ref{1.8}) is meaningful if $k\alpha<\nu$, because 
this integrability condition is satisfied in this case. On 
the other hand, this integral cannot be defined if 
$k\alpha\ge\nu$, because in this case this integrability 
condition is violated. 

\subsection{A more detailed description of the results.}

Next I give a more detailed overview about the results of
this paper.

First I characterize the distribution 
of the vector valued Gaussian stationary random fields 
$X(p)=(X_1(p),\dots,X_d(p))$, $p\in{\mathbb Z}^\nu$, 
with expectation zero. This is the subject of the 
second section of this work. Because of the Gaussian 
and stationary property of such a random field its 
distribution is determined by the correlation 
function $r_{j,j'}(p)=EX_j(0)X_{j'}(p)$ for all 
$1\le j,j'\le d$ and $p\in{\mathbb Z}^\nu$. We are 
interested in the description of those functions 
$r_{j,j'}(p)$ that can appear as the correlation 
function of a vector valued stationary random field.

In the scalar valued case a well-known result solves
this problem. The correlation function $r(p)=EX(0)X(p)$,
$p\in{\mathbb Z}^\nu$, of a stationary field $X(p)$, 
$p\in{\mathbb Z}^\nu$, can be represented in a unique
way as the Fourier transform of a spectral measure,
and the spectral measures can be characterized. Namely,
we call the finite (non negative), even measures on the 
torus $[-\pi,\pi)^\nu$ spectral measures. For any 
correlation function $r(p)$ of a stationary field 
there is a unique spectral measure $\mu$ such
that $r(p)=\int e^{i(p,x)}\mu(\,dx)$ for all 
$p\in{\mathbb Z}^\nu$, and for all spectral
measures~$\mu$ there is a (Gaussian) stationary random
field whose correlation function equals the Fourier
transform of this spectral measure~$\mu$.

In Section~2 we prove a similar result for vector 
valued stationary random fields. In the case of a 
vector valued Gaussian stationary random field 
$X(p)=(X_1(p),\dots,X_d(p))$, $p\in{\mathbb Z}^\nu$, 
we have for all pairs of indices $(j,j')$, 
$1\le j,j'\le d$, a unique complex  
measure $G_{j,j'}$ on the torus $[-\pi,\pi)^\nu$ with 
finite total variation such that 
$r_{j,j'}(p)=EX_j(0)X_{j'}(p)=\int e^{i(p,x)}G_{j,j'}(\,dx)$
for all $p\in{\mathbb Z}^\nu$. This can be interpreted
so that the correlation function $r_{j,j'}(p)$, 
$1\le j,j'\le d$, $p\in{\mathbb Z}^\nu$, is the 
Fourier transform of a matrix valued measure 
$(G_{j,j'})$, $1\le j,j'\le d$, on the torus 
$[-\pi,\pi)^\nu$. We want to give, similarly to the
scalar valued case, a complete description of those
matrix valued measures on the torus $[-\pi,\pi)^\nu$ 
for which the correlation function of a vector valued  
Gaussian stationary random field can be represented 
as its Fourier transform. Such matrix valued measures
will be called matrix valued spectral measures.

As I have mentioned, the coordinates of a matrix
valued spectral measure are complex
measures with finite total variation. 
The scalar valued counterpart of this condition
is the condition that the spectral measure of
a scalar valued stationary random field must be 
finite. Another important property of a matrix 
valued spectral measure is that it must be 
positive semidefinite. The meaning of this 
property is explained before the formulation
of Theorem~2.2, and Lemma~2.3 gives a different,
equivalent characterization of this property. 
Let me remark that in the scalar valued case 
the spectral measure must be a measure (and 
not only a complex measure),
and this fact corresponds to the above property 
of matrix valued spectral measures. Finally, a 
matrix valued spectral measure must be even. 
This means that its coordinates are even, i.e., 
for all $1\le j,j'\le d$ and measurable sets 
$A$ on the torus 
$G_{j,j'}(-A)=\overline{G_{j,j'}(A)}$, where
the overline indicates complex conjugate. 

Theorem~2.2 states that the above properties 
characterize the matrix valued spectral measures.
Let me remark that there are papers (see for example 
\cite{3}, \cite{7} or \cite{14}) 
containing the above 
results, although in a slightly different
formulation, at least in the case $\nu=1$. 
Nevertheless, I worked out their proof, 
since I applied a different method, which 
is used also in the later part of the paper.

In Section 3, I introduce the vector valued random 
spectral measures corresponding to a matrix valued 
spectral measure $(G_{j,j'})$, $1\le j,j'\le d$.
To do this first I consider a vector valued stationary 
Gaussian random field  $X(p)=(X_1(p),\dots,X_d(p))$, 
$p\in{\mathbb Z}^\nu$, with  spectral measure 
$(G_{j,j'})$, $1\le j,j'\le d$, and show that a
vector valued random measure 
$Z_G=(Z_{G_1},\dots,Z_{G_d})$ can be defined on the
measurable subsets $A\subset[-\pi,\pi)^\nu$ of 
the torus which have some nice properties.
A random integral can be defined with respect
to the coordinates of this random measure, and the 
coordinates $X_j(p)$, $1\le j\le d$,  
$p\in\mathbb Z^\nu$, of the random field $X(p)$ 
can be expressed as the Fourier transforms of 
the appropriate coordinate $Z_{G_j}$ of this
random measure. More explicitly, 
$X_j(p)=\int e^{i(p,x)}Z_{G,j}(\,dx)$ for all 
$p\in{\mathbb Z}^\nu$ and $1\le j\le d$.
I remark that the random variables $Z_{G,j}(A)$, 
$1\le j\le d$, $A\subset[-\pi,\pi)^\nu$, are 
complex valued. 

I have listed some properties of this random
measure $(Z_{G,1},\dots,Z_{G,d})$. These properties
determine its distribution, and they depend 
only on the spectral measure $(G_{j.j'})$, 
$1\le j,j'\le d$, of the underlying random field
$X(p)$, $p\in\mathbb Z^\nu$. We shall call the 
vector valued random measures with these properties 
a vector valued random spectral measure 
corresponding to the matrix valued spectral measure 
$(G_{j,j'})$, $1\le j,j'\le d$. We can prove 
that the Fourier transform of all vector 
valued random spectral measures corresponding 
to a matrix valued spectral measure can be defined,
and it is a vector valued Gaussian stationary random 
field with this matrix valued spectral measure.

Besides the above results  I also proved some 
important properties of the random integrals 
with respect to a vector valued spectral 
measure in Section~3. I characterized those
functions which can be integrated with respect 
to these random spectral measure, and also 
described those functions whose integrals are 
real valued random variables. In particular, I 
proved that if a vector valued Gaussian stationary 
random field $X(p)=(X_1(p),\dots,X_d(p))$, 
$p\in{\mathbb Z}^\nu$, is given, we fix some 
parameter $1\le j\le d$, and take the real Hilbert 
space consisting of the closure of finite linear 
combinations $\sum_k c_kX_j(p_k)$ with real number 
valued coefficient $c_k$ in the Hilbert space of 
square integrable random variables, then
each element of this Hilbert space can be 
expressed as the integral of a function on
the torus~$[-\pi,\pi)^\nu$ with respect to 
the random spectral measure $Z_{G,j}$. The
functions taking part in the representation
of this Hilbert space also constitute a real 
Hilbert space. A more detailed formulation of 
this result is given in Lemma~3.2.

It may be worth discussing the relation of 
the results in Section~3 to their scalar 
valued correspondents. The results about
the existence of random spectral measures 
for scalar valued Gaussian stationary 
random fields give a great help in proving 
the results in Section~3. In particular, 
these results provide the definition of the 
random spectral measures $Z_{G,j}$, and 
determine their distribution for all 
$1\le j\le d$. The definition of $Z_{G,j}$, 
and the properties determining its 
distribution depend only on the 
measure~$G_{j,j}$. On the other hand,
we had to carry out some additional work 
to prove those properties of a vector 
valued spectral random measure which 
determine the joint distribution of their 
coordinates. The non-diagonal elements  
$G_{j,j'}$ with $j\neq j'$ of the matrix 
valued spectral measure $(G_{j,j'})$,
$1\le j,j'\le d$, appear at this point of 
the investigation.

The fourth section deals with a special 
subject, and our motivation to study it 
demands some explanation. Here we consider 
vector valued Gaussian stationary generalized 
random fields.

We could have considered the continuous time version of 
vector valued stationary random fields where the parameter 
set is $t\in {\mathbb R}^\nu$ and not $p\in{\mathbb Z}^\nu$. 
Here we did not discuss such models, we have considered 
instead vector valued Gaussian stationary generalized 
random fields. This means a set of random vectors  
$(X_1(\varphi),\dots,X_d(\varphi))$ with some nice
properties which are indexed by an appropriately chosen
class of functions. The precise definition of this notion 
is given in Section~4. We have constructed a large 
class of Gaussian stationary generalized random fields, 
presented their matrix valued spectral measures, and 
constructed the vector valued random spectral measures
corresponding to them. In \cite{9} the notion of Gaussian 
stationary generalized random fields was introduced
and investigated in the scalar valued case. Some
useful results were proved there. It was shown, with 
the help of some important results of Laurent Schwartz 
about distributions (generalized functions), that in 
the scalar valued case the class of Gaussian, stationary 
generalized random fields constructed in such a way as 
it was done in the present paper contains all Gaussian 
stationary generalized random fields. (Here I 
consider two random fields the same if their finite 
dimensional distributions agree.) Similarly, it is very 
likely that also in the multivariate case all 
stationary generalized Gaussian random fields can be 
constructed by the method described in this paper. But  
I did not study this question, because I was interested 
in a different problem.

Although the theory of generalized random fields is 
an interesting subject in itself, I investigated 
it for a different reason. I was interested in 
the matrix valued spectral measures of vector
valued Gaussian stationary generalized random
fields and the vector valued random spectral 
measures corresponding to them and not in the
Gaussian, stationary generalized random fields 
which were needed for their construction. They 
behave similarly to the analogous objects
corresponding to (non-generalized) Gaussian
stationary random fields. We can work with them
in the same way. Nevertheless, there is a difference
between these new spectral and random spectral 
measures and their previously defined counterparts 
which is very important for us. Namely, the coordinates 
of a matrix valued spectral measure corresponding 
to a non-generalized random field are complex  
measures with finite total variation, while 
in the case of generalized random fields the matrix 
valued spectral measures need not satisfy this condition. 
It is enough to demand that the corresponding matrix 
valued measures have locally finite total variation, 
and the matrix valued spectral measures are semidefinite
matrix valued measures with moderately increasing
distribution at infinity. (The definition of these
notions is contained in Section~4.)

The above facts mean that we can work with a much
larger class of random spectral measures after the
introduction of Gaussian stationary generalized random
fields and random spectral measures corresponding
to them. This is important for us, because in the limit
theorems we are interested in the limit can be expressed
by means of multiple Wiener--It\^o integrals with
respect to random spectral measures  constructed with 
the help of vector valued Gaussian stationary 
generalized random fields. Theorem~1.2 discussed in this
introduction is an example for such a limit theorem. 

Sections 2---4 contain the main results about the  
linear functionals of vector valued Gaussian stationary 
random fields. They are also needed in the study of 
their non-linear functionals , and this is the subject 
of Sections~5 and~6. The results of these sections help 
us to work out some tools which are useful in the 
study of limit theorems with a new type of 
non-Gaussian limit.

In Section 5 multiple Wiener--It\^o integrals are defined
with respect to the coordinates of a vector valued random
spectral measure $(Z_{G,1},\dots,Z_{G,d})$. We define for
all numbers $n=1,2,\dots$, and parameters 
$j_1,\dots,j_n$ such that $1\le j_k\le d$ for all 
$1\le k\le n$ and all functions $f\in{\cal K}_{n,j_1,\dots,j_n}$,
where ${\cal K}_{n,j_1,\dots,j_n}$ is a real Hilbert space 
defined in Section~5, an $n$-fold Wiener--It\^o integral
$$
I_n(f|j_1,\dots,j_n)=
\int f(x_1,\dots,x_n)Z_{G,{j_1}}(\,dx_1)\dots Z_{G,{j_n}}(\,dx_n),
$$
and prove some of its basic properties. The definition and 
proofs are very similar to the definition and proofs in 
scalar valued case, only we have to apply the properties 
of vector valued random spectral measures.

There is one point where we have a weaker estimate than
in the scalar valued case. We can give an upper bound on 
the second moment of a multiple Wiener--It\^o integral 
with the help of the $L_2$ norm of the kernel function 
of this integral in the way as it is formulated in
formula~(\ref{5.6}), but we can state here only an inequality 
and not an equality. The behaviour of Wiener--It\^o 
integrals with respect to a scalar valued random spectral 
measure is different. If we integrate in this case a 
symmetric function, and we may restrict our attention to 
such integrals, then we have equality in the corresponding 
relation. This weaker form of the estimate~(\ref{5.6})  
has the consequence that in certain problems we can get 
only weaker results for Wiener--It\^o integrals with 
respect to the coordinates of a vector valued random 
spectral measure than for Wiener--It\^o integrals with  
respect to scalar valued random spectral measures. But 
this will cause no serious problem in our study about
multiple Wiener--It\^o integrals with respect to vector
valued random spectral measures. 

Multiple Wiener--It\^o integrals were introduced
in order to express a large class of random 
variables with their help. More precisely, we are 
interested in the following problem. Let us have a 
vector valued Gaussian stationary random field 
$X(p)=(X_1(p),\dots,X_d(p))$, $p\in{\mathbb Z}^\nu$. 
Their elements can be expressed as the  Fourier 
transforms of a vector valued random spectral 
measure $Z_G=(Z_{G,1},\dots,Z_{G,d})$. Let us consider
the real Hilbert space ${\cal H}$ defined in the
second paragraph of Section~5 with the help
of this vector valued stationary Gaussian random field.
We would like to express the elements of this Hilbert 
space in the form of a sum of multiple Wiener--It\^o 
integrals with respect to the coordinates of the 
vector valued spectral measure~$Z_G$. This problem
together with the study of a theory useful in the
investigation of limit theorems for non-linear 
functionals of vector valued stationary Gaussian 
random fields will be the subject of the second 
part of this work~\cite{11}. But to carry out
this program we still need the proof of an important
result about multiple Wiener--It\^o integrals
discussed in Section~6 of this work. 

In Section~6 I formulate and prove the multivariate 
version of a classical result. I describe the product 
of two multiple Wiener--It\^o integrals as the sum of 
multiple Wiener--It\^o integrals with respect to the
coordinates of a vector valued random spectral measure.
The formulation and proof of this result is similar
to that of the corresponding result in the scalar 
valued case. In this result we define the kernel 
functions of the Wiener--It\^o integrals appearing 
in the sum expressing the product of two Wiener--It\^o 
integrals with the help of some diagrams. Hence this 
result got the name diagram formula. I wrote down 
the formulation of the diagram formula in the case 
of vector valued random spectral measures in detail. 
On the other hand, I gave only a sketch of its 
proof, because it is actually an adaptation of 
the original proof with a rather unpleasant notation. 
I concentrated on the points which explain why the 
diagram formula has such a form as we claim. 
Besides, I tried to explain those steps of the 
proof where we have to apply some new ideas. I 
hope that the interested reader can reconstruct 
the proof on the basis of these explanations by 
looking at the original proof.

Section 6 also contains a corollary of the diagram
formula, where I formulate this result in a special
case. I formulated this corollary, because in this
work we need only this corollary of the diagram formula.

\section{Spectral representation of vector valued 
stationary random fields}

Let $X(p)=(X_1(p),\dots,X_d(p))$, $p\in{\mathbb Z}^\nu$, where 
${\mathbb Z}^\nu$ denotes the lattice of points with integer 
coordinates in the $\nu$-dimensional Euclidean space 
${\mathbb R}^\nu$, be a $d$-dimensional real valued Gaussian 
stationary random field with expected value $EX(p)=0$, 
$p\in{\mathbb Z}^\nu$. Let us first characterize the 
covariance matrices $R(p)=(r_{j,j'}(p))$, $1\le j,j'\le d$, 
$p\in{\mathbb Z}^\nu$, of this $d$-dimensional stationary 
random field, where
$r_{j,j'}(p)=EX_j(0)X_{j'}(p)=EX_j(m)X_{j'}(p+m)$, 
$1\le j,j'\le d$, $p,m\in{\mathbb Z}^\nu$.

In the case $d=1$ we can characterize the function $R(p)=EX(0)X(p)$,
(in this case $j=j'=1$, so we can omit these indices) as the Fourier
transform of an even, finite (and positive) measure $G$ on the torus 
$[-\pi,\pi)^\nu$, called the spectral measure. We are looking for 
the vector valued version of this result. Before discussing this
problem I recall the definition of the torus $[-\pi,\pi)^\nu$. 

The points of the torus $[-\pi,\pi)^\nu$ are those points 
$x=(x_1,\dots,x_\nu)\in {\mathbb R}^\nu$ for which 
$-\pi\le x_j\le\pi$ for all $1\le j\le\nu$. But if a coordinate 
of $x$ in this set equals $\pi$, then we consider this point 
the same if we replace this coordinate by~$-\pi$. In such 
a way we can identify all points of this set by a point of 
the set $[-\pi,\pi)^\nu\subset {\mathbb R}^\nu$. We define 
the topology on the torus on $[-\pi,\pi)^\nu$ as the topology 
induced by the metric 
$\rho(x,y)=\sum\limits_{j=1}^\nu(|x_j-y_j|\;\textrm{mod}\;2\pi)$
if $x=(x_1,\dots,x_\nu)\in[-\pi,\pi)^\nu$ and 
$y=(y_1,\dots,y_\nu)\in[-\pi,\pi)^\nu$. These properties of the 
torus $[-\pi,\pi)^\nu$ must be taken into account when we 
speak of the set $-A=\{-x\colon\; x\in A\}$ for a set 
$A\subset [-\pi,\pi)^\nu$ or of a continuous function on the 
torus $[-\pi,\pi)^\nu$.

Later we shall speak also about the torus $[-A,A)^\nu$ for arbitrary
$A>0$. This is defined in the same way, only the number $\pi$ is 
replaced by~$A$ in the definition.

\medskip\noindent

It is natural to expect that there is a natural definition of 
even positive semidefinite matrix valued measures also in the 
$d$-dimensional case, $d\ge2$, and this takes the role of the 
spectral measure in the vector valued case. To define this 
notion first I prove a lemma. Before formulating it I recall 
the definition of a complex measure with finite 
total variation, since this notion appears in the formulation 
of the lemma. We say that a complex measure on
a measurable space has finite total variation if both its real
and imaginary part can be represented as the difference of 
two finite measures. I also recall Bochner's theorem, more 
precisely the version of this result that we shall apply in 
the proof.

\medskip\noindent
{\bf Bochner's theorem.} {\it Let $f(p)$, $p\in{\mathbb Z}^\nu$,
be a positive definite function on ${\mathbb Z}^\nu$, i.e., such 
a function for which the inequality
$\sum\limits_{j=1}^N\sum\limits_{j'=1}^N z_j\bar z_{j'}f(p_j-p_{j'})\ge0$
holds for any set of points $p_j\in{\mathbb Z}^\nu$, and complex 
numbers $z_j$, $1\le j\le N$, with some number $N\ge1$. Then 
there exists a unique finite measure $G$ on the torus 
$[-\pi,\pi)^\nu$ such that
$$
f(p)=\int_{[-\pi,\pi)^\nu}e^{i(p,x)}G(\,dx) \quad \textrm{for all }
p\in{\mathbb Z}^\nu.
$$
If the function $f$ is real valued, then the measure $G$ is
even, i.e., $G(-A)=G(A)$ for all $A\subset[-\pi,\pi)^\nu$.}

\medskip
Next I formulate the following lemma.

\medskip\noindent
{\bf Lemma 2.1.} {\it Let $X(p)=(X_1(p),\dots,X_d(p))$, 
$p\in{\mathbb Z}^\nu$, be a $d$-dimensional stationary 
Gaussian random field with expectation zero. Then for 
all pairs $1\le j,j'\le d$ the correlation function 
$r_{j,j'}(p)=EX_j(0)X_{j'}(p)$, $p\in{\mathbb Z}^\nu$, 
can be written in the form
\begin{equation}
r_{j,j'}(p)=EX_j(0)X_{j'}(p)=EX_j(m)X_{j'}(m+p)
=\int_{[-\pi,\pi)^\nu} e^{i(p,x)}G_{j,j'}(\,dx) \label{2.1}
\end{equation}
with a complex measure $G_{j,j'}$ on the torus 
$[-\pi,\pi)^\nu$ with finite total variation. The function 
$r_{j,j'}(p)$, $p\in{\mathbb Z}^\nu$, uniquely determines this 
complex measure $G_{j,j'}$ with finite total 
variation. It is even, i.e., $G_{j,j'}(-A)=\overline{G_{j,j'}(A)}$ 
for all measurable sets $A\subset[-\pi,\pi)^\nu$. The relation 
$G_{j',j}(A)=\overline{G_{j,j'}(A)}$ also holds for all 
$1\le j,j'\le d$ and $A\subset[-\pi,\pi)^\nu$.}

\medskip\noindent
{\it Remark.}\/ Let us remark that given a $d$-dimensional 
stationary random field with expectation zero, there exist
also such $d$-dimensional stationary random fields with
expectation zero which are Gaussian and have the same 
correlation function. As a consequence, in Lemma~2.1 we 
could drop the condition that the stationary random field
we are considering is Gaussian. The same can be told about
the other results of Section~2. I imposed this condition,
because later, as we work with random spectral measures and
random integrals with respect to them the Gaussian property
of the underlying random field is important.

\medskip\noindent
{\it Proof of Lemma 2.1.}
By Bochner's theorem we may write
$$
r_{j,j}(p)=\int_{[-\pi,\pi)^\nu} e^{i(p,x)}G_{j,j}(\,dx), \quad 
p\in{\mathbb Z}^\nu, 
$$ 
for all $1\le j\le d$ with some finite measure $G_{j,j}$ on 
$[-\pi,\pi)^\nu$. We find a good representation for $r_{j,j'}(n)$ 
if $j\neq j'$ with the help of following argument. 

The function 
\begin{eqnarray*}
q_{j,j'}(p)&=&E[X_j(0)+iX_{j'}(0)][X_j(p)-iX_{j'}(p)] \\
&=&E[X_j(0)+iX_{j'}(0)]\overline{[X_j(p)+iX_{j'}(p)]},
\end{eqnarray*}
$p\in{\mathbb Z}^\nu$, is positive definite, hence it can be 
written in the form 
$$
E[X_j(0)+iX_{j'}(0)][X_j(p)-iX_{j'}(p)]=\int_{[-\pi,\pi)^\nu}
e^{i(p,x)}H_{j,j'}(\,dx)
$$ 
with some finite measure $H_{j,j'}$ on $[-\pi,\pi)^\nu$.
Similarly,
$$
E[X_j(0)+X_{j'}(0)][X_j(p)+X_{j'}(p)]=\int_{-[\pi,\pi)^\nu} e^{i(p,x)}K_{j,j'}(\,dx)
$$ 
with some finite measure $K_{j,j'}$ on $[-\pi,\pi)^\nu$. Hence
\begin{eqnarray*}
EX_j(0)X_{j'}(p)&=&
\frac i2E[X_j(0)+iX_{j'}(0)][X_j(p)-iX_{j'}(p)] \\
&&\quad +\frac12E[X_j(0)+X_{j'}(0)][X_j(p)+X_{j'}(p)]\\
&&\quad\;-\frac{(1+i)}2 [EX_j(0)X_j(p)+EX_{j'}(0)X_{j'}(p)] \\
&=&\int_{[-\pi,\pi)^\nu} e^{i(p,x)}G_{j,j'}(\,dx)
\end{eqnarray*}
with $G_{j,j'}(\,dx)=\frac12[i H_{j,j'}(\,dx)+K_{j,j'}(\,dx)]
-\frac{(1+i)}2[G_{j,j}(\,dx)+G_{j',j'}(\,dx)]$.

In such a way we have found complex measures 
$G_{j,j'}$ with finite total variation which satisfy relation~(\ref{2.1}). 
Since this relation holds for all $p\in{\mathbb Z}^\nu$, the 
function $r_{j,j'}(p)$, $p\in{\mathbb Z}^\nu$, determines the 
measure $G_{j,j'}$ uniquely.

Since $r_{j,j'}(p)$ is real valued, i.e., 
$r_{j,j'}(p)=\overline{r_{j,j'}(p)}$, it can be written both in 
the form 
$$
r_{j,j'}(p)=\int_{[-\pi,\pi)^\nu} e^{i(p,x)}G_{j,j'}(\,dx)
$$
and 
$$
r_{j,j'}(p)=\int_{[-\pi,\pi)^\nu} e^{-i(p,x)}\overline{G_{j,j'}(\,dx)}=
\int_{[-\pi,\pi)^\nu} e^{i(p,x)}\overline{G_{j,j'}(-\,dx)}.
$$ 
Comparing 
these relations we get that $G_{j,j'}(A)=\overline{G_{j,j'}(-A)}$ 
for all measurable sets $A\subset[-\pi,\pi)^\nu$. Similarly, 
the relation $r_{j',j}(p)=r_{j,j'}(-p)$ implies that 
$G_{j',j}(A)=G_{j,j'}(-A)=\overline{G_{j,j'}(A)}$ for all measurable 
sets $A\subset[-\pi,\pi)^\nu$. Lemma~2.1 is proved.  

\medskip
Since all complex measures $G_{j,j'}$, $1\le j,j'\le d$, have 
finite total variation by Lemma~2.1, there is a finite measure 
$\mu$ on the torus $[-\pi,\pi)^\nu$ such that all these 
complex measures $G_{j,j'}$ are absolutely continuous with 
respect to $\mu$, and the absolute value of the Radon--Nikodym 
derivatives $g_{j,j'}(x)=\frac{dG_{j,j'}}{d\mu}(x)$ is integrable 
with respect to $\mu$. The properties of the measures $G_{j,j'}$ 
proved in Lemma~2.1 imply that the $d\times d$ matrix 
$(g_{j,j'}(x))$, $1\le j,j'\le d$, is Hermitian for almost all 
$x\in [-\pi,\pi)^\nu$ with respect to the measure $\mu$. We 
shall call the matrix valued measure 
$(G_{j,j'}(A))$, $A\subset[-\pi,\pi)^\nu$, positive semidefinite 
if the matrix $(g_{j,j'}(x))$, $1\le j,j'\le d$, is positive 
semidefinite for almost all $x\in [-\pi,\pi)^\nu$ with respect
to $\mu$. More precisely, we introduce the following definition.

\medskip\noindent
{\bf Definition of positive semidefinite matrix valued, even 
measures on the torus.} {\it Let us have some complex  
measures $G_{j,j'}$, $1\le j,j'\le d$, with finite total 
variation on the $\sigma$-algebra of the Borel measurable sets of 
the torus $[-\pi,\pi)^\nu$. Let us consider the matrix valued measure 
$(G_{j,j'})$, $1\le j,j'\le d$. We call this matrix valued measure 
positive semidefinite if there exists a (finite) positive 
measure $\mu$ on $[-\pi,\pi)^\nu$  such that all complex  
measures $G_{j,j'}$, $1\le j,j'\le d$, are absolutely 
continuous with respect to it, and their Radon--Nikodym derivatives 
$g_{j,j'}(x)=\frac{dG_{j,j'}}{d\mu}(x)$, $1\le j,j'\le d$, constitute 
a positive semidefinite matrix $(g_{j,j'}(x))$, $1\le j,j'\le d$
for almost all $x\in{\mathbb Z}^\nu$ with respect to the 
measure~$\mu$. We call this positive semidefinite matrix valued 
measure $(G_{j,j'})$, $1\le j,j'\le d$, on the torus even if 
$G_{j,j'}(-A)=\overline{G_{j,j'}(A)}$ for all measurable sets 
$A\subset[-\pi,\pi)^\nu$ and $1\le j,j'\le d$. 

Later we shall speak also of  positive semidefinite matrix 
valued even measures on a torus $[-A,A)^\nu$ for arbitrary 
$A>0$ which is defined in the same way, only the complex 
measures $G_{j,j'}$ and the dominating measure 
$\mu$ are defined on $[-A,A)^\nu$.}

\medskip\noindent
{\it Remark.}\/ Here I am speaking about measures with finite
total variation, although such (complex) measures are called
generally bounded measures in the literature. 
Actually, we know by Stone's theorem that any bounded signed 
measure can be represented as the difference of two bounded 
measures (with disjoint support). Nevertheless, I shall remain 
at this name, because actually we prove directly the finite 
total variation of the measures we shall work with in this paper. 
Besides, (in Section~4) I shall
define complex measures on ${\mathbb R}^\nu$ with locally 
finite total variation, and I prefer such a name which refers to
the similarity of these objects. (The complex measures with locally
finite total variation are not measures in the original meaning
of this word, only their restrictions to compact sets are complex 
measures.)

\medskip
The next theorem about the characterization of the correlation 
function of a $d$-di\-men\-sio\-nal stationary Gaussian random 
field with zero expectation states that the correlation 
functions $r_{j,j'}(p)$, $1\le j,j'\le d$, $p\in{\mathbb Z}^\nu$, 
can be given in the form~(\ref{2.1}) with the help of a positive 
semidefinite matrix valued, even measure $(G_{j,j'})$, 
$1\le j,j'\le d$, on the torus $[-\pi,\pi)^\nu$. Moreover, it 
will be shown that we have somewhat more freedom when we choose
a dominating measure~$\mu$ in the definition of positive 
semidefinite matrix valued measures on the torus. If the 
coordinates of a matrix valued measure $(G_{j,j'})$, 
$1\le j,k\le d$, are complex measures with finite total 
variation, and this matrix valued measure satisfies the 
definition of the positive semidefinite property with some 
measure~$\mu$, then this measure $\mu$ can be replaced in 
the definition by any such finite measure on the torus with 
respect to which the complex measures $G_{j,j'}$ are absolutely 
continuous. More explicitly, the following result holds. 

\medskip\noindent
{\bf Theorem 2.2.} {\it The covariance matrices of a $d$-dimensional
stationary random field $X(p)=(X_1(p),\dots,X_d(p))$, 
$p\in{\mathbb Z}^\nu$,  with expectation zero can be given 
in the following form. For all $1\le j,j'\le d$ there exists 
a complex measure $G_{j,j'}$ with finite total 
variation on the $\nu$-dimensional torus $[-\pi,\pi)^{\nu}$ 
in such a way that for all  $1\le j,j'\le d$ the correlation 
function $r_{j,j'}(p)=EX_j(0)X_{j'}(p)$, $p\in{\mathbb Z}^\nu$, 
is given by formula~(\ref{2.1}) with this complex 
measure~$G_{j,j'}$. The $d\times d$ matrix $G=(G_{j,j'})$, 
$1\le j,j'\le d$, whose coordinates are the complex  
measures $G_{j,j'}$ has the following properties. This 
matrix is Hermitian, i.e., the measures $G_{j,j'}$ satisfy the 
relation $G_{j',j}(A)=\overline{G_{j,j'}(A)}$ for all pairs of indices 
$1\le j,j'\le d$ and measurable sets $A\subset [-\pi,\pi)^\nu$,
and the measures $G_{j,j'}$ are even, i.e.,
$G_{j,j'}(-A)=\overline{G_{j,j'}(A)}$ for all $1\le j,j'\le d$ and 
$A\subset[-\pi,\pi)^\nu$. For all pairs $(j,j')$, $1\le j,j'\le d$, 
the function $r_{j,j'}(p)$, $p\in{\mathbb Z}^\nu$, defined by 
formula~(\ref{2.1}) uniquely determines the complex  
measure $G_{j,j'}$ with finite total variation. Besides, $G_{j,j'}$
has the following property.

Let us take a finite measure $\mu$ on the torus $[-\pi,\pi)^\nu$ 
such that all complex measures $G_{j,j'}$ are absolutely 
continuous with respect to it (because of the finite total 
variation of the complex measures $G_{j,j'}$ there 
exist such measures), and put 
$g_{j,j'}(x)=g_{j,j',\mu}(x)=\frac{dG_{j,j'}}{d\mu}(x)$. Then the 
matrix $(g_{j,j'}(x))$, $1\le j,j'\le d$, is positive 
semidefinite for almost all $x\in[-\pi,\pi)^\nu$ with respect
to the measure~$\mu$.

Conversely, if a class of complex measures 
$G_{j,j'}$ on $[-\pi,\pi)^\nu$, $1\le j,j'\le d$, have finite 
total variation, and $(G_{j,j'})$, $1\le j,j'\le d$, is a positive 
semidefinite matrix valued, even measure on the torus, then 
there exists a $d$-dimensional stationary Gaussian random field 
$X(p)=(X_1(p),\dots,X_d(p))$, $p\in{\mathbb Z}^\nu$, with 
expectation $EX_j(p)=0$ and covariance 
$EX_j(p)X_{j'}(q)=r_{j,j'}(p-q)$, where the function $r_{j,j'}(p)$ 
is defined in~(\ref{2.1}) with the complex  
measure $G_{j,j'}$ for all parameters $1\le j,j'\le d$ and 
$p,q\in{\mathbb Z}^\nu$.}

\medskip\noindent
{\it Remark.} We shall call the positive semidefinite matrix
valued, even measure $(G_{j,j'})$, $1\le j,j'\le d$, on the torus 
$[-\pi,\pi)^\nu$ with coordinates $G_{j,j'}$ satisfying 
relation~(\ref{2.1}) the matrix valued spectral measure of the 
correlation function $r_{j,j'}(p)$, $1\le j,j'\le d$, 
$p\in{\mathbb Z}^\nu$. In general, we shall call an
arbitrary positive semidefinite matrix valued, even measure 
on the torus $[-\pi,\pi)^\nu$ a matrix valued spectral
measure on the torus $[-\pi,\pi)^\nu$. (More generally, later 
we shall call for any $A>0$ a positive semidefinite matrix 
valued, even measure on the torus $[-A,A)^\nu$ a matrix 
valued spectral measure on this torus.) We have the 
right for such a terminology, since by Theorem~2.2 for an 
arbitrary positive semidefinite matrix valued, even measure 
on the torus $[-\pi,\pi)^\nu$ there exists a vector valued 
stationary Gaussian random field on ${\mathbb Z}^\nu$ such 
that this positive semidefinite matrix valued, even measure 
is the spectral measure of its correlation function.

\medskip\noindent
{\it Proof of Theorem 2.2.} The statements formulated in the
first paragraph of Theorem~2.2 follow from Lemma~2.1.
Next we prove that the matrix $(g_{j,j'}(x))$, $1\le j,j'\le d$, 
whose elements are defined as the Radon--Nikodym
derivatives of the complex measures $G_{j,j'}$ 
with respect to a measure $\mu$ satisfying the conditions 
of Theorem~2.2 is positive semidefinite for $\mu$~almost all~$x$.

We prove this by first showing with the help of Weierstrass' 
second approximation theorem that
\begin{equation}
\int_{[-\pi,\pi)^\nu} v(x) g(x) v^*(x)\mu(\,dx)\ge0 \label{2.2}
\end{equation}
for any continuous $d$-dimensional vector valued function
\hfill\break 
$v(x)=(v_1(x),\dots,v_d(x))$ on the $\nu$-dimensional torus 
$[-\pi,\pi)^\nu$, where $g(x)$ denotes the $d\times d$ matrix 
$(g_{j,j'}(x))$, $1\le j,j'\le d$, and $v^*(x)$ is the conjugate 
of the vector $v(x)$. 

To prove (\ref{2.2}) let us first observe that by Weierstrass' second 
approximation theorem for all $\varepsilon>0$ there exists a 
number $N=N(\varepsilon)$ and $d$ trigonometrical polynomials 
of order~$N$ 
$$
v_{N,j}(x)=\sum_{\substack{s=(s_1,\dots,s_\nu)\\ 
-N\le s_k< N,\;1\le k\le\nu}}  
a_{j,s_1,\dots,s_\nu}e^{i(s,x)}, \quad 1\le j\le d, \quad
x\in[-\pi,\pi)^\nu
$$
for which 
$$
\sup _{x\in[-\pi,\pi)^\nu}|v_{N,j}(x)-v_j(x)|\le\varepsilon 
\quad\textrm{for all }1\le j\le d.
$$
Let us also define the random vector $Y_N=(Y_{N,1},\dots,Y_{N,d})$ with
coordinates
$$
Y_{N,j}=\sum_{\substack{s=(s_1,\dots,s_\nu)\\ 
-N\le s_k< N, \;1\le k\le\nu}} 
a_{j,s_1,\dots,s_\nu}X_j(s), \quad 1\le j\le d,
$$
Then we have because of the relation 
$EX_j(s)X_{j'}(s')=\int e^{i(s-s',x)}g_{j,j'}(x)\mu(\,dx)$
$$
0\le E\left(\sum_{j=1}^d Y_{N,j}\right)
\overline{\left(\sum_{j=1}^d Y_{N,j}\right)}
=\sum_{j=1}^d\sum_{j'=1}^d \int_{[-\pi,\pi)^\nu}
g_{j,j'}(x)v_{N,j}(x)\overline{v_{N,j'}(x)}\mu(\,dx).
$$
Hence
$$
\int_{[-\pi,\pi)^\nu} v_N(x) g(x) v_N^*(x)\mu(\,dx)\ge0,
$$
and we get relation~(\ref{2.2}) from it with the  help of the limiting 
procedure $N\to\infty$.

Let us choose a vector $a=(a_1,\dots,a_d)\in {\mathbb R}^d$ and a non-negative
continuous function $u(x)$ on the torus $[-\pi,\pi)^\nu$. Let us apply
formula (\ref{2.2}) with the choice of the function 
$v(x)=(a_1\sqrt{u(x)},\dots,a_d\sqrt{u(x)})$. With this choice 
formula (\ref{2.2}) yields that 
$$
0\le\int_{[-\pi,\pi)^\nu} v(x) g(x) v^*(x)\mu(\,dx)
=\int_{[-\pi,\pi)} u(x) h_a(x)\mu(\,dx)
$$
with the function $h_a(x)=a g(x)a^*$. Since this inequality 
holds for all non-negative continuous functions this implies 
that $h_a(x)\ge0$ for almost all $x$ with respect to the 
measure~$\mu$. Moreover, since $h_a(x)=a g(x)a^*$ is a 
continuous function of the parameter~$a$ for a fixed number 
$x\in[-\pi,\pi)^\nu$ this also implies that $g(x)$ is a 
positive semidefinite matrix for almost all $x$ with respect
to the measure~$\mu$. We have proved that the covariance 
matrix of a vector valued stationary field has the properties 
stated in Theorem~2.2.

Next I show that if we have a class of complex 
measures $G_{j,j'}$ with finite total variation such that $(G_{j,j'})$ 
is a positive semidefinite matrix valued even measure on the 
torus, and the functions $r_{j,j'}(p)$, $p\in{\mathbb Z}^\nu$, 
are defined by formula~(\ref{2.1}) with these complex  
measures~$G_{j,j'}$, then there exists a vector valued stationary 
Gaussian field $X(p)=(X_1(p),\dots,X_d(p))$ with expectation 
zero and covariance function $EX_j(0)X_{j'}(p)=r_{j,j'}(p)$. 

First I show that for all $N\ge1$ there is a set of Gaussian 
random vectors $X(p)=(X_1(p),\dots,X_d(p))$, with parameters 
$p=(p_1,\dots,p_\nu)$, $-N\le p_j\le N$ for all $j=1,\dots,d$, such 
that $EX_j(p)X_{j'}(q)=r_{j,j'}(p-q))$ for all $1\le j,j'\le d$ and 
$p=(p_1,\dots,p_\nu)$, $q=(q_1,\dots,q_\nu)$ with 
$-N\le p_s,q_s\le N$, $1\le s\le\nu$. 

Let us observe that the covariances $r_{j,j'}(p)$ defined by~(\ref{2.1}) 
are real-valued, since $G_{j,j'}(A)=\overline{G_{j,j'}(-A)}$. To 
show that there exists a set of Gaussian random vectors with the 
desired covariance we have to check that the covariance matrix 
determined by the coordinates of these random vectors is positive 
semidefinite. This means that for all sets of complex numbers 
$$
{\cal A}_N=\{a_{j,p}=a_{j,p_1,\dots,p_\nu}\colon\;
1\le j\le d, \; -N\le p_s\le N, \textrm{ for all } 1\le s\le\nu\}
$$
$$
I({\cal A}_N)=\sum_{j=1}^d\sum_{j'=1}^d \sum_{\substack{p=(p_1,\dots,p_\nu)\\
-N\le p_s\le N,\; 1\le s\le \nu }}
\sum_{\substack{q=(q_1,\dots,q_\nu)\\
-N\le q_s\le N,\; 1\le s\le \nu}}
a_{j,p}\overline{a_{j',q}} r_{j,j'}(p-q)\ge0.
$$
This inequality holds since
\begin{eqnarray*}
I({\cal A}_N)&=&\int \sum_{j=1}^d\sum_{j'=1}^d 
\left(
\sum_{\substack{p=(p_1,\dots,p_\nu)\\
-N\le p_s\le N,\; 1\le s\le \nu}} 
a_{j,p} e^{i(p,x)}\right) g_{j,j'}(x) \\
&&\qquad\qquad  
\overline{\left(
\sum_{\substack{p=(p_1,\dots,p_\nu)\\
-N\le p_s\le N,\; 1\le s\le \nu}} 
a_{j',p} e^{i(p,x)}\right)}\mu(\,dx)\\
&&\qquad\qquad\qquad =\int\left( \sum_{j=1}^d\sum_{j'=1}^d 
b_j(x)g_{j,j'}(x)\overline{b_{j'}(x)}\right)\mu(\,dx)\ge0,
\end{eqnarray*}
where $b_j(x)=\sum\limits_{\substack{p=(p_1,\dots,p_\nu)\\
-N\le p_s\le N,\; 1\le s\le \nu}} a_{j,p} e^{i(p,x)}$. This
expression is really non-negative, since the matrix $g_{j,j'}(x)$
is positive semidefinite for $\mu$-almost all $x$, and this implies
that the integrand at the right-hand side of this expression is 
non-negative for $\mu$-almost all~$x$.

Since the distribution of the above sets of Gaussian random vectors
are consistent for different parameters~$N$ it follows from 
Kolmogorov's existence theorem for random processes with consistent 
finite distributions that there exists a Gaussian random field $X(p)$, 
$p\in{\mathbb Z}^\nu$, with $EZ_p=0$, $EX_j(p)X_{j'}(q)=r_{j,j'}(p-q)$,
where $r_{j,j'}(p)$ is defined by formula~(\ref{2.1}) with our matrix 
valued spectral measure $G=(G_{j,j'})$, $1\le j,j'\le d$. In such 
a way we constructed a stationary Gaussian random field with the 
desired properties. Theorem~2.2 is proved.   

\medskip
In the next lemma I give a different characterization of 
positive semidefinite matrix valued, even measures on the torus 
$[-\pi,\pi)^\nu$. 

\medskip\noindent
{\bf Lemma 2.3.} {\it Let us have a class of complex 
measures $G_{j,j'}$, $1\le j,j'\le d$, with finite total variation 
on the torus $[-\pi,\pi)^\nu$.  Let us define with their help the 
following $\sigma$-additive matrix valued function on the  
measurable subsets of the torus $[-\pi,\pi)^\nu$. Define for all 
measurable sets $A\subset[-\pi,\pi)^\nu$ the $d\times d$ matrix 
$G(A)=(G_{j,j'}(A))$, $1\le j,j'\le d$. This matrix valued function
is a positive semidefinite matrix valued, even measure on the torus 
$[-\pi,\pi)^\nu$ if and only if the matrix $(G_{j,j'}(A))$, 
$1\le j,j'\le d$, is positive semidefinite, and 
$G_{j,j'}(-A)=\overline{G_{j,j'}(A)}$ for all measurable sets 
$A\subset[-\pi,\pi)^\nu$ and $1\le j,j'\le d$.}

\medskip\noindent
{\it Proof of Lemma 2.3.}\/ It is clear that if $(G_{j,j'})$ is
a positive semidefinite matrix valued, even measure, then the 
matrix $(G_{j,j'}(A))$ with 
$$
G_{j,j'}(A)=\int_A g_{j,j'}(x)\mu(\,dx), \quad 1\le j,j'\le d,
$$ 
is a positive semidefinite matrix, and 
$G_{j,j'}(-A)=\overline{G_{j,j'}(A)}$ for all measurable sets
$A\subset[-\pi,\pi)^\nu$ and $1\le j,j'\le d$.

On the other hand, it is not difficult to see that if the
above properties hold, then 
$\sum\limits_{j=1}^d\sum\limits_{j'=1}^d
\int v_j(x)\overline{v_{j'}(x)} G_{j,j'}(\,dx)\ge0$
 for all vectors $v(x)=(v_1(x),\dots,v_d(x))$, where $
v_j(\cdot)$, $1\le j\le d$, is a continuous function on the 
torus $[-\pi,\pi)^\nu$. If $\mu$ is a finite measure on 
$[-\pi,\pi)^\nu$ such that  all complex measures 
$G_{j,j'}$, $1\le j,j'\le d$, are absolutely continuous with 
respect to it with Radon--Nikodym derivative $g_{j,j'}(x)$, 
and we denote the matrix $(g_{j,j'}(x))$, $1\le j,j'\le d$, 
by $g(x)$, then the above inequality can be rewritten in the 
form $\int v(x) g(x)v^*(x)\mu(\,dx)\ge0$. In the proof of 
Theorem~2.2 we have seen that this implies that $g(x)$ is
a positive semidefinite matrix for $\mu$ almost all 
$x\in[-\pi,\pi)^\nu$. Lemma~2.3 is proved.

\medskip
Let me also remark that the proof of Lemma~2.3 also implies 
that if the definition of positive semidefinite matrix valued, 
even measures holds with some finite measure $\mu$ on the
torus with the property that each complex  measure $G_{j,j'}$, 
$1\le j,j'\le d$, is absolutely continuous with respect to 
it, then the conditions of this definition also hold  with 
any measure $\mu$ on the torus with the same properties.

Given a positive semidefinite matrix valued even measure 
$G=(G_{j,j'})$, $1\le j,j'\le d$, on the torus $[-\pi,\pi)^\nu$, 
there is a natural candidate for the choice of the measure $\mu$
on the torus $[-\pi,\pi)^\nu$ with respect to which all measures 
$G_{j,j'}$, $1\le j,j'\le d$, are absolute continuous. We shall prove
an estimate in formula~(\ref{3.2}) which implies that the measure 
$\mu=\sum_{j=1}^d G_{j,j}$, i.e., the trace of the matrix valued measure
$G$ has this property. Later this measure will be our choice for the
measure $\mu$.

Let me remark that the proof of Lemma~2.3 yields another characterization
of positive semidefinite matrix valued measures on the torus. I
present it, although I shall not use it later.

A matrix valued measure $G=(G_{j,j'})$, $1\le j,j'\le d$, on the torus
such that $G_{j,j'}(A)=\overline{G_{j',j}(A)}$ for all $1\le j,j'\le d$
and measurable sets $A\subset [-\pi,\pi)^\nu$ is positive semidefinite
if and only if
$$
\sum_{j=1}^d\sum_{j'=1}^d \int_{[-\pi,\pi)^\nu} u_j(x)\overline{u_{j'}(x)} 
G_{j,j'}(\,dx)\ge0
$$
for all vectors $u(x)=(u_1(x),\dots,u_d(x))$ whose coordinates are
continuous functions on the torus $[-\pi,\pi)^\nu$.

\section{Random spectral measures in the multi-di\-men\-sio\-nal case}

If $X(p)=(X_1(p),\dots,X_d(p))$, $p\in{\mathbb Z}^\nu$, is a 
$d$-dimensional stationary Gaussian  random field with 
expectation zero, then its distribution is determined by its 
correlation functions 
$r_{j,j'}(p)=EX_j(0)X_{j'}(p)$, $1\le j,j'\le d$,
$p\in{\mathbb Z}^\nu$. In Theorem~2.2 we described 
this correlation function as the Fourier transform of a 
matrix valued spectral measure $G=(G_{j,j'})$, $1\le j,j'\le d$. 
In the case of scalar valued stationary random fields this 
result has a continuation. A so-called random spectral 
measure $Z_G$ can be constructed, and the elements of the 
stationary random field can be represented as an appropriate 
random integral with respect to it. This result can be 
interpreted so that the elements of a scalar valued 
stationary random field can be represented as the 
Fourier transforms of a random spectral measure. We want 
to find the multi-dimensional version of this result.

The results about scalar valued stationary random fields also 
help in the study of vector valued stationary random 
fields. Indeed, since the $j$-th coordinates $X_j(p)$, 
of the random vectors $X(p)$, $p\in{\mathbb Z}^\nu$, 
define a scalar valued stationary random field we can 
apply for them the results known in the scalar valued 
case. This enables us to construct such a random spectral 
measure $Z_{G.j}$ for all $1\le j\le d$ for which the 
identity $X_j(p)=\int_{[-\pi,\pi)^\nu} e^{i(p,x)}Z_{G,j}(\,dx)$ 
holds for all $p\in{\mathbb Z}^\nu$. The distribution 
of the random spectral measure $Z_{G,j}$ depends on the
coordinate $G_{j,j}$ of the matrix valued spectral 
measure $G$, which is the spectral measure of the
stationary random field $X_j(p)$, $p\in{\mathbb Z}^\nu$.  
For a fixed number $1\le j\le d$ the  properties of 
the random spectral measure $Z_{G,j}$ and the definition 
of the random integral with respect to it is worked out 
in the literature. I shall refer to my lecture notes~\cite{9}, 
where I described this theory. 

Nevertheless, the results obtained in such a way are not 
sufficient for us. They describe the distribution of the  
random spectral measure $Z_{G,j}$ for each $1\le j\le d$, 
but we need some additional results about their joint 
distribution. To get them I recall the results in~\cite{9} 
which led to the construction of the random spectral 
measures $Z_{G,j}$, and then I extend them in order to get 
the results we need to describe their joint distribution. 

I explain how we define simultaneously all random 
spectral measures $Z_{G,j}$, $1\le j\le d$, by recalling 
the method of~\cite{9} with some necessary modifications in 
the notation to adapt this method to our case. 

We construct the random spectral measure $Z_{G,j}$ for 
all $1\le j\le d$ in the following way. First we introduce 
two Hilbert spaces ${\cal K}^c_{1,j}$ and ${\cal H}^c_{1,j}$, 
and define an appropriate norm-preserving invertible 
linear transformation $T_j$ from ${\cal K}^c_{1,j}$ to 
${\cal H}^c_{1,j}$. (Here, and in the subsequent discussion 
I apply the superscript $^c$ in the notation to emphasize 
that we are working in a complex, and not in a real 
Hilbert space.) The Hilbert space ${\cal K}^c_{1,j}$ 
consists of those complex  valued functions 
$u(x)$ on the torus $[-\pi,\pi)^\nu$ for which 
$\int_{[-\pi,\pi)^\nu}|u(x)|^2G_{j,j}(\,dx)<\infty$, and the 
norm is defined in this space by the formula 
$\|u\|_{0,j}^2=\int_{[-\pi,\pi)^\nu} |u(x)|^2G_{j,j}(\,dx)$.
The Hilbert space ${\cal H}^c_{1,j}$ is defined as the closure
of the linear space consisting of the linear combinations
$\sum c_{p_s}X_j(p_s)$ with some (complex valued) coefficients
$c_{p_s}$ and parameters $p_s\in{\mathbb Z}^\nu$ in the
Hilbert space ${\cal H}^c$. The Hilbert space ${\cal H}^c$
consists of the complex valued random variables 
with finite second moment, measurable with respect to 
the $\sigma$-algebra generated by the random variables 
$X_j(p)$, $1\le j\le d$, $p\in{\mathbb Z}^\nu$, and 
the norm $\|\cdot\|_{1,j}$ in it is determined by 
the scalar product defined by the formula 
$\langle\xi,\eta\rangle=E\xi\bar\eta$, 
$\xi,\eta\in{\cal H}^c$. First we define the
transformation $T_j$ only for finite trigonometrical 
sums in ${\cal K}^c_{1,j}$. We define it by the formula 
$T_j(\sum c_{p_s}e^{i(p_s,x)})=\sum c_{p_s}X_j(p_s)$. We showed 
in~\cite{9} that we have defined in such a way a norm-preserving
linear transformation from an everywhere dense subspace
of ${\cal K}^c_{1,j}$ to an everywhere dense subspace of 
${\cal H}^c_{1,j}$. This can be extended to a norm-preserving 
invertible linear transformation $T_j$ from ${\cal K}^c_{1,j}$ to 
${\cal H}^c_{1,j}$ in a unique way. We define the random spectral 
measure $Z_{G,j}(A)$ for a measurable set $A\subset[-\pi,\pi)^\nu$ 
by the formula 
$Z_{G,j}(A)=T_j({\mathbb I}_A(\cdot))$, where ${\mathbb I}_A(\cdot)$ 
denotes the indicator function of the set~$A$.

It follows from the results of~\cite{9} that for any 
$1\le j\le d$ the measure $G_{j,j}$ determines the 
distribution of the random spectral measure $Z_{G,j}$, 
(i.e., the joint distribution of the random variables 
$Z_{G,j}(A_1),\dots Z_{G,j}(A_N)$ for all $N\ge1$ and
measurable sets $A_k\subset [-\pi,\pi)^\nu$, 
$1\le k\le N$). Next we shall study the joint 
distribution of the random fields $Z_{G,j}$ for all
$1\le j\le d$, i.e., the joint distribution of the
random variables $Z_{G,j}(A_1),\dots Z_{G,j}(A_N)$ for all 
$N\ge1$, measurable sets $A_k\subset [-\pi,\pi)^\nu$, 
$1\le k\le N$ and $1\le j\le d$. In particular, we shall
show that the joint distribution of the random fields
$Z_{G,j}$, $1\le j\le d$, are determined by the matrix 
valued spectral measure $G=(G_{j,j'})$, $1\le j,j'\le d$.
The joint distribution of these random fields are 
determined by the matrix valued measure $G$, and not 
only by their diagonal elements $G_{j,j}$, $1\le j\le d$.

To investigate the joint behaviour of the random spectral
measures $Z_{G,j}$, $1\le j\le d$, first we define two
Hilbert spaces ${\cal K}^c_1$ and ${\cal H}^c_1$ together with 
a norm-preserving  and invertible transformation between 
them. The elements of the Hilbert space ${\cal K}^c_1$ are 
the vectors $u=(u_1(x),\dots,u_d(x))$ with 
$u_j(x)\in{\cal K}^c_{1,j}$, $1\le j\le d$. To define the 
(semi)-norm in ${\cal K}^c_1$ we introduce a positive 
semidefinite bilinear form $\langle \cdot,\cdot\rangle_0$ 
on it. To make some subsequent discussions simpler I 
make the following convention in the rest of the paper. 
Given a positive semidefinite matrix valued measure 
$(G_{j,j'})$, $1\le j,j'\le d$, on the torus 
$[-\pi,\pi)^\nu$, I fix a finite and even measure 
$\mu$ on $[-\pi,\pi)^\nu$ such that all complex 
measures $G_{j,j'}$ are absolutely 
continuous with respect to it, and I denote by 
$g_{j,j'}(x)$ their Radon--Nikodym derivative with 
respect to $\mu$. With the help of this notation 
we define $\langle \cdot,\cdot\rangle_0$ in the 
following way. If 
$u(x)=(u_1(x),\dots,u_d(x))\in{\cal K}^c_1$  and 
$v(x)=(v_1(x),\dots,v_d(x))\in{\cal K}^c_1$, then
\begin{eqnarray}
\langle u(x), v(x)\rangle_0
&=&\sum_{j=1}^d\sum_{j'=1}^d\int u_j(x)\overline{v_{j'}(x)}G_{j,j'}(\,dx) 
\label{3.1} \\ 
&=&\sum_{j=1}^d\sum_{j'=1}^d\int g_{j,j'}(x)u_j(x)
\overline{ v_{j'}(x)}\mu(\,dx) \nonumber \\ 
&=&\int_{[-\pi,\pi)^\nu} u(x) g(x) v(x)^*\mu(\,dx) \nonumber
\end{eqnarray}
with the matrix $g(x)=(g_{j,j'}(x))$, $1\le j,j'\le d$, where 
$v^*(x)$ denotes the column vector whose elements are
the functions $\overline {v_k(x)}$, $1\le k\le d$. 

To show that the integral in the definition of
$\langle u(x), v(x)\rangle_0$ is convergent let us observe 
that 
\begin{equation}
|g_{j,j'}(x)|^2\le g_{j,j}(x)g_{j',j'}(x) \textrm{ for almost all } x 
\textrm{ with respect to the measure }\mu \label{3.2}
\end{equation}
for all $1\le j,j'\le d$, because $g(x)$ 
is a positive semidefinite matrix for almost all~$x$. This fact 
together with the Schwarz inequality imply that
\begin{eqnarray*}
&&   \!\!\!\!\!\!\!
\left|\int_{[-\pi,\pi)^\nu} u_j(x)g_{j,j'}(x)\overline{v_{j'}(x)}
\mu(\,dx)\right| \\ 
&&\!\le\int_{[-\pi,\pi)^\nu} |u_j(x)|
\sqrt{g_{j,j}(x)g_{j',j'}(x)}|v_{j'}(x)|\mu(\,dx) \\
&&\le\! \left(\int_{[-\pi,\pi)^\nu} |u_j(x)|^2g_{j,j}(x)\mu(\,dx)\right)^{1/2}
\left(\int_{[-\pi,\pi)^\nu} |v_{j'}(x)|^2g_{j',j'}(x)\mu(\,dx)\right)^{1/2} \\
&&\!<\infty 
\end{eqnarray*}
for all pairs $1\le j,j'\le d$ and $u_j\in{\cal K}^c_{1,j}$
and $v_{j'}\in{\cal K}^c_{1,k}$. This implies that the integral
in~(\ref{3.1}) is finite. Moreover, the last inequality implies that
\begin{eqnarray}
\langle u(x),u(x)\rangle_0&\le&
\left(\sum_{j=1}^d
\left(\int_{[-\pi,\pi)^\nu} |u_j(x)|^2G_{j,j}(\,dx)\right)^{1/2} \right)^2
\nonumber \\
&\le& d\sum_{j=1}^d \int_{[-\pi,\pi)^\nu}|u_j(x)|^2G_{j,j}(\,dx)
=d\sum_{j=1}^d \|u_j\|^2_{0,j} \label{3.3}
\end{eqnarray}
for all $u(x)=(u_1(x),\dots,u_d(x))\in{\cal K}^c_1$. 

Observe that $\langle u(x),u(x)\rangle_0\ge 0$, because $g(x)$ 
is a positive semidefinite matrix, which implies that 
$u(x) g(x)u^*(x)\ge0$ for almost all $x$ with respect to 
the measure $\mu$. In such a way we can define
the norm $\|\cdot\|_0$ in ${\cal K}^c_1$ by the formula
$\|u\|_0=\langle u(x),u(x)\rangle_0$. We identify two elements
$u$ and $v$ in ${\cal K}^c_1$ if $\|u-v\|_0=0$.

Next we define the Hilbert space ${\cal H}^c_1$ with the norm
$\|\cdot\|_1$ on it. The elements of ${\cal H}^c_1$ are the
vectors $\xi=(\xi_1,\dots,\xi_d)$, where $\xi_j\in{\cal H}^c_{1,j}$,
$1\le j\le d$, and we define the norm on it by the formula
$\|\xi\|_1^2=E\left|\sum_{j=1}^d\xi_j\right|^2$ if 
$\xi=(\xi_1,\dots,\xi_d)\in{\cal H}^c_1$. It is the norm
induced by the scalar product 
$\langle\xi,\eta\rangle_1=E\left(\sum_{j=1}^d\xi_j\right)
\overline{\left(\sum_{j=1}^d\eta_j\right)}$ for
$\xi=(\xi_1,\dots,\xi_d)\in{\cal H}^c_1$ and
$\eta=(\eta_1,\dots,\eta_d)\in{\cal H}^c_1$. We identify
two elements $\xi\in{\cal H}^c_1$ and $\eta\in{\cal H}^c_1$
if $\|\xi-\eta\|_1=0$.

Observe that
\begin{eqnarray}
\|\xi\|_1^2 \! &=&   \! E\left(\sum_{j=1}^d\xi_j\right)
\overline{\left(\sum_{j'=1}^d\xi_{j'}\right)}
\le\sum_{j=1}^d\sum_{j'=1}^d (E|\xi_j|^2)^{1/2}(E|\xi_{j'}|^2)^{1/2}  
\label{3.4} \\
\! &=&\!   \left(\sum_{j=1}^d (E|\xi_j|^2)^{1/2}\right)
\left(\sum_{j'=1}^d (E|\xi|_{j'}^2)^{1/2}\right)
\le d\sum_{j=1}^d E|\xi|_j^2=d\sum_{j=1}^k\|\xi_j\|^2_{1,j} \nonumber
\end{eqnarray}
for a vector $\xi=(\xi_1,\dots,\xi_d)\in{\cal H}^c_1$ 

We define the operator $T$ mapping from ${\cal K}^c_1$ to ${\cal H}^c_1$
by the formula 
$$
Tu=T(u_1,\dots,u_d)=(T_1u_1,\dots,T_du_d)
$$  
for 
$u=(u_1,\dots,u_d)$, $u_j\in{\cal K}^c_{1,j}$, with the help of 
the already defined operators $T_j$, $1\le j\le d$. We show 
that $Tu=T(u_1,\dots,u_d)=(T_1u_1,\dots,T_du_d)$ for 
$u=(u_1,\dots,u_d)\in{\cal K}^c_1$ is a norm preserving and 
invertible transformation from ${\cal K}^c_1$ to ${\cal H}^c_1$. To
prove this let us first observe that because of inequality~(\ref{3.3})
and Weierstrass' second approximation theorem the finite linear
combinations
$$
\left(\sum_{p\in A_N}c_{1,p}e^{i(p,x)},\dots,
\sum_{p\in A_N} c_{d,p}e^{i(p,x)}\right),
$$
where $A_N=\{p=(p_1,\dots,p_\nu)\colon\; -N\le p_s\le N, 
\textrm{ for all }1\le s\le \nu\}$,
constitute an everywhere dense linear subspace in ${\cal K}^c_1$,
and because of the inequality~(\ref{3.4}) the finite linear combinations
\begin{eqnarray}
&&\left(\sum_{p\in A_N} c_{1,p}X_1(p),\dots,\sum_{p\in A_N} c_{d,p}X_d(p)\right)
\nonumber \\
&&\qquad =T\left(\sum_{p\in A_N}c_{1,p}e^{i(p,x)},\dots,
\sum_{p\in A_N}c_{d,p}e^{i(p,x)}\right)
\label{3.5}
\end{eqnarray}
constitute an everywhere dense linear subspace in ${\cal H}^c_1$
if $N=1,2,\dots$, and the coefficients $c_{j,p}$, $1\le j\le d$,
$p\in A_N$, are arbitrary complex numbers. Hence the following
calculation implies that $T$ is a norm preserving and invertible
transformation from ${\cal K}^c_1$ to ${\cal H}^c_1$. 

If
$$
u(x)=\left(\sum_{p\in A_N}c_{1,p}e^{i(p,x)},\dots,
\sum_{p\in A_N}c_{d,p}e^{i(p,x)}\right)
$$
and 
$$
v(x)=\left(\sum_{p\in A_N}c'_{1,p}e^{i(p,x)},\dots,
\sum_{p\in A_N}^Nc'_{d,p}e^{i(p,x)}\right),
$$
then
\begin{eqnarray*}
&&\langle u(x),v(x)\rangle_0
=\left\langle\left(\sum_{p\in A_N} c_{1,p}e^{i(p,x)},\dots,
\sum_{p\in A_N}c_{d,p}e^{i(p,x)}\right),\right.  \\
&&\qquad \qquad \qquad\qquad \left.\left(\sum_{p\in A_N}c'_{1,p}e^{-i(p,x)},\dots,
\sum_{p\in A_N} c'_{d,p}e^{-i(p,x)}\right)
\right\rangle_0 \\
&&\qquad =\sum_{j=1}^d\sum_{j'=1}^d \sum_{s\in A_N}\sum_{t\in A_N} 
c_{j,s}\bar {c'}_{j',t}
\int_{[-\pi,\pi)^\pi} g_{j,j'}(x) e^{i(s-t,x)}\mu(\,dx) \\
&&\qquad=E\left(\sum_{j=1}^d\sum_{s\in A_N} c_{j,s}X_j(s)\right)
\overline{\left(\sum_{j'=1}^d\sum_{t\in A_N} c'_{j',t}X_{j'}(t)\right)}
=\langle Tu(x), Tv(x)\rangle_1.
\end{eqnarray*}

We shall define the random variables $Z_{G,j}(A)$ for all
indices $1\le j\le d$ and measurable sets 
$A\subset [-\pi,\pi)^\nu$, by the formula 
$Z_{G,j}(A)=T_j({\mathbb I}_A(x))$ with the above defined 
operators $T_j$, $1\le j\le d$, where ${\mathbb I}_A(\cdot)$ 
denotes the indicator function of the set 
$A\subset[-\pi,\pi)^\nu$. Next I formulate some properties 
of this class of random variables. These properties will 
appear in the definition of random spectral measures. All 
sets appearing in the next statements are measurable 
subsets of the torus $[-\pi,\pi)^\nu$. 

\medskip
\begin{description}
\item[{\rm(i)}] The random variables  $Z_{G,j}(A)$ are complex 
valued, and their real and imaginary parts are jointly Gaussian, 
i.e., for any positive integer $N$ and sets $A_s$, $1\le s\le N$, 
the random variables $\textrm{Re}\, Z_{G,j}(A_s)$, 
$\textrm{Im}\, Z_{G,j}(A_s)$, $1\le s\le N$, 
$1\le j\le d$, are jointly Gaussian.
\item[{\rm(ii)}] $EZ_{G,j}(A)=0$ for all $1\le j\le d$ and $A$,
\item[{\rm(iii)}] $EZ_{G,j}(A)\overline {Z_{G,{j'}}(B)}=G_{j,j'}(A\cap B)$
for all $1\le j,j'\le d$ and sets $A,B$.
\item[{\rm(iv)}] $\sum\limits_{s=1}^nZ_{G,j}(A_s)
=Z_{G,j}\left(\bigcup\limits_{s=1}^n A_s\right)$ if
$A_1,\dots,A_n$ are disjoint sets, $1\le j\le d$.
\item[{\rm(v)}] $Z_{G,j}(A)=\overline{Z_{G,j}(-A)}$ for all 
$1\le j\le d$ and sets $A$.
\end{description}

\medskip
Properties (i)--(v) were proved in the one-dimensional case, 
for example, in~\cite{9}. The only difference in checking its several 
dimensional version is that we have to apply the 
multi-dimensional operator $T$ from ${\cal K}^c_1$ to 
${\cal H}^c_1$ to prove property~(i), and to apply the same 
mapping~$T$ in proving~Property ~(iii). Here we exploit 
that $\langle u,v\rangle_0=\langle Tu,Tv\rangle_1$. We apply 
this identity with the vector $u\in{\cal K}^c_1$ whose 
$j$-th coordinate is ${\mathbb I}_A(x)$, and the other coordinates 
are zero and  the vector $v\in{\cal K}^c_1$ whose $k$-th 
coordinate is ${\mathbb I}_B(x)$ and the other coordinates are zero. 
Property~(v) can be proved as the special case of the 
following more general relation.

\medskip
\begin{description}
\item[{\rm(v$'$)}] $T_j(u)=\overline{T_j(u_-)}$ for all $1\le j\le d$ 
and $u\in{\cal K}^c_j$, where $u_-(x)=\overline{u(-x)}$.
\end{description}
\medskip

Property~(v$'$) can be proved by first proving it in the special
case when $u(x)$ is a trigonometrical polynomial, and then 
applying a limiting procedure.

Next we define the vector valued random spectral measures 
corresponding to a matrix valued spectral measure. 

\medskip\noindent
{\bf Definition of vector valued random spectral measures on 
the torus.} {\it Let a matrix valued spectral measure 
$G=(G_{j,j'})$, $1\le j,j'\le d$, be given on the torus 
$[-\pi,\pi)^\nu$ together with a set of complex valued 
random variables indexed by pairs $(j,A)$, where $1\le j\le d$, 
and $A$ is an element of the $\sigma$-algebra~${\cal A}$ 
$$
{\cal A}=\{ A\colon\; A\subset [-\pi,\pi)^\nu 
\textrm{ is a Borel measurable set}\}
$$ 
of the Borel measurable sets of the torus 
whose joint distribution depends on the matrix valued spectral 
measure~$G$. To recall this dependence we denote the random
variable indexed by a pair $(j,A)$, $1\le j\le d$, $A\in {\cal A}$,
by $Z_{G,j}(A)$. We call the set of random variables $Z_{G,j}(A)$,
$1\le j\le d$, $A\in{\cal A}$, a $d$-dimensional vector valued 
random spectral measure corresponding to the matrix valued 
spectral measure~$G$ on the torus $[-\pi,\pi)^\nu$ if this set 
of random variables satisfies properties~(i)--(v) defined above. 
Given a fixed parameter $1\le j\le d$ we call the set of random 
variables $Z_{G,j}(A)$, $A\in{\cal A}$, the $j$-th coordinate of 
this $d$-dimensional vector valued random spectral measure, and 
we denote it by $Z_{G,j}$. We denote the vector valued random 
spectral measure $Z_{G,j}(A)$, $1\le j\le d$, $A\in{\cal A}$, 
by $Z_G=(Z_{G,1},\dots,Z_{G,d})$. 

More generally, if a matrix valued spectral measure $G$ is given 
on the torus $[-B,B)^\nu$ with some number $B>0$ together with a
set of complex valued random variables $Z_{G,j}(A)$, where
$1\le j\le d$, and $A$ is a Borel measurable set on the torus
$[-B,B)^\nu$ which satisfies properties (i)--(v) defined above,
then we call this set of random variables a $d$-dimensional vector
valued random spectral measure corresponding to the spectral 
measure~$G$. We call the set of random variables 
$Z_{G,j}(A)$, $A\in{\cal A}$, for a fixed $1\le j\le d$ the $j$-th
coordinate of this vector valued spectral measure, and denote it
by~$Z_{G,j}$. We denote the vector valued spectral measure by
$Z_G=(Z_{G,1},\dots,Z_{G,d})$.}

\medskip\noindent
{\it Remark:}\/ If $G=(G_{j,j'})$, $1\le j,j'\le d$, is a 
matrix valued spectral measure, $Z_G=(Z_{G,1},\dots,Z_{G,d})$ 
is a vector valued spectral measure corresponding to it, 
then $G_{j,j}$ is a scalar valued spectral measure for any 
$1\le j\le d$, and $Z_{G,j}$ is a scalar valued random 
spectral measure corresponding to it. As we shall see in 
Lemma~3.3 the spectral measure $G$ determines the distribution
of the random spectral measure~$Z_G$.

\medskip
It follows from the above considerations that for any 
$d$-dimensional matrix valued spectral measure there exists 
a $d$-dimensional vector valued random spectral measure 
corresponding to it. We can define the random integral with 
respect to it by means of the method applied in the scalar 
valued case.

We shall define the random integrals of the functions
$f\in{\cal K}^c_{1,j}$ with respect to the random spectral
measure $Z_{G,j}$, $1\le j\le d$. First  we define these
integrals for elementary functions. They are finite sums of 
the form $\sum_{s=1}^N c_s {\mathbb I}_{A_s}(x)$, where 
$A_1,\dots,A_N$ are disjoint sets in $[-\pi,\pi)^\nu$, and 
$c_s$, $1\le s\le N$, are arbitrary complex numbers. Their 
integrals with respect to the random spectral measure 
$Z_{G,j}$, $1\le j\le d$, are defined as
$$
\int\left(\sum_{s=1}^N c_s {\mathbb I}_{A_s}(x)\right)Z_{G,j}(\,dx)
=\sum_{s=1}^N c_s Z_{G,j}(A_s).
$$
As it is remarked in~\cite{9}, property~(iv) implies that this definition
is meaningful, the integral of an elementary function does not depend
on its representation.  Then a simple calculation with the help of~(iii) 
shows that for two elementary functions $u$ and $v$
\begin{equation}
E\left(\int u(x)Z_{G,j}(\,dx)\overline{\int v(x)Z_{G,j}(\,dx)}\right)
=\int u(x)\overline{v(x)}G_{j,j}(\,dx), \quad  1\le j\le d. \label{3.6}
\end{equation}
This implies that the integral of the elementary functions with 
respect to the random spectral measure $Z_{G,j}$ define a norm 
preserving transformation from an everywhere dense subspace of 
the Hilbert space of ${\cal K}^c_{1,j}$ to an everywhere dense 
subspace of the Hilbert space of ${\cal H}^c_{1,j}$. This can be 
extended to a unitary transformation from ${\cal K}^c_{1,j}$ to 
${\cal H}^c_{1,j}$ in a unique way, and this extension defines 
the integral of a function $u\in{\cal K}^c_{1,j}$. It is clear 
that relation~(\ref{3.6}) remains valid for general functions 
$u,v\in{\cal K}^c_{1,j}$. Moreover, it is not difficult to see 
with the help of~(iii) that it can be generalized to the formula
\begin{equation}
E\left(\int u(x)Z_{G,j}(\,dx)\overline{\int v(x)Z_{G,{j'}}(\,dx)}\right)
=\int u(x)\overline{v(x)}G_{j,j'}(\,dx) \label{3.7}
\end{equation}
if $u\in{\cal K}^c_{1,j}$ and $v\in{\cal K}^c_{1,j'}$, $1\le j,j'\le d$.

It is clear that
\begin{equation}
E\int u(x) Z_{G,j}(\,dx)=0 \quad\textrm{for all } u\in{\cal K}_{1,j},
\quad 1\le j\le d.  \label{3.8}
\end{equation}

Another important property of the random integrals with respect to
$Z_{G,j}$ is that for all $1\le j\le d$
\begin{equation}
\int u(x)Z_{G,j}(\,dx) \quad\textrm {is real valued if } 
u(-x)=\overline{u(x)} \textrm{ for $\mu$ almost all }x\in[-\pi,\pi)^\nu. 
\label{3.9}
\end{equation}
This relation holds, since
$\int u(x)Z_{G,j}(\,dx) =\overline{\int u(x)Z_{G,j}(\,dx)}$
if $u(-x)=\overline{u(x)}$. We get this identity by means of the
change of variables $x\to -x$ with the help of relation~(v).

In the next theorem, I formulate the results we have about random
spectral measures and random integrals with respect to them.

\medskip\noindent
{\bf Theorem 3.1.} {\it Given a positive semidefinite matrix 
valued, even measure $G=(G_{j,j'})$, $1\le j,j'\le d$, on 
the torus $[-\pi,\pi)^\nu$ there exists a vector valued 
random spectral measure $Z_G=(Z_{G,1},\dots,Z_{G,d})$ 
corresponding to it. We have defined the random integrals 
$\int u(x)Z_{G,j}(\,dx)$ for all $1\le j\le d$ and 
$u\in{\cal K}^c_{1,j}$. This is a linear operator which 
satisfies relations~(\ref{3.7}), (\ref{3.8}), (\ref{3.9}), 
and the formula
\begin{equation}
X_j(p)=\int_{[-\pi,\pi)^\nu} e^{i(p,x)}Z_{G,j}(\,dx),  
\quad 1\le j\le d, \;\; p\in{\mathbb Z}^\nu, \label{3.10}
\end{equation}
defines a $d$-dimensional vector valued Gaussian stationary 
field whose matrix valued spectral measure is $G=(G_{j,j'})$, 
$1\le j,j'\le d$. Moreover, if a $d$-dimensional vector 
valued Gaussian stationary random field is given with this 
matrix valued spectral measure, then the random integrals 
in formula~(\ref{3.10}) taken with respect to the random spectral 
measure that we have constructed with its help through 
an operator $T$ in this section equals this vector valued 
Gaussian stationary random field.}

\medskip\noindent
{\it Proof of Theorem 3.1.}\/ We have already proved the 
existence  of the vector valued random spectral measure, 
and we constructed the random integral with respect to it. 
It satisfies formulas~(\ref{3.7}) and~(\ref{3.8}).  The random 
variables $X_j(p)$ defined in~(\ref{3.10}) are real valued 
by~(\ref{3.9}) and Gaussian with expectation zero. Hence we can 
show that they define a Gaussian stationary sequence with spectral
measure $G=(G_{j.j'})$, $1\le j,j'\le d$, by calculating their 
correlation function. We get by formula~(\ref{3.7}) that 
$EX_j(p)X_{j'}(q)=\int_{[-\pi,\pi)^\nu} e^{i(p-q,x)}G_{j,j'}(\,dx)$,
and this had to be checked. If the random spectral measure is 
constructed in the way as we have done in this section, then a 
comparison of the random integral we have defined with its help 
and of the operator~$T$ shows that 
$\int u(x)Z_{G,j}(\,dx)=T_j(u(x))$ for all $u\in{\cal K}^c_{1,j}$.
In particular,
$\int_{-[\pi,\pi)^\nu} e^{i(p,x)}Z_{G,j}(\,dx)=T_j(e^{i(p,x)})=X_j(p)$. 
This identity implies the last statement of Theorem~3.1. 
Theorem~3.1 is proved.

\medskip
Formula (\ref{3.9}) and Theorem~3.1 make possible to define for all
$1\le j\le d$ a real Hilbert space ${\cal K}_{1,j}$ consisting
of appropriate elements of ${\cal K}^c_{1,j}$ for which the operator
$T_j$ is a norm preserving invertible transformation from
${\cal K}_{1,}j$ to the real Hilbert space ${\cal H}_{1,j}$ consisting
of the real valued functions of the Hilbert space ${\cal H}^c_{1,j}$.
More precisely, the following statement holds.

\medskip\noindent
{\bf Lemma 3.2.} {\it Let $(G_{j,j'})$, $1\le j,j'\le d$, be a
matrix valued spectral measure on the torus $[-\pi,\pi)^\nu$,
and let $(Z_{G,1},\dots,Z_{G,d})$ be a vector valued spectral
measure corresponding to it. Define the $d$-dimensional vector
valued Gaussian stationary field $(X_1(p),\dots,X_p(d))$ by
formula~(\ref{3.10}) with the help of this vector valued random 
spectral measure. Define for all $1\le j\le d$ the  set of 
complex valued functions ${\cal K}_{1,j}$ on the torus 
$[-\pi,\pi)^\nu$ as
$$
{\cal K}_{1,j}=\left\{u\colon \int |u(x)|^2G_{j,j}(\,dx)<\infty, \quad 
u(-x)=\overline{u(x)}\textrm{ for all } x\in[-\pi,\pi)^\nu\right\}.
$$
Then ${\cal K}_{1,j}$ is a real Hilbert space with the scalar product
$$
\langle u,v\rangle=\int u(x)\overline{v(x)}G_{j,j}(\,dx), 
\quad u,v\in{\cal K}_{1,j}.
$$ 
Let ${\cal H}_{1,j}$ be the real Hilbert 
space consisting of the closure of the finite linear combinations 
$\sum_{k=1}^Nc_kX_j(p_k)$, $p_k\in{\mathbb Z}^\nu$, with real
coefficients~$c_k$ in the Hilbert space ${\cal H}$ of  random 
variables with finite second moments in the probability space 
where the random spectral measures $Z_{G,j}$ exists. (We define 
the scalar product in ${\cal H}$ in the usual way.) Then the 
map $T_j(u)=\int u(x)Z_{G,j}(\,dx)$, $u\in{\cal K}_{1,j}$, is a 
norm preserving, invertible linear transformation from 
the real Hilbert space ${\cal K}_{1,j}$ to the real 
Hilbert space~${\cal H}_{1,j}$.}

\medskip\noindent
{\it Proof of Lemma 3.2.}\/ The space ${\cal K}_{1,j}$ is a 
real Hilbert space, since the change of variable $x\to -x$ 
in the integral
$\langle u,v\rangle=\int u(x)\overline{v(x)}G_{j,j}(\,dx)$ 
implies that
$\langle u,v\rangle=\overline{\langle u,v\rangle}$ for all 
$u,v\in{\cal K}_{1,j}$ because of the evenness of the measure 
$G_{j,j}$. Clearly $e^{i(p,x)}\in{\cal K}_{1,j}$ for all 
$p\in{\mathbb Z}^\nu$. The class of functions 
${\cal K}_{1,j}$ agrees with the class of functions which 
have the form $u(x)=\frac{v(x)+\overline{v(-x)}}2$
with some $v\in{\cal K}^c_{1,j}$. As a consequence the set 
of finite trigonometrical polynomials $\sum c_ke^{i(p_k,x)}$, 
$p_k\in{\mathbb Z}^\nu$, with real valued coefficients $c_k$ 
is an everywhere dense subspace of ${\cal K}_{1,j}$. Since 
$T_j(\sum c_ke^{i(p_k,x)})=\sum c_k X_j(p_k)$, the transformation 
$T_j$ maps an everywhere dense subspace of ${\cal K}_{1,j}$ to
an everywhere dense subspace of ${\cal H}_{1,j}$. Because of
formulas (\ref{3.7}) and~(\ref{3.9}) $T_j$ is a norm preserving 
transformation in ${\cal K}_{1,j}$. Hence $T_j$ is an invertible, 
norm preserving transformation from ${\cal K}_{1,j}$ to 
${\cal H}_{1,j}$. Lemma~3.2 is proved.

\medskip
I would remark that the transformation $T_j$ on ${\cal K}_{1,j}$ 
defined in Lemma~3.2 is the restriction of the previously 
defined transformation $T_j$ on ${\cal K}^c_{1,j}$ to its subset 
${\cal K}_{1,j}$. I make also the following remark.

\medskip\noindent
{\bf Lemma 3.3.} {\it The positive semidefinite matrix valued,
even  measure $G(A)=(G_{j,j'}(A))$, $1\le j,j'\le d$, 
$A\in[-\pi,\pi)^\nu$, determines the distribution of a 
vector valued spectral random measure $Z_{G,j}$, $1\le j\le d$, 
corresponding to it.}

\medskip\noindent
To prove this lemma we have to show that for any collection of 
measurable sets $A_1$,\dots, $A_N$, the matrix valued measure 
$G(A)$ determines the joint distribution of the random vector 
consisting of the elements 
$\textrm{Re}\,Z_{G,j}(A_s)$, $\textrm{Im}\,Z_{G,j}(A_s)$, $1\le s\le N$,
$1\le j\le d$. Since this is a Gaussian random vector with 
expectation zero, it is enough to check that the covariance 
of these random variables can be expressed by means of the 
matrix valued measure $G(A)$. Since 
$\textrm{Re}Z_{G,j}(A)=\frac{Z_{G,j}(A)+\overline{Z_{G,j}(A)} }2$ 
and $\textrm{Im }Z_{G,j}(A)=\frac{Z_{G,j}(A)-\overline{Z_{G,j}(A)} }{2i}$ 
we can calculate these covariances with the help of 
properties~(iii) and~(v) of vector valued random spectral measures.

\medskip
Finally I prove an additional property of the vector valued 
random spectral measures which will be useful in Section~5,
in the study of multiple Wiener--It\^o integrals.

\medskip
\begin{description}
\item[{\rm(vi)}] The random variables of the form $Z_{G,j}(A\cup(-A))$ 
are real valued. Let a set $A\cup(-A)$ be disjoint from some sets
$B_1\cup(-B_1)$,\dots, $B_n\cup(-B_n)$. Then for any indices 
$1\le j,j'\le d$ the (complex valued) random vector 
$(Z_{G,j}(A), Z_{G,{j'}}(A))$, is independent of the random vector 
consisting of the elements $Z_{G,k}(B_s)$, $1\le s\le n$, $1\le k\le d$.
\end{description}

\medskip\noindent
{\it Proof of property (vi).}\/ It follows from property~(v) 
that $Z_{G,j}(A\cup(-A))=\overline{Z_{G,j}(A\cup(-A))}$, hence 
$Z_{G,j}(A\cup(-A))$ is real valued. To prove the second 
statement of~(vi) it is enough to check that under its conditions 
the (real valued) random variables $\textrm{Re}\,Z_{G,j}(A)$ 
and $\textrm{Im}\,Z_{G,j}(A)$ are uncorrelated to all random 
variables $\textrm{Re}\,Z_{G,k}(B_s)$, $\textrm{Im}\,Z_{G,k}(B_s)$,
$1\le s\le n$, $1\le k\le d$. This relation holds, since 
by the conditions of~(vi) $(\pm A)\cap(\pm B_s)=\emptyset$, 
hence relation~(iii) implies that
$EZ_{G,j}(\pm A)\overline{Z_{G,{j'}}(\pm B_s)}=0$ for all sets $B_s$, 
$1\le s\le n$, and indices $1\le j,j'\le d$. On the other hand, 
all covariances can be expressed as a linear combination of 
such expressions, since by relation~(v) $\textrm{Re}\,Z_{G,j}(\pm A)= 
\frac{Z_{G,j}(\pm A)+\overline{Z_{G,j}(\pm A)}}2
=\frac{Z_{G,j}(\pm A)+ Z_{G,j}(\mp A)}2$, and a similar relation 
holds also for $\textrm{Im}\,Z_{G,j}(\pm A)$,
$\textrm{Re}\,Z_{G,{j'}}(\pm B_s)$ and $\textrm{Im}\,Z_{G,{j'}}(\pm B_s)$,
$1\le s\le n$, $1\le j'\le d$. 

\section{Spectral representation of vector valued stationary 
generalized random fields}

In Sections 2 and 3 we discussed the properties of vector valued 
Gaussian stationary random fields with discrete parameters, which 
means a class of Gaussian random vectors $X(p)$, $p\in{\mathbb Z}^\nu$, 
with some nice properties. Similarly, we could have defined and 
investigated vector valued Gaussian stationary random fields 
with continuous parameters, where we consider a set of random 
vectors $X(t)$ indexed by  $t\in {\mathbb R}^\nu$ which have some nice 
properties. But we do not discuss this topic here. Here we 
define and investigate instead so-called vector valued 
Gaussian stationary generalized random fields  
$X(\varphi)=(X_1(\varphi),\dots,X_d(\varphi))$, parametrized 
with a nice linear space of functions~$\varphi$. 

Actually I am interested here in the vector valued
Gaussian stationary generalized random fields not for 
their own sake. We shall construct a class of vector 
valued Gaussian stationary generalized random fields. 
We shall show that their distribution can be described 
by means of a matrix valued spectral measure. We can 
also construct a vector valued random spectral measure 
in such a way that the elements of our vector valued 
generalized random field can be expressed in a form that 
can be considered as the Fourier transform of 
this random spectral measure. These matrix valued 
spectral measures and vector valued random spectral 
measures slightly differ from those defined in 
Sections~2 and~3, but since they are very similar to 
the corresponding objects defined for stationary 
random fields with discrete parameters it is natural 
to give them the same name.

The results that we shall prove are very similar to 
the results we got about vector valued random fields with 
discrete parameters. The main difference is that we can 
construct a larger class of matrix valued spectral 
measures and vector valued random spectral measures by 
means of generalized random fields. We shall need them, 
because in our later investigations we shall deal with 
such limit theorems where we can express the limit by 
means of these new, more general objects. On the other 
hand, these new vector valued random spectral measures 
behave similarly to the previous ones. In particular, 
the later results of this paper about multiple
Wiener--It\^o integrals also hold for this more general 
class of vector valued random spectral measures. Let 
me remark that we met a similar picture in the study 
of scalar valued Gaussian random fields in~\cite{9}, 
so that here we actually generalize the results in 
that work to the multi-dimensional case.

In the definition of vector valued generalized random 
fields we shall choose the functions of the Schwartz 
space for the class of parameter set. So to define the 
vector valued generalized random fields first I recall 
the definition of the Schwartz space, (see~\cite{6}).

We define the Schwartz space ${\cal S}$ of real valued 
functions on ${\mathbb R}^\nu$ together with its version 
${\cal S}^c$ consisting of complex valued functions on 
${\mathbb R}^\nu$. The space ${\cal S}^c=({\cal S}^\nu)^c$  
consists of those complex valued functions of $\nu$ 
arguments which decrease at infinity, together with 
their derivatives, faster than any polynomial. 
More explicitly, $\varphi\in{\cal S}^c$ for a complex 
valued function $\varphi$ defined on ${\mathbb R}^\nu$ if
$$
\left|x_1^{k_1}\cdots x_\nu^{k_\nu}\frac{\partial^{q_1+\cdots+q_\nu}}
{\partial x_1^{q_1}\dots \partial x_\nu^{q_\nu}}
\varphi(x_1,\dots,x_\nu)\right|
\le C(k_1,\dots,k_\nu,q_1,\dots,q_\nu)
$$
for all points $x=(x_1,\dots,x_\nu)\in {\mathbb R}^\nu$ and  vectors
$(k_1,\dots,k_\nu)$, $(q_1,\dots,q_\nu)$ with non-negative
integer coordinates and with some constant
$C(k_1,\dots,k_\nu,q_1,\dots,q_\nu)$ which may depend on the
function~$\varphi$. The elements of the space ${\cal S}$ 
are defined similarly, with the only difference that they 
are real valued functions.

To complete the definition of the spaces ${\cal S}$ 
and ${\cal S}^c$ we still have to define the topology
in them. We introduce the following topology in these
spaces. 

Let a basis of neighbourhoods of the origin consist of 
the sets
$$
U(k,p,\varepsilon)=\left\{\varphi\colon\; \varphi\in{\cal S}, \;\;
\max_{\substack{q=(q_1,\dots,q_\nu)\\
0\le q_s\le p,\textrm{ for all } 1\le s\le \nu}}
\sup_x(1+|x|^2)^k |D^q\varphi(x)|<\varepsilon\right\}
$$
with $k=0,1,2,\dots$, $p=1.2,\dots$ and $\varepsilon>0$, 
where $|x|^2=x_1^2+\cdots+x_\nu^2$, and 
$D^q=\frac{\partial^{q_1+\cdots+q_\nu}}
{\partial x_1^{q_1}\dots \partial x_\nu^{q_\nu}}$
for $q=(q_1,\dots,q_\nu)$. A basis of neighbourhoods of an 
arbitrary function $\varphi\in{\cal S}^c$ (or 
$\varphi\in{\cal S}$) consists of sets of the form 
$\varphi+U(k,q,\varepsilon)$, where the class of sets 
$U(k,q,\varepsilon)$ is a basis of neighbourhood of the 
origin. Actually we shall use only the following property 
of this topology. A sequence of functions 
$\varphi_n\in{\cal S}^c$ (or $\varphi_n\in{\cal S}$) 
converges to a function $\varphi$ in this topology if and 
only if
$$
\lim_{n\to\infty}\sup_{x\in {\mathbb R}^\nu}
(1+|x|^2)^k|D^q\varphi_n(x)-D^q\varphi(x)|=0
$$
for all $k=1,2,\dots$ and $q=(q_1,\dots,q_\nu)$.
The limit function $\varphi$ is also in the
space~${\cal S}^c$ (or in the space ${\cal S}$).

I shall define the notion of vector valued generalized 
random fields together with some related notions with 
the help of the notion of Schwartz spaces. A 
$d$-dimensional generalized random field is a random 
field whose elements are $d$-dimensional random vectors
$$
(X_1(\varphi),\dots,X_d(\varphi))=
(X_1(\varphi,\omega),\dots,X_d(\varphi,\omega))
$$ 
defined for all functions $\varphi\in {\cal S}$,
where ${\cal S}={\cal S}^\nu$ is the Schwartz space.
Before defining vector valued generalized
random fields I write down briefly the idea of 
their definition. This is explained in~\cite{9} 
and~\cite{10} in more detail.

Given a vector valued Gaussian stationary random 
field 
$$
X(t)=(X_1(t),\dots,X_d(t)), \quad t\in {\mathbb R}^\nu,
$$ 
we can define with its help the random field 
$X(\varphi)=(X_1(\varphi),\dots,X_d(\varphi))$,
$\varphi\in{\cal S}^\nu$, 
$X_j(\varphi)=\int\varphi(t)X_j(t)\,dt$, $1\le j\le d$, 
indexed by the elements of the Schwartz space, and 
this determines the original random field. We define 
generalized random fields with elements indexed by 
$\varphi\in{\cal S}$ as such random fields which
behave similarly to the random fields defined by
means of such integrals.

\medskip\noindent
{\bf Definition of vector valued generalized random fields.} 
{\it We say that the set of random vectors 
$(X_1(\varphi),\dots,X_d(\varphi))$,  $\varphi\in{\cal S}$, 
is a $d$-dimensional vector valued generalized random field 
over the Schwartz space ${\cal S}={\cal S}^\nu$ of rapidly
decreasing smooth functions if:

\medskip
\begin{description}
\item[{\rm(a)}]  $X_j(a_1\varphi+a_2\psi)=a_1X_j(\varphi)+
a_2X_j(\psi)$ with probability 1 for the $j$-th coordinate of 
the random vectors $(X_1(\varphi),\dots,X_d(\varphi))$ and  
$(X_1(\psi),\dots,X_d(\psi))$. This relation holds for each 
coordinate $1\le j\le d$, all real numbers $a_1$ and $a_2$,
 and pair of functions $\varphi$, $\psi$  from the Schwartz 
space ${\cal S}$. (The exceptional set of probability~0 
where this identity does not hold may depend on $a_1$, 
$a_2$, $\varphi$, and $\psi$.)

\item[{\rm(b)}] $X_j(\varphi_n)\Rightarrow X_j(\varphi)$ 
stochastically for any $1\le j\le d$ if 
$\varphi_n\to\varphi$ in the topology of ${\cal S}$.
\end{description}
}

\medskip
We also introduce the following definition. In its formulation 
we use the notation $\stackrel{\Delta}{=}$ for equality in 
distribution. 

\medskip\noindent
{\bf Definition of stationarity and Gaussian property for a vector 
valued generalized random field.} 
{\it The $d$-dimensional vector valued generalized random field
$X=\{(X_1(\varphi)\dots,X_d(\varphi)),\;\varphi\in {\cal S}\}$ 
is stationary if 
$$
(X_1(\varphi)\dots,X_d(\varphi))\stackrel{\Delta}{=}
(X_1(T_t\varphi)\dots,X_d(T_t\varphi))
$$
for all $\varphi\in{\cal S}$ and $t\in {\mathbb R}^\nu$, where 
$T_t\varphi(x)=\varphi(x-t)$. This field is called Gaussian
if $(X_1(\varphi),\dots,X_d(\varphi))$ is a Gaussian random 
vector for all $\varphi\in{\cal S}$. We call a vector valued
generalized random field a vector valued generalized random field
with zero expectation if $EX_j(\varphi)=0$ for all 
$\varphi\in{\cal S}$ and coordinates $1\le j\le d$.} 

\medskip
In the definition of stationarity and Gaussian property we 
imposed a condition for a single random vector. But because 
of the linearity property of generalized random fields 
formulated in property~(a) of their definition and the fact 
that if we have $N$ random vectors $\xi_1,\dots,\xi_N$ and 
$\eta_1,\dots,\eta_N$ such that the linear combinations 
$\sum\limits_{k=1}^Na_k\xi_k$ and $\sum\limits_{k=1}^Na_k\eta_k$ 
have the same distribution for any coefficients $a_k$, 
$1\le k\le N$, then the joint distribution of the random 
vectors $\xi_1,\dots,\xi_N$ and $\eta_1,\dots,\eta_N$ 
agree imply that an analogous statement holds about the 
properties of the joint distribution of several random 
vectors in a vector valued stationary random field. 
Indeed, if we take $N$ random vectors 
$(X_1(\varphi_k),\dots,X_d(\varphi_k))$, $1\le k\le N$, 
then their joint distribution agrees with the joint 
distribution of their shifts 
$(X_1(T_t\varphi_k),\dots,X_d(T_t\varphi_k))$, 
$1\le k\le N$, for any $t\in {\mathbb R}^\nu$. This follows 
from the fact that
$$
\sum_{k=1}^N a_k (X_1(\varphi_k),\dots,X_d(\varphi_k))
\stackrel{\Delta}{=} \sum_{k=1}^N a_k 
(X_1(T_t\varphi_k),\dots,X_d(T_t\varphi_k))
$$ 
for all $t\in {\mathbb R}^\nu$ and coefficients $a_k$, 
$1\le k\le N$, for a $d$-dimensional vector valued 
stationary generalized random field because of the 
linearity property of the generalized random fields 
and the properties of the operator $T_t$. A similar 
argument shows that the joint distribution of some 
vectors
$(X_1(\varphi_k),\dots,X_d(\varphi_k))$, $1\le k\le N$, 
in a
vector valued Gaussian generalized random field 
is Gaussian. 

\medskip
I shall construct a large class of $d$-dimensional vector
valued Gaussian stationary generalized random fields with 
expectation zero. I shall construct them with the help 
of positive semidefinite matrix valued  even measures on 
${\mathbb R}^\nu$. In the next step I write down this definition. 
The main difference between the definition of this notion 
and its counterpart defined on the torus $[-\pi,\pi)^\nu$ 
is that now we consider such complex  
measures which may have non-finite total variation. We 
impose instead a less restrictive condition. We shall 
work with complex measures on ${\mathbb R}^\nu$ which 
have locally finite total variation. For the sake of 
completeness I give their definition.

\medskip
\noindent
{\bf Definition of complex measures on ${\mathbb R}^\nu$
with locally finite total variation. The definition of their
evenness property.} {\it A complex measure on 
${\mathbb R}^\nu$ with locally finite total variation is such a complex 
valued function on the bounded, Borel measurable 
subsets of ${\mathbb R}^\nu$ whose restrictions to the measurable 
subsets of a cube $[-T,T]^\nu$ are complex  
measures with finite total variation for all $T>0$. We say
that a complex measure $G$ on ${\mathbb R}^\nu$ with 
locally finite total variation is even, if $G(-A)=\overline{G(A)}$ 
for all bounded and measurable sets $A\subset {\mathbb R}^\nu$.} 

\medskip
Let me remark that not all complex measures with 
locally finite total variation can be extended to a complex 
measure on all measurable subsets of ${\mathbb R}^\nu$. On the other 
hand, this can be done if we are working with a (real, 
positive number valued) measure. Next I formulate the 
definition we need in our discussion.

\medskip\noindent
{\bf Definition of positive semidefinite matrix valued 
measures on ${\mathbb R}^\nu$ with moderately increasing distribution 
at infinity. The definition of their evenness property.} 
{\it A Hermitian matrix valued measure on ${\mathbb R}^\nu$ is a 
class of such Hermitian matrices $(G_{j,j'}(A))$, 
$1\le j,j'\le d$, defined for all bounded, measurable sets 
$A\subset {\mathbb R}^{\nu}$ for which all coordinates $G_{j,j'}(\cdot)$,
$1\le j,j'\le d$, are complex measures on ${\mathbb R}^\nu$ 
with locally finite total variation. We call a Hermitian matrix 
valued measure $(G_{j,j'}(\cdot))$, $1\le j,j'\le d$, on ${\mathbb R}^\nu$ 
positive semidefinite if there exists a ($\sigma$-finite) 
positive measure $\mu$ on ${\mathbb R}^\nu$ such that for all numbers 
$T>0$ and indices  $1\le j,j'\le d$ the restriction of the 
complex measures $G_{j,j'}$ to the cube 
$[-T,T]^\nu$ is absolutely continuous with respect to $\mu$, 
and the matrices  $(g_{j,j'}(x))$, $1\le j,j'\le d$, defined 
with the help of the Radon--Nikodym  derivatives 
$g_{j,j'}(x)=\frac{dG_{j,j'}}{d\mu}(x)$, $1\le j,j'\le d$,
are Hermitian, positive semidefinite matrices for almost 
all $x\in {\mathbb R}^\nu$ with respect to the measure~$\mu$. We call
this Hermitian matrix valued measure $(G_{j,j'}(\cdot))$,
$1\le j,j'\le d$, on ${\mathbb R}^\nu$ even if the complex  
measures $G_{j,j'}$ with locally finite variation 
are even for all $1\le j,j'\le d$.

We shall say that the distribution of a positive semidefinite 
matrix valued measure $(G_{j,j'}(\cdot))$, $1\le j,j'\le d$, 
on ${\mathbb R}^\nu$ is moderately increasing at infinity if
\begin{equation}
\int (1+|x|)^{-r}G_{j,j}(\,dx)<\infty \quad\textrm{ for all }
1\le j\le d \textrm{ with some number } r>0.\label{4.1}
\end{equation}
}

\medskip\noindent
{\it Remark.} We can give, similarly to Lemma~2.3, a different 
characterization of positive semidefinite matrix valued, even 
measures on ${\mathbb R}^\nu$. Let us have some complex  
measures $G_{j,j'}$, $1\le j,j'\le d$, on the 
$\sigma$-algebra of the Borel measurable sets of ${\mathbb R}^\nu$ such 
that their restrictions to any cube $[-T,T]^\nu$, $T>0$, have 
finite total variation. Let us consider the matrix valued 
measure $(G_{j,j'}(A))$, $1\le j,j'\le d$ on ${\mathbb R}^\nu$ for all 
bounded, measurable sets $A\subset {\mathbb R}^\nu$. This matrix 
valued measure is positive semidefinite and even if and only
if it satisfies the following two conditions.

\medskip
\begin{description}
\item[{\rm(i.)}] The $d\times d$ matrix $(G_{j,j'}(A))$, 
$1\le j,j'\le d$, is Hermitian, positive semidefinite 
for all bounded, measurable sets $A\subset {\mathbb R}^\nu$.
\item[{\rm(ii.)}] $G_{j,j'}(-A)=\overline{G_{j,j'}(A)}$, for
all $1\le j,j'\le d$ and bounded, measurable sets 
$A\subset {\mathbb R}^\nu$.
\end{description}

\medskip
This statement has almost the same proof as Lemma~2.3.
The only difference in the proof is that now we have to
work with such vectors $v(x)=(v_1(x),\dots,v_d(x))$ whose 
coordinates $v_j(x)$ are continuous functions on ${\mathbb R}^\nu$ 
with bounded support, $1\le j\le d$. Let me also remark
that the following statement also follows from this proof.
If a matrix valued measure $(G_{j,j'}(A))$, $1\le j,j'\le d$,
on ${\mathbb R}^\nu$ satisfies the conditions in the definition of
positive semidefinite matrices with some $\sigma$-finite
measure $\mu$ on ${\mathbb R}^\nu$ with respect to which all 
complex measures $G_{j,j}$ are absolutely 
continuous, then it satisfies these conditions with any
 $\sigma$-finite measure $\mu$ on ${\mathbb R}^\nu$ with the 
same property. 

\medskip
Before constructing a large class of vector valued Gaussian 
stationary generalized random fields I recall an important 
property of the Fourier transform of the functions in the 
Schwartz spaces $\cal S$ and ${\cal S}^c$ (see for example~\cite{6}). 
Actually this property of the Schwartz spaces made useful 
their choice in the definition of generalized fields. 

The Fourier transform $f\to\tilde f$ is a bicontinuous map 
from ${\cal S}^c$ to~${\cal S}^c$. (This means that this 
transformation is invertible, and both the Fourier transform 
and its inverse are continuous maps from ${\cal S}^c$ to 
${\cal S}^c$.) (The restriction of the Fourier transform to 
the space ${\cal S}$ of real valued functions is a bicontinuous 
map from ${\cal S}$ to the subspace of ${\cal S}^c$ consisting 
of those functions $f\in{\cal S}^c$ for which
$f(-x)=\overline{f(x)}$ for all $x\in {\mathbb R}^\nu$.) 

Next I formulate the following result.

\medskip\noindent
{\bf Theorem 4.1 about the construction of vector valued Gaussian 
stationary generalized random fields with zero expectation.}
{\it Let $(G_{j,j'})$, $1\le j,j'\le d$, be a positive semidefinite 
matrix valued even measure on ${\mathbb R}^\nu$ whose distribution is 
moderately increasing at infinity.

Then there exists a vector valued Gaussian stationary 
generalized random field $(X_1(\varphi),\dots,X_d(\varphi))$,  
$\varphi\in{\cal S}$, such that $EX_j(\varphi)=0$ for all 
$\varphi\in{\cal S}$, and given two Shwartz functions 
$\varphi\in{\cal S}$ and $\psi\in{\cal S}$, the covariance
function $r_{j,j'}(\varphi,\psi)=EX_j(\varphi)X_{j'}(\psi)$ is 
given by the formula
\begin{equation}
r_{j,j'}(\varphi,\psi)=EX_j(\varphi)X_{j'}(\psi)
=\int\tilde\varphi(x)\,\bar{\!\tilde\psi}(x)G_{j,j'}(\,dx)  
\quad \textrm{for all } \varphi,\psi\in{\cal S}, \label{4.2}
\end{equation}
where $\,\tilde{}\,$ denotes Fourier transform, and $\,\bar{}\,$ is
complex conjugate. 

Formula~(\ref{4.2}) and the identity $EX_j(\varphi)=0$ for all 
$\varphi\in{\cal S}$ determine the distribution
of the vector valued, Gaussian stationary random field \hfill\break
$(X_1(\varphi),\dots,X_d(\varphi))$. 

Contrariwise, for 
all $1\le j,j'\le d$ the covariance function 
$EX_j(\varphi)X_{j'}(\psi)$, $\varphi,\psi\in{\cal S}$, 
determines the coordinate $G_{j,j'}$ of the positive 
semidefinite, even matrix $(G_{j,j'})$. $1\le j,j'\le d$, 
with moderately increasing distribution at infinity 
for which identity~(\ref{4.2}) holds.}

\medskip
Let me remark that the moderate decrease of the distribution of
the positive semidefinite matrix $(G_{j,j'})$, $1\le j,j'\le d$,
together with inequality~(\ref{3.2}) and the fast decrease of the 
functions $\varphi\in{\cal S}$ at infinity guarantee that the 
integral in~(\ref{4.2}) is convergent.

Condition~(\ref{4.1}) which we wrote in the definition of moderately
increasing positive semidefinite matrix valued measures appears
in the theory of distributions in a natural way. Such a
condition characterizes those measures which are distributions, 
i.e., continuous linear maps in the Schwartz space.

\medskip
In~\cite{9} we have proved with the help of some important results 
of Laurent Schwartz about distributions that in the 
case of scalar valued models, i.e., if $d=1$ the covariance 
function of every Gaussian stationary generalized random field 
with expectation zero agrees with the covariance function of 
a Gaussian stationary generalized random field constructed 
in the same way as we have done in Theorem~4.1. (In the case
$d=1$ the formulation of this result is simpler.) It seems 
very likely that a refinement of that argument would give 
the proof of an analogous statement in the general case. 
I did not investigate this question, because in the 
present paper we do not need such a result. 

\medskip\noindent
{\it Remark.}\/ Similarly to the case of vector valued 
stationary fields with discrete parameter we shall 
introduce the following terminology. If $(G_{j,j'})$, 
$1\le j,j'\le d$, is a positive semidefinite, matrix 
valued even measure with moderately increasing distribution 
at infinity, and there is a stationary generalized random field 
$(X_1(\varphi),\dots,X_d(\varphi))$, $\varphi\in{\cal S}$,
whose covariance function 
$$
r_{j,j'}(\varphi,\psi)=EX_j(\varphi)X_{j'}(\psi), \quad 
1\le j,j'\le d, \;\; \varphi,\psi\in{\cal S},
$$ 
satisfies relation~(\ref{4.2}) with this matrix valued 
measure $G$, then we call $G$ the matrix valued spectral 
measure of this covariance function $r_{j,j'}(\varphi,\psi)$. 
In general, we shall call a positive semidefinite matrix 
valued even measure on ${\mathbb R}^\nu$ with moderately increasing 
distribution at infinity a matrix valued spectral 
measure on ${\mathbb R}^\nu$. We have the right for such a 
terminology, because by Theorem~4.1 for any such matrix 
valued measure there exists a Gaussian stationary 
generalized random field such that this matrix valued 
measure is the matrix valued spectral measure of its 
covariance function. 

\medskip
Let me remark that the diagonal elements $G_{j,j}$ of the 
matrix valued spectral measure of the correlation function 
$r_{j,j'}(\varphi,\psi)$ of a vector valued stationary random 
field may be non finite measures on ${\mathbb R}^\nu$, they have 
to satisfy only relation~(\ref{4.1}). As a consequence, we can 
find a much richer class of matrix valued spectral measures by 
working with generalized random fields than by working only 
with classical stationary random fields. As we shall see, 
also vector valued random spectral measures corresponding 
to these matrix valued spectral measures can be constructed.
Actually we discussed vector valued stationary generalized 
random fields in this paper in order to construct this 
larger class of matrix valued spectral and vector valued
random spectral measures. 

\medskip\noindent
{\it Proof of Theorem 4.1.}\/ Let us observe that the function 
$r_{j,j'}(\varphi,\psi)$ defined in~(\ref{4.2}) is real valued. This can
be seen  by applying the change of variables $x\to -x$ in this
integral and by exploiting that $G_{j,j'}(-A)=\overline{G_{j,j'}(A)}$,
and  $\tilde\varphi(-x)=\bar{\tilde\varphi}(x)$, 
$\tilde\psi(-x)=\bar{\tilde\psi}(x)$, since this calculation 
yields that $r_{j,j'}(\varphi,\psi)=\overline{r_{j,j'}(\varphi,\psi)}$. 
Let us also remark that 
$r_{j,j'}(\varphi,\psi)=r_{j',j}(\psi,\varphi)$, since by 
formula~(\ref{4.2})
and the property $G_{j,j'}(A)=\overline{G_{j',j}(A)}$ of the matrix 
$(G_{j,j'}(A))$, $1\le j,j'\le d$, for all measurable sets 
$A\subset {\mathbb R}^\nu$ we have
$r_{j,j'}(\varphi,\psi)=\overline{r_{j',j}(\psi,\varphi)}$, and
we know that both side of this identity is real valued.

First we show that for all positive integers $N$ and 
functions $\varphi_k\in{\cal S}$, $1\le k\le N$, there 
are some Gaussian random vectors
$(X_1(\varphi_k),\dots,X_d(\varphi_k))$, $1\le k\le N$,
with expectation zero and 
covariances $EX_j(\varphi_k)X_{j'}(\varphi_{k'})
=r_{j,j'}(\varphi_k,\varphi_{k'})$
for all $1\le j,j'\le d$, $1\le k,k'\le N$, 
on an appropriate probability space, where 
$r_{j,j'}(\varphi_k,\varphi_{k'})$ is defined at the 
right-hand side of formula~(\ref{4.2}) with our matrix
valued measure $(G_{j,j'})$, $1\le j,j'\le d$, and
with the choice $\varphi=\varphi_k$, $\psi=\varphi_{k'}$.

We prove this statement if we show that the matrix with 
elements 
$$
d_{(j,k),(j',k')}=r_{j,j'}(\varphi_k,\varphi_{k'}), \quad 
1\le j,j'\le d, \;\; 1\le k,k'\le N,
$$ 
is positive semidefinite. To prove this result take any vector 
$(a_{j,k},\;1\le j\le d, 1\le k\le N)$, and observe that
\begin{eqnarray*}
&&\sum_{j=1}^d\sum_{j'=1}^d\sum_{k=1}^N\sum_{k'=1}^N a_{j,k}\overline{a_{j',k'}}
r_{j,j'}(\varphi_k,\varphi_{k'}) \\
&&\qquad=\sum_{j=1}^d\sum_{j'=1}^d\sum_{k=1}^N\sum_{k'=1}^N 
\int (a_{j,k}\tilde\varphi_k(x))
(\overline{{a_{j',k'}\tilde\varphi}_{k'}(x)})
g_{j,j'}(x)\mu(\,dx)  \\
&&\qquad=\sum_{j=1}^d\sum_{j'=1}^d\int \psi_j(x)\overline{\psi_{j'}(x)}
g_{j,j'}(x)\mu(\,dx)
=\int \psi(x)g(x)\overline{\psi(x)}\,\mu(\,dx)\ge0,
\end{eqnarray*}
where $\psi_j(x)=\sum\limits_{k=1}^N a_{j,k}\tilde\varphi_k(x)$,
$1\le j\le d$, $\psi(x)=(\psi_1(x),\dots,\psi_d(x))$, and $g(x)$
denotes the matrix $(g_{j,j'}(x))$, $1\le j,j'\le d$. In this calculation
we applied formula~(\ref{4.2}), the representation 
$G_{j,j'}(\,dx)=g_{j,j'}(x)\mu(\,dx)$ and finally the fact that $g(x)$
is a semidefinite matrix for $\mu$ almost all~$x$. 

Then it follows from Kolmogorov's existence theorem for random 
processes with consistent finite distributions that there is 
a Gaussian random field 
$$
(X_1(\varphi),\dots,X_d(\varphi)), \quad \varphi\in{\cal S},
$$ 
with zero expectation such that
$EX_j(\varphi)X_{j'}(\psi)=r_{j,j'}(\varphi,\psi)$ for all 
functions $\varphi\in{\cal S}$, $(\psi\in{\cal S}$ and 
$1\le j,j'\le d$. Besides, the finite dimensional distributions 
of this random field are determined because of the Gaussian 
property. Next we show that this random field is a vector 
valued generalized random field.

Property~(a) of the vector valued generalized random fields
follows from the following calculation.
\begin{eqnarray*}
&&E[a_1X_j(\varphi)+a_2X_j(\psi)
-X_j(a_1\varphi+a_2\psi)]^2\\
&&\qquad =\int\left(a_1\tilde\varphi(x)+a_2\tilde\psi(x)
-(\widetilde{a_1\varphi+a_2\psi)}(x)\right) \\
&&\qquad\qquad\qquad \times \left(\overline{a_1\tilde\varphi(x)}
+\overline{a_2\tilde\psi(x)}
-\overline{(\widetilde{a_1\varphi+a_2\psi)}(x)}\right)
G_{j,j}(\,dx)=0
\end{eqnarray*}
by formula~(\ref{4.2}) for all real numbers $a_1$, $a_2$, 
$1\le j\le d$ and $\varphi,\psi\in{\cal S}$.

Property~(b) of the vector valued generalized random fields
also holds for this model. Actually it is proved in~\cite{9} that 
if $\varphi_n\to\varphi$ in the topology of the space 
${\cal S}$, then 
$E[X_j(\varphi_n)-X_j(\varphi)]^2
=\int|\tilde\varphi_n(x)-\tilde\varphi(x)|^2 G_{j,j}(\,dx)\to0$ 
as $n\to\infty$, hence property~(b) also holds. (The proof is
not difficult. It exploits that for a sequence of functions
$\varphi_n\in{\cal S}^c$, $n=0,1,2,\dots$, $\varphi_n\to\varphi_0$
as $n\to\infty$ in the topology of ${\cal S}^c$ if and only if
$\tilde \varphi_n\to\tilde\varphi_0$ in the same topology.  
Besides, the measure $G_{j,j}$ satisfies inequality~(\ref{4.1}).)

It is also clear that the Gaussian random field constructed in 
such a way is stationary.

It remained to show that the covariance function 
$r_{j,j'}(\varphi,\psi)=EX_j(\varphi)X_{j'}(\psi)$
determines the complex measure $G_{j,j'}$. To show 
this we have to observe that inequality~(\ref{3.2}) holds also in 
this case, hence the Schwarz inequality implies that
$$
\int (1+|x|)^{-r}|g_{j,j'}(x)|\mu(\,dx)<\infty \quad\textrm{ for all }
1\le j,j'\le d 
$$
for a positive semidefinite matrix valued measure with moderately 
increasing distribution, i.e., this inequality holds not only for
$j=j'$. Then it follows from the standard theory of Schwartz
spaces that the class of Schwartz functions is sufficiently rich
to guarantee that the function $r_{j,j'}(\varphi,\psi)$ determines
the complex measure $G_{j,j'}$. Theorem~4.1 is proved.

\medskip
Next we construct a vector valued random spectral measure 
corresponding to a matrix valued spectral measure $(G_{j,j'})$,
$1\le j,j'\le d$, on~${\mathbb R}^\nu$. We argue similarly to Section~3, 
where the vector valued random spectral measures corresponding
to matrix valued spectral measures on $[-\pi,\pi)^\nu$ were 
considered. In the construction we shall also refer to some 
results in~\cite{9}.

Let us have a vector valued Gaussian stationary generalized
random field $X=(X_1(\varphi),\dots,X_d(\varphi))$, 
$\varphi\in{\cal S}$, $1\le j\le d$, with a matrix 
valued spectral measure $(G_{j,j'})$, $1\le j,j'\le d$.
First we define for all $1\le j\le d$ some (complex) 
Hilbert spaces ${\cal K}^c_{1,j}$, ${\cal H}^c_{1,j}$ and 
a norm preserving, invertible linear transformation 
$T_j$ between them in the following way. ${\cal K}^c_{1,j}$ 
consists of those complex valued functions $u(x)$ on 
${\mathbb R}^\nu$ for which $\int |u(x)|^2G_{j,j}(\,dx)<\infty$ 
with the scalar product 
$\langle u(x),v(x)\rangle=\int u(x)\overline{v(x)}G_{j,j}(\,dx)$.
To define the Hilbert space ${\cal H}^c_{1,j}$ let us first 
introduce the Hilbert space ${\cal H}={\cal H}^c$ of (complex 
valued) random variables with finite second moment on the 
probability space $(\Omega,\cal A,P)$ where our  stationary 
generalized random field is defined. We define the Hilbert
space ${\cal H}^c$ in the space consisting of these random
variables with the usual scalar product 
$\langle\xi,\eta\rangle=E\xi\bar\eta$ in ${\cal H}^c$. The 
Hilbert space ${\cal H}^c_{1,j}$ is defined  as the closure 
of the linear subspace of ${\cal H}^c$ consisting of the 
complex valued random variables $X_j(\varphi)+iX_j(\psi)$, 
$\varphi,\psi\in{\cal S}$. 

First we define the operator $T_j$ for functions of the 
form $\widetilde{\varphi+i\psi}$, 
$\varphi,\psi\in{\cal S}$. We define it by the formula
\begin{equation}
T_j(\widetilde{\varphi+i\psi)}=X_j(\varphi)+iX_j(\psi),
\quad \varphi,\,\psi\in{\cal S}. \label{4.3}
\end{equation}
A calculation, which was actually carried out in~\cite{9} 
shows that the set of functions
$\widetilde{\varphi+i\psi}$, $\varphi,\psi\in{\cal S}$, is 
dense in ${\cal K}^c_{1,j}$, and the transformation $T_j$, 
defined in~(\ref{4.3}) can be extended to a norm preserving, 
invertible linear transformation from ${\cal K}^c_{1,j}$ 
to ${\cal H}^c_{1,j}$. (In the calculation leading to this 
statement we apply formula~(\ref{4.2}) with the choice $j'=j$.)

Then we can define the random spectral measure $Z_{G,j}(A)$, 
similarly to the case discussed in Section~3, by the formula 
$Z_{G,j}(A)=T_j{\mathbb I}_A(\cdot))$ for all bounded 
measurable sets $A\subset {\mathbb R}^\nu$. To determine 
the joint distribution of the spectral measures $Z_{G,j}$ 
we make the following version of the corresponding 
argument in Section~3.

We define the following two Hilbert spaces ${\cal K}^c_1$
and ${\cal H}^c_1$ together with a norm preserving linear
transformation $T$ between them.

The elements of the Hilbert space ${\cal K}^c_1$ are 
the vectors $u=(u_1(x),\dots,u_d(x))$ with 
$u_j(x)\in{\cal K}^c_{1,j}$, $1\le j\le d$. We define
the scalar product on ${\cal K}^c_1$ with the 
help of the following positive semidefinite bilinear 
form $\langle \cdot,\cdot\rangle_0$. If 
$u(x)=(u_1(x),\dots,u_d(x))\in{\cal K}^c_1$  and 
$v(x)=(v_1(x),\dots,v_d(x))\in{\cal K}^c_1$, then
\begin{eqnarray*}
\langle u(x), v(x)\rangle_0
\!&=& \! \sum_{j=1}^d\sum_{j'=1}^d\int u_j(x)
\overline{ v_{j'}(x)}G_{j,j'}(\,dx) \\
\! &=& \! \sum_{j=1}^d\sum_{j'=1}^d\int g_{j,j'}(x)u_j(x)
\overline{ v_{j'}(x)}\mu(\,dx) 
=\int u(x) g(x) v(x)^*\mu(\,dx) 
\end{eqnarray*}
with the matrix $g(x)=(g_{j,j'}(x))$, $1\le j,j'\le d$, where 
$v^*(x)$ denotes the column vector whose elements are
the functions $\overline {v_{j'}(x)}$, $1\le j'\le d$. 
Actually, here we simply copied the corresponding
definition in Section~3 for the discrete time model,
and we can also prove that ${\cal K}^c_1$ is a Hilbert
space with the scalar $\langle\cdot,\cdot\rangle_0$ in the
same way as it was done in Section~3. 

The construction ${\cal H}^c_1$, and the proof of its
properties is again a simple copying of argument made
in Section~3. The elements of ${\cal H}^c_1$ are the
vectors $\xi=(\xi_1,\dots,\xi_d)$, where 
$\xi_j\in{\cal H}^c_{1,j}$, $1\le j\le d$, and we 
define the norm on it by means of the scalar product 
$\langle\xi,\eta\rangle_1=E\left(\sum_{j=1}^d\xi_j\right)
\overline{\left(\sum_{j=1}^d\eta_j\right)}$ for
$\xi=(\xi_1,\dots,\xi_d)\in{\cal H}^c_1$ and
$\eta=(\eta_1,\dots,\eta_d)\in{\cal H}^c_1$. We identify
two elements $\xi\in{\cal H}^c_1$ and $\eta\in{\cal H}^c_1$
if $\|\xi-\eta\|_1=0$. Then the argument of Section~3
yields that  ${\cal H}^c_1$ is a Hilbert space with the
scalar product $\langle\cdot,\cdot\rangle_1$.

We define the operator $T$ from ${\cal K}^c_1$ to ${\cal H}^c_1$ 
again in the same way as in Section~3. We define it by the
formula
$$
Tu=T(u_1,\dots,u_d)=(T_1u_1,\dots,T_du_d)
$$  
for 
$u=(u_1,\dots,u_d)$, $u_j\in{\cal K}^c_{1,j}$, with the help 
of the already defined operators $T_j$, $1\le j\le d$. We 
want to show that it is a norm preserving and invertible 
transformation from ${\cal K}^c_1$ to ${\cal H}^c_1$. Here
again we apply a similar, but sightly different argument 
from that in Section~3. We exploit that if we take the 
class of vectors 
$$
W=\{w=(u_1+iv_1,\dots,u_d+iv_d)\colon\; u_j\in{\cal S},\; 
v_j\in{\cal S} \textrm{ for all } 1\le j\le d\}
$$
then the class of vectors 
$$
\tilde W=\{(\widetilde{u_1+iv_1},\dots,\widetilde{u_d+iv_d})\colon\;
(u_1+iv_1,\dots, u_d+iv_d)\in W\}
$$ 
is an everywhere dense subspace of ${\cal K}^c_1$. and 
the class of vectors
$$
W(X)=\{((X_1(u_1+iv_1),\dots,X_d(u_d+iv_d)) \colon\; 
(u_1+iv_1,\dots,u_d+iv_d)\in W\}
$$ 
is an everywhere dense subspace of ${\cal H}^c_1$. (Here 
again the sign $\tilde{\ }$ denotes Fourier transform.)

Take two vectors 
$(u_{1,1}+iv_{1,1},\dots, u_{d,1}+iv_{d,1})\in W$ 
and $(u_{1,2}+iv_{1,2},\dots, u_{d,2}+iv_{d,2})\in W$.
The desired property of the operator~$T$ will follow from 
the following calculation:
\begin{eqnarray*}
&&\langle (\widetilde{u_{1,1}+iv_{1,1}},\dots,\widetilde{u_{d,1}+v_{d,1}}),
(\widetilde{u_{1,2}+iv_{1,2}},\dots,\widetilde{u_{d,2}+v_{d,2}})\rangle_0 \\
&&\qquad=\sum_{j=1}^d\sum_{j'=1}^d\int 
\widetilde{(u_{j,1}(x)+iv_{j,1}(x))} 
\overline{\widetilde{(u_{j',2}(x)+iv_{j',2}(x))}}G_{j,j'}(\,dx)\\
&&\qquad
=\sum_{j=1}^d\sum_{j'=1}^d E[X_j(u_{j,1})+iX_j(v_{j,1})]
[X_{j'}(u_{j',2})-iX_j(u_{j',2})] \\
&&\qquad=\langle (X_1(u_{1,1})+iX_1(v_{1,1}),\dots, 
X_d(u_{d,1})+iX_d(v_{d,1})), \\
&&\qquad\qquad\qquad 
(X_1(u_{1,2})+iX_1(v_{1,2}),\dots, X_d(u_{d,2})+iX_d(v_{d,2}))\rangle_1,
\end{eqnarray*}
i.e.,
\begin{eqnarray*}
&&\langle (\widetilde{u_{1,1}+iv_{1,1}},\dots,\widetilde{u_{d,1}+v_{d,1}}),
(\widetilde{u_{1,2}+iv_{1,2}},\dots,\widetilde{u_{d,2}+v_{d,2}})\rangle_0 \\
&&\qquad=\langle (T_1(u_{1,1}+iv_{1,1}),\dots, 
T_d(u_{d,1}+iv_{d,1})), \\
&& \qquad \qquad \qquad \qquad
(T_1(u_{1,2}+iv_{1,2}),\dots, T_d(u_{d,2}+iv_{d,2}))\rangle_1.
\end{eqnarray*}
This means that the operator $T$ maps the everywhere dense
subspace $\tilde W$ of ${\cal K}^c_1$ to the everywhere 
dense subspace $W(X)$ of ${\cal H}^c_1$ in a norm preserving
form. This implies that $T$ is a norm preserving, invertible 
transformation from ${\cal K}^c_1$ to ${\cal H}^c_1$.

Now we turn to the definition of the vector valued random
spectral measures corresponding to a matrix valued spectral
measure on ${\mathbb R}^\nu$. 

Let a vector valued, Gaussian stationary generalized 
random field 
$$
X(\varphi)=(X_1(\varphi),\dots,X_d(\varphi)), 
\quad\varphi\in{\cal S},
$$ 
be given with a matrix valued spectral measure 
$(G_{j,j'})$, $1\le j,j'\le d$, on ${\mathbb R}^\nu$. (We 
take such generalized, stationary random fields which 
were constructed in Theorem~4.1.) Let us consider the 
operators $T_j$, $1\le j\le d$, and $T$ constructed
above with the help of these quantities.  We define, 
similarly to the case of Gaussian stationary random 
fields with discrete parameters discussed in Section~3 
the random variables  $Z_{G,j}(A)=T_j({\mathbb I}_A(x))$ 
for all $1\le j\le d$ and bounded, measurable sets 
$A\subset {\mathbb R}^\nu$. (These functions 
${\mathbb I}_A(\cdot)$ are clearly elements of the 
Hilbert space ${\cal K}^c_{1,j}$ for all $\le j\le d$). 
It can be proved with the help of the properties of 
the operator~$T$ that these random functions satisfy 
properties (i)--(v) formulated in the definition of 
random spectral measures on the torus, considered in 
Section~3. The argument applied in Section~3 holds 
also in in this case. In particular, property~(v) 
can be proved with the help of property~(v$'$). 
Property~(v$'$) can be proved with some work, and 
actually this was done in~\cite{9}. We prove~(v$'$) 
by checking it first for functions $u\in{\cal S}^c$.

The above result makes natural the following definition of 
vector valued random spectral measures corresponding to a 
matrix valued spectral measure on ${\mathbb R}^\nu$. This is very
similar to the definition of vector valued random spectral
measures on the torus.

\medskip\noindent
{\bf Definition of vector valued random spectral measures on 
${\mathbb R}^\nu$.} {\it Let $G=(G_{j,j'})$, $1\le j,j'\le d$, be 
a matrix valued spectral measure on ${\mathbb R}^\nu$. We call 
a set of complex valued random variables $Z_{G,j}(A)$
depending on pairs $(j,A)$, where $1\le j\le d$,  
$A\in{\cal A}$, and ${\cal A}$ is the algebra  
$$
{\cal A}=\{ A\colon\; A 
\textrm{ is a bounded Borel measurable set in } {\mathbb R}^\nu\},
$$ 
a $d$-dimensional vector valued random spectral measure 
corresponding to the matrix valued spectral measure~$G$ 
on ${\mathbb R}^\nu$ if this set of random variables 
$Z_{G,j}(A)$, $1\le j\le d$, $A\in{\cal A}$, satisfies 
properties~(i)--(v) introduced in Section~3 in the 
definition of vector valued random spectral measures on 
the torus. Given a fixed index $1\le j\le d$, we call the
set of random variables $Z_{G,j}(A)$, $A\in{\cal A}$,
with this index $j$ the $j$-th coordinate of this matrix
valued spectral measure, and we denote it by~$Z_{G,j}$.
We denote a $d$-dimensional vector valued random spectral 
measure corresponding to the matrix valued spectral 
measure~$G$ by $Z_G=(Z_{G,1},\dots,Z_{G,d})$. }

\medskip
We can show with the help of the arguments applied in
Section~3 that for any $d$-dimensional matrix valued spectral 
measure on ${\mathbb R}^\nu$ there exists a $d$-dimensional vector 
valued random spectral measure corresponding to it. 

We can define the random integral $\int f(x)Z_{G,j}(\,dx)$ of 
the functions $f\in{\cal K}^c_{1,j}$ with respect to the random 
spectral measure $Z_{G,j}$, $1\le j\le d$, corresponding to 
the matrix valued spectral measure $(G_{j,j'})$, $1\le j,j'\le d$, 
of a Gaussian stationary generalized field in the same way as we
defined these random  integrals with respect to random spectral
measures corresponding to a spectral measures on the torus
$[-\pi,\pi)^\nu$ in Section~3. First we define these integrals for 
elementary functions which are defined in the same way as
it was done in Section~3. Then following the calculation
of that section we can define these integrals for a general
function $f\in{\cal K}^c_{1,j}$, and it can be seen that
formulas (\ref{3.7}), (\ref{3.8}) and (\ref{3.9}) remain 
valid for them. In particular, the random integrals 
$\int \tilde \varphi(x)Z_{G,j}(\,dx)$ are (meaningful and) real 
valued random variables for all $\varphi\in{\cal S}$, and
$$
E \left(\int \tilde \varphi(x)Z_{G,j}(\,dx)
 \int\bar{\tilde \psi}(x)Z_{G,{j'}}(\,dx)\right)
=\int \tilde\varphi(x)\bar{\tilde \psi}(x) G_{j,j'}(\,dx)
$$
for all $\varphi,\psi\in {\cal S}$ and $1\le j,j'\le d$. This 
identity together with relation~(\ref{3.7}) and the fact that the 
above considered random integrals are linear operators imply 
that the set of random variables
\begin{equation}
X_j(\varphi)=\int\tilde\varphi (x) Z_{G,j}(\,dx), 
\quad \varphi\in{\cal S}, \;\; 1\le j\le d, \label{4.4}
\end{equation}
constitute a vector valued Gaussian, stationary generalized
random field with spectral measure $(G_{j,j'})$, $1\le j,j'\le d$.

This implies that the natural version of Theorem~3.1 remains
valid if we consider a matrix valued spectral measure~$(G_{j,j'})$,
$1\le j,j'\le d$, on ${\mathbb R}^\nu$. Then there exists a random spectral
measure $Z_G=(Z_{G,1},\dots,Z_{G,d})$ corresponding to it, and
we have defined the random integrals $\int u(x)Z_{G,j}(\,dx)$,
$1\le j\le d$, with respect to it for all $u\in{\cal K}^c_{1,j}$. 
The class of random variables, $X_j(\varphi)$, 
$\varphi\in{\cal S}$, $1\le j\le d$, defined in~(\ref{4.4})
constitute a vector valued, Gaussian stationary generalized 
random field with matrix valued spectral measure $(G_{j,j'})$, 
$1\le j,j'\le d$. Moreover, if a $d$-dimensional vector 
valued Gaussian stationary random field is given with 
spectral measure $(G_{j,j'})$, $1\le j,j'\le d$, then we can 
consider  the random spectral measure $(Z_{G,1},\dots,Z_{G,d})$ 
constructed in this section with the help of this random 
field. This random spectral measure has the property that 
the random field given by the random integrals defined 
in formula~(\ref{4.4}) with their help agrees with the original 
vector valued Gaussian stationary generalized random field.

We can formulate a natural version of Lemma~3.2 where we
consider a matrix valued spectral measure $(G_{j,j'})$, 
$1\le j,j'\le d$, on ${\mathbb R}^\nu$ instead of a matrix valued 
spectral measure on the torus $[-\pi,\pi)^\nu$. In this
version of Lemma~3.2 we define ${\cal K}_{1,j}$ as
$$
{\cal K}_{1,j}=\left\{u\colon \int |u(x)|^2G_{j,j}(\,dx)<\infty, \quad 
u(-x)=\overline{u(x)}\textrm{ for all } x\in {\mathbb R}^\nu\right\},
$$
with the scalar product
$\langle u,v\rangle=\int u(x)\overline{v(x)}G_{j,j}(\,dx)$, 
$u,v\in{\cal K}_{1,j}$, and ${\cal H}_{1,j}$ as the closure of 
the linear space consisting of the finite linear combination 
of the random variables $X_j(\varphi)$, $\varphi\in{\cal S}$, 
with real coefficients in the Hilbert space ${\cal H}$. 
This version of Lemma~3.2 states that ${\cal K}_{1,j}$ and 
${\cal H}_{1,j}$ are real Hilbert spaces, and 
$T_j(u)=\int u(x)Z_{G,j}(\,dx)$ is a norm preserving and
invertible transformation from ${\cal K}_{1,j}$ to 
${\cal H}_{1,j}$.

The proof of this version of Lemma~3.2 is very similar to the proof 
of the original lemma. The main difference is that now we show that
the class of functions $\tilde\varphi$ with $\varphi\in{\cal S}$ is
a dense linear subspace of ${\cal K}_{1,j}$, and the transformation
$T_j(\tilde{\varphi})=\int \tilde{\varphi}(x)Z_{G,j}(\,dx)=X_j(\varphi)$,
$\varphi\in{\cal S}$, is a norm preserving transformation from 
an everywhere dense subspace of ${\cal K}_{1,j}$ to an everywhere
dense subspace of ${\cal H}_{1,j}$.

The natural version of Lemma~3.3 also holds. It states that a
matrix valued spectral measure $(G_{j,j'})$, $1\le j,j'\le d$, 
on ${\mathbb R}^\nu$ determines the distribution of a vector 
valued random spectral measure $Z_{G,j}$, $1\le j\le d$, 
corresponding to it. The proof of this version is the same as 
the proof of the original lemma. The only difference is that 
now we consider the random spectral measure $Z_{G,j}(A)$ for 
all measurable, bounded sets $A\subset {\mathbb R}^\nu$.

Finally I would remark that property (vi) of the random
spectral measures also remains valid for this new class of 
random spectral measures, because its proof applies only 
properties (i)--(v) of random spectral measures.

\vfill\eject

\section{Multiple Wiener--It\^o integrals with respect to vector 
valued random spectral measures}

Next we want to rewrite the random variables with finite second
moments which are measurable with respect the $\sigma$-algebra
generated by the elements of a vector valued Gaussian stationary 
random field  in an appropriate form, which enables us  to rewrite 
also the random sums defined in~(\ref{1.1}) in a form that 
helps in the study of their limit behaviour. In the scalar 
valued case, i.e., when $d=1$ we could do this with the help of 
multiple Wiener--It\^o integrals. We could rewrite the random 
sums~(\ref{1.1}) with their help in such a form that provided 
great help in the study of the limit theorems we were interested 
in. Next we show that a similar method can be applied also in 
the case of vector valued Gaussian stationary random fields. To 
do this first we have to define the multiple  Wiener--It\^o 
integrals also in the vector valued case. We start the 
definition of multiple  Wiener--it\^o integrals in this case 
with the introduction of the following notation.

\medskip
Let $X(p)=(X_1(p),\dots,X_d(p))$, $EX(p)=0$, 
$p\in{\mathbb Z}^\nu$, be a vector valued 
stationary Gaussian random field with some matrix 
valued spectral measure 
$G=(G_{j,j'})$, $1\le j,j'\le d$. Let
$Z_G=(Z_{G,1},\dots,Z_{G,d})$ be a
vector valued random spectral measure 
corresponding to it which is chosen in such
a way that $X_j(p)=\int e^{i(p,x)}Z_{G,j}(\,dx)$ for all 
$p\in {\mathbb Z}^\nu$ and $1\le j\le d$. 
Let us consider 
the (real) Hilbert space ${\cal H}$ of square integrable 
random variables measurable with respect to the 
$\sigma$-algebra generated by the random vectors $X(p)$, 
$p\in{\mathbb Z}^\nu$. More generally, let us consider a 
(possibly generalized) matrix valued spectral measure 
$G=(G_{j,j'})$, $1\le j,j'\le d$, and a vector valued random 
spectral measure $Z_G=(Z_{G,1},\dots,Z_{G,d})$ corresponding 
to it, where the matrix valued spectral measures $G_{j,j'}$ 
and vector valued random spectral measures $Z_{G,j}$ are 
defined  either on the torus $[-\pi,\pi)^\nu$ or on ${\mathbb R}^\nu$, 
and consider the (real) Hilbert space ${\cal H}$ of the 
square integrable (real valued) random variables, 
measurable with respect to the $\sigma$-algebra generated 
by the random variables of the vector valued random 
spectral measures $Z_G$ with the usual scalar product 
in this space. We would like to write the elements of 
the Hilbert space ${\cal H}$ in the form of a sum of 
multiple Wiener--It\^o integrals with respect to the 
vector valued random spectral measure $Z_G$. I shall 
construct these Wiener--It\^o integrals in this section, 
and I prove some of their important properties. 

As a discussion in Section~2 of \cite{11} will show we 
cannot write all elements of ${\cal H}$ in the form of 
a sum of Wiener--It\^o integrals, but we can do this for 
the elements of an everywhere dense subspace of ${\cal H}$.
In particular, if we consider finitely many random variables
$X_j(p)$, $1\le j\le d$, $p\in{\mathbb Z}^\nu$ of a discrete
or $X_j(\varphi)$, $1\le j\le d$, $\varphi\in{\cal S}^\nu$,
of a generalized vector valued stationary Gaussian 
random field, then all polynomials of these random
variables can be written as the sum of Wiener--It\^o
integrals. Such a result will be sufficient for our
purposes. In the subsequent discussion I impose a 
technical condition about the properties of the matrix 
valued spectral measure $G=(G_{j,j'})$ I shall be working 
with. I assume that it is non-atomic. More precisely, 
I assume that we are working with such a dominating 
measure $\mu$ for the coordinates of the matrix valued 
spectral measures $G_{j,j'}$ for which 
$\mu(\{x\})=0$ for all $x\in\ {\mathbb R}^\nu$.

First I define for all $n=1,2,\dots$ and $1\le j_s\le d$ 
for the indices $1\le s\le n$ the $n$-fold multiple 
Wiener--It\^o integral
$$
I_n(f|j_1,\dots,j_n)=\int f(x_1,\dots,x_n)
Z_{G,{j_1}}(\,dx_1)\dots Z_{G,{j_n}}(d\,x_n)
$$
with respect to the coordinates of a vector valued random 
spectral measure $Z_G=(Z_{G,1},\dots,Z_{G,d})$, corresponding 
to a matrix valued spectral measure $G=(G_{j,j'})$,
$1\le j,j'\le d$. I shall define these Wiener--It\^o 
integrals with kernel functions $f\in{\cal K}_{n,j_1,\dots,j_n}$ 
in a (real) Hilbert space 
${\cal K}_{n,j_1,\dots,j_n}
={\cal K}_{n,j_1,\dots,j_n}(G_{j_1,j_1},\dots,G_{j_n,j_n})$ 
defined below.

We define ${\cal K}_{n,j_1,\dots,j_n}
={\cal K}_{n,j_1,\dots,j_n}(G_{j_1,j_1}.\dots,G_{j_n,j_n})$ as the 
Hilbert space consisting of those complex valued functions
$f(x_1,\dots,x_n)$ on ${\mathbb R}^{n\nu}$ which satisfy the following 
relations~(a) and~(b):

\medskip
\begin{description}
\item[{\rm(a)}] $f(-x_1,\dots,-x_n)=\overline{f(x_1,\dots,x_n)}$ 
for all $(x_1,\dots,x_n)\in {\mathbb R}^{n\nu}$,
\item[{\rm(b)}]
$\|f\|^2=\int|f(x_1,\dots,x_n)|^2
G_{j_1,j_1}(\,dx_1)\dots G_{j_n,j_n}(\,dx_n)<\infty$.
\end{description}

\medskip\noindent
We define the scalar product in ${\cal K}_{n,j_1,\dots,j_n}$ in the 
following way. If $f,\,g\in{\cal K}_{n,j_1,\dots,j_n}$, then
\begin{eqnarray*}
\langle f,g\rangle&=&\int f(x_1,\dots,x_n)\overline{g(x_1,\dots,x_n)}
G_{j_1,j_1}(\,dx_1)\dots G_{j_n,j_n}(\,dx_n)\\
&=&\int f(x_1,\dots,x_n)g(-x_1,\dots,-x_n)
G_{j_1,j_1}(\,dx_1)\dots G_{j_n,j_n}(\,dx_n).
\end{eqnarray*}

Because of the symmetry $G_{j_s,j_s}(A)=G_{j_s,j_s}(-A)$ of the 
spectral measure
$\langle f,g\rangle=\overline{\langle f,g\rangle}$, i.e., the 
scalar product $\langle f,g\rangle$ is a real
number for all $f,\,g\in{\cal K}_{n,j_1,\dots,j_n}$. This means that
${\cal K}_{n,j_1,\dots,j_n}$ is a real Hilbert space, as I claimed.
We also define the real Hilbert space ${\cal K}_0$ for $n=0$ as
the space of real constants with the norm $\|c\|=|c|$. 

\medskip\noindent
{\it Remark.}\/ In the case $n=1$ the above defined real Hilbert 
space ${\cal K}_{1,j}$ agrees with the real Hilbert 
space~${\cal K}_{1,j}$ introduced in Lemma~3.2.

\medskip
Similarly to the scalar valued case, first we introduce so-called
simple functions and define the multiple integrals for them. We
prove some properties of this integral which enable us to extend
its definition by means of an $L_2$ extension for all functions 
$f\in{\cal K}_{j_1,\dots,j_n}$. We define the class of simple 
functions together with the notion of regular systems.

\medskip\noindent
{\bf Definition of regular systems and of the class of simple 
functions.} {\it Let 
$$
{\cal D}=\{\Delta_k,\;k=\pm1,\pm2,\dots,\pm N\}
$$ 
be a finite collection of bounded, measurable sets in ${\mathbb R}^\nu$  
indexed by the integers $\pm1$,\dots, $\pm N$ with some positive
integer~$N$. We say that ${\cal D}$ is a regular system if
$\Delta_k=-\Delta_{-k}$, and $\Delta_k\cap\Delta_l=\emptyset$
if $k\neq l$ for all $k,l=\pm1,\pm2,\dots,\pm N$. A function
$f\in{\cal K}_{n,j_1,\dots,j_n}$ is adapted to this system 
${\cal D}$ if $f(x_1,\dots,x_n)$ is constant on the sets
$\Delta_{k_1}\times\Delta_{k_2}\times\cdots\times\Delta_{k_n}$, \
$k_l=\pm1,\dots,\pm N$, $l=1,2,\dots,n$,  it vanishes outside these
sets, and it also vanishes on those sets of the above form for which 
$k_l=\pm k_{l'}$ for some $l\neq l'$. 

A function $f\in{\cal K}_{n,j_1,\dots,j_n}$ is in the class
$\hat{\cal K}_{n,j_1,\dots,j_n}$ of simple functions if it is adapted to
some regular system ${\cal D}=\{\Delta_k,\;k=\pm1,\dots,\pm N\}$.} 

\medskip\noindent
{\bf Definition of Wiener--It\^o integrals of simple functions.}
{\it Let a simple function $f\in \hat{\cal K}_{n,j_1,\dots,j_n}$ 
be adapted  to some regular system 
$$
{\cal D}=\{\Delta_k,\;k=\pm1,\dots,\pm N\}.
$$ 
Its $n$-fold Wiener--It\^o integral with respect to 
$Z_G=(Z_{G,1},\dots,Z_{G,d})$ with parameters 
$j_1,\dots,j_n$, $1\le j_k\le d$ for all 
$1\le k\le n$, is defined as
\begin{eqnarray}
&&\int f(x_1,\dots,x_n)Z_{G,{j_1}}(\,dx_1)\dots Z_{G,{j_n}}(\,dx_n) 
\label{5.1} \\
&&\qquad =I_n(f|j_1,\dots,j_n) \nonumber  \\
&&\qquad=\sum_{\substack{k_l=\pm1,\dots,\pm N\\l=1,2,\dots,n}}
f(u_{k_1},\dots,u_{k_n})
Z_{G,{j_1}}(\Delta_{k_1})\cdots Z_{G,{j_n}}(\Delta_{k_n}),
\nonumber
 \end{eqnarray}
where $u_k\in\Delta_k$, $k=\pm1,\dots,\pm N$.} 

\medskip\noindent
Although the regular system ${\cal D}$ to which $f$ is adapted is 
not uniquely determined (for example the elements of ${\cal D}$ can 
be divided to smaller sets), the integral defined in~(\ref{5.1}) 
is meaningful, i.e., its value does not depend on the choice of 
${\cal D}$. This can be proved with the help of property~(iv) 
of vector valued random spectral measures defined in Section~3 
in the same way as it was done in the scalar valued 
case in~\cite{9}. (Let me also remark that here I defined the random 
integral $I_n(f|j_1,\dots,j_n)$ with a normalization different 
from the normalization of the corresponding expression $I_G(f)$ 
introduced in~\cite{9}. Here I omitted the norming term $\frac1{n!}$.)

\medskip
Because of the definition of simple functions the sum in~(\ref{5.1})
does not change if we allow in it summation only for such sequences
$k_1,\dots,k_n$ for which $k_l\neq\pm k_{l'}$ if $l\neq l'$. This
fact will be exploited in the subsequent considerations.

\medskip
Next I formulate some important properties about the 
Wiener--It\^o integrals of simple functions. Later we shall 
see that these properties remain valid in the general case.

\begin{equation}
I_n(f|j_1,\dots,j_n) 
\textrm{ is a real valued random variable for all } 
f\in\hat{\cal K}_{n,j_1,\dots,j_n}. \label{5.2}
\end{equation}
Indeed, $I_n(f|j_1,\dots,j_n)=\overline{I_n(f|j_1,\dots,j_n)}$
by Property~(a) of the functions in ${\cal K}_{n,j_1,\dots,j_n}$ 
and property~(v) of the random spectral measures defined in 
Section~3, hence~(\ref{5.2}) holds. It is also clear that  
$\hat{\cal K}_{n,j_1,\dots,j_n}$ is a linear space, and the mapping
$f\to I_n(f|j_1,\dots,j_n)$ is a linear transformation on it.

The relation
\begin{equation}
EI_n(f|j_1,\dots,j_n)=0 \quad 
\textrm{for }f\in\hat{\cal K}_{n,j_1,\dots,j_k} \quad
\textrm { if } n\neq0 \label{5.3}
\end{equation}
also holds. (In the non-zero terms of the sum in~(\ref{5.1}) we have 
the product of independent random variables with expectation 
zero by property~(vi) of the random spectral measures described
also in Section~3.) Next I express the covariance 
between random variables of the form  $I_n(f|j_1,\dots,j_n)$. 
To do this first I introduce the following notation. Let 
$\Pi(n)$ denote the set of all permutations of the set 
$\{1,\dots,n\}$, and let $\pi=(\pi(1),\dots,\pi(n))$ denote 
one of its element. 

Let us have a positive integer $n\ge1$, and two sequences
$j_1,\dots,j_n$ and $j'_1,\dots,j'_n$, $1\le j_s,j'_s\le d$ for all 
$1\le s\le d$. Let $f\in\hat{\cal K}_{n,j_1,\dots,j_n}$ and 
$h\in\hat{\cal K}_{n,j'_1,\dots,j'_n}$. I shall show that
\begin{eqnarray}
&&EI_n(f|j_1,\dots,j_n)I_n(h|j'_1,\dots,j'_n) \label{5.4} \\
&&\qquad =\sum_{\pi\in\Pi(n)} \int f(x_1,\dots x_n) 
\overline{h(x_{\pi(1)},\dots,x_{\pi(n)})} \nonumber \\
&&\qquad\qquad\qquad 
G_{j_1,j'_{\pi^{-1}(1)}}(\,dx_1)\dots G_{j_n,j'_{\pi^{-1}(n)}}(\,dx_n). 
\nonumber
\end{eqnarray}
On the other hand, if $n\neq n'$, and 
$f\in\hat{\cal K}_{n,j_1,\dots,j_n}$,
$h\in\hat{\cal K}_{n',j'_1,\dots,j'_{n'}}$, then
\begin{equation}
EI_n(f|j_1,\dots,j_n)I_{n'}(h|j'_1,\dots,j'_{n'})=0. \label{5.5}
\end{equation}

Next I show the following inequality with the help of formula~(\ref{5.4}):
\begin{eqnarray}
E|I_n(f|j_1,\dots,j_n)|^2 &\le& n! \int |f(x_1,\dots x_n)|^2
G_{j_1,j_1}(\,dx_1)\dots G_{j_n,j_n}(\,dx_n) \nonumber \\
&=&n!\|f_{n,j_1,\dots,j_n}\|^2 \label{5.6}
\end{eqnarray}
for all $f\in\hat{\cal K}_{n,j_1,\dots,j_n}$.

Indeed we get by applying (\ref{5.4}) for 
$f=h\in\hat{\cal K}_{n,j_1,\dots,j_n}$
together with relation~(\ref{3.2}) that
\begin{eqnarray}
&&\!\!\!\!\!\! E|I_n(f|j_1,\dots,j_n)|^2 
 \le\sum_{\pi\in\Pi(n)} \int |f(x_1,\dots x_n)|
|f(x_{\pi(1)},\dots,x_{\pi(n)})|  \label{5.7} \\
&&\qquad\quad \times \prod_{s=1}^n 
\left(g_{j_s,j_s}(x_s) g_{j_{\pi^{-1}(s)},j_{\pi^{-1}(s)}}(x_s)\right)^{1/2}
\mu(\,dx_1)\dots\mu(\,dx_n). \nonumber
\end{eqnarray}
On the other hand, we get with the help of the Schwarz inequality that
\begin{eqnarray}
&&\!\!\!\!\!\!\!\!\!\!\!\!\!\!
\int |f(x_1,\dots x_n)||f(x_{\pi(1)},\dots,x_{\pi(n)})|
\prod_{s=1}^n \left(g_{j_s,j_s}(x_s) g_{j_{\pi^{-1}(s)},j_{\pi^{-1}(s)}}(x_s)\right)^{1/2}
\nonumber  \\
&&\qquad\qquad\qquad\qquad\qquad\qquad \qquad\qquad\qquad
\times \mu(\,dx_1)\dots \mu(\,dx_n) 
\label{5.8} \\
&&\!\!\!\!\!\!\le \left(\int |f(x_1,\dots x_n)|^2
\prod_{s=1}^n g_{j_s,j_s}(x_s)\mu(\,dx_1)\dots \mu(\,dx_n)\right)^{1/2} 
\nonumber \\ 
&&\times \left(\int |f(x_{\pi(1)},\dots,x_{\pi(n)})|^2
\prod_{s=1}^n g_{j_{\pi^{-1}(s)},j_{\pi^{-1}(s)}}(x_s)
\times \mu(\,dx_1)\dots \mu(\,dx_n)\right)^{1/2}  
\nonumber 
\end{eqnarray}
for all $\pi\in\Pi(n)$. Let us also observe that the map $T$ 
from ${\mathbb R}^{n\nu}$ to ${\mathbb R}^{n\nu}$, defined as 
$$
T(x_1,\dots,x_n)=(x_{\pi(1)},\dots,x_{\pi(n)})
$$ 
is a bijection, and it is a measure preserving transformation from 
$$
({\mathbb R}^{n\nu}, G_{j_1,j_1}\times\cdots\times G_{j_n,j_n})=
({\mathbb R}^{n\nu}, g_{j_1,j_1}(x_1)\cdots g_{j_n,j_n}(x_n)\mu(\,dx_1)
\dots\mu(\,dx_n)\,)
$$
to 
\begin{eqnarray*}
&&({\mathbb R}^{n\nu}, G_{j_{\pi^{-1}(1)} ,j_{\pi^{-1}(1) } }\times\cdots\times 
G_{j_{\pi^{-1}(n)},j_{\pi^{-1}(n)}})\\
&&\qquad=({\mathbb R}^{n\nu}, g_{j_{\pi^{-1}(1)},j_{\pi^{-1}(1)}}(x_1)
\cdots g_{j_{\pi^{-1}(n)},j_{\pi^{-1}(n)}}(x_n)\mu(\,dx_1)\dots\mu(\,dx_n)\,).
\end{eqnarray*}
To see this it is enough to check that if $A=A_1\times\cdots\times A_n$,
then 
$$
(G_{1,1}\times \cdots\times G_{n,n})(A)=\prod_{l=1}^nG_{l,l}(A_l),
$$
$TA=A_{\pi^{-1}(1)}\times\cdots\times A_{\pi^{-1}(n)}$,
\begin{eqnarray*}
&&(G_{j_{\pi^{-1}(1)} ,j_{\pi^{-1}(1) } }\times\cdots\times 
G_{j_{\pi^{-1}(n)},j_{\pi^{-1}(n)}})(TA) \\
&&\qquad=\prod_{l=1}^n
G_{j_{\pi^{-1}(l)},j_{\pi^{-1}(l)}}(A_{\pi^{-1}(l)})=
(G_{1,1}\times \cdots\times G_{n,n})(A).
\end{eqnarray*}
The last identity together with the bijective property of $T$
imply that it is measure preserving.

Because of the measure preserving property of the operator~$T$ 
we can write that
\begin{eqnarray}
&&\int |f(x_1,\dots x_n)|^2 \prod_{s=1}^n g_{j_s,j_s}(x_s)
\mu( \,dx_1)\dots \mu(\,dx_n)  
\label{5.9} \\
&&\qquad=  \int |f(x_{\pi(1)},\dots,x_{\pi(n)})|^2
\prod_{s=1}^n g_{j_{\pi^{-1}(s)},j_{\pi^{-1}(s)}}(x_s)\mu( \,dx_1)\dots\mu( \,dx_n).
\nonumber 
\end{eqnarray}
Relation (\ref{5.6}) follows from relations 
(\ref{5.7}), (\ref{5.8}) and (\ref{5.9}).

To prove formulas (\ref{5.4}) and (\ref{5.5}) first we 
prove the following relations. Let a regular system
${\cal D}=\{\Delta_k,\;k=\pm1,\pm2,\dots,\pm N\}$ be given, choose an
integer $n\ge 1$, some numbers $j_1,\dots,j_n$ and $j_1'\dots,j'_n$ 
such that $1\le j_s,j'_s\le d$, $1\le s\le d$, together with two 
sequences of numbers $k_1,\dots,k_n$ and $l_1,\dots,l_n$ such that
$k_s,l_s\in\{\pm1,\dots,\pm N\}$ for all $1\le s\le n$, and they also
satisfy the relation $k_s\neq\pm k_{s'}$, and $l_s\neq\pm l_{s'}$ if 
$s\neq s'$. I claim that under these conditions
\begin{equation}
EZ_{G,{j_1}}(\Delta_{k_1})\cdots Z_{G,{j_n}}(\Delta_{k_n})
\overline{Z_{G,{j'_1}}(\Delta_{l_1})\cdots Z_{G,{j'_n}}(\Delta_{l_n})}=0 
\label{5.10}
\end{equation}
if $\{k_1,\dots,k_n\}\neq\{l_1,\dots,l_n\}$. On the other hand,
if 
\begin{equation}
l_p=k_{\pi(p)} \;\;\textrm{for all } 1\le p\le n \label{5.11} 
\end{equation}
with some permutation $\pi\in\Pi(n)$, then
\begin{eqnarray}
&&EZ_{G,{j_1}}(\Delta_{k_1})\cdots Z_{G,{j_n}}(\Delta_{k_n})
\overline{Z_{G,{j'_1}}(\Delta_{l_1})\cdots Z_{G,{j'_n}}(\Delta_{l_n})} 
\nonumber \\
&&\qquad= G_{j_1,j'_{\pi^{-1}(1)}}(\Delta_{k_1})\cdots G_{j_n,j'_{\pi^{-1}(n)}}
(\Delta_{k_n}).
\label{5.12}
\end{eqnarray}
Let me remark that there cannot be two different permutations
$\pi\in\Pi(n)$ satisfying relation~(\ref{5.11}), since by our
assumption also elements of the set $\{k_1,\dots,k_n\}$ are 
different, and the same relation holds for the set
$\{1_1,\dots,l_n\}$.

To prove (\ref{5.10}) we show that under its conditions the product
$$
Z_{G,{j_1}}(\Delta_{k_1})\cdots Z_{G,{j_n}}(\Delta_{k_n})
\overline{Z_{G,{j'_1}}(\Delta_{l_1})\cdots Z_{G,{j'_n}}(\Delta_{l_n})}
$$
can be written in the form of a product of two independent 
terms in such a way that one of them has expectation zero. 

Indeed, since $\{k_1,\dots,k_n\}\neq\{l_1,\dots,l_n\}$, there 
is such an element $k_s$ for which $k_s\neq l_t$ for all 
$1\le t\le n$, and also the relation $k_s\neq \pm k_t$ if 
$s\neq t$, holds. If the relation $k_s\neq \pm l_t$ also
holds for all $1\le t\le n$, then $Z_{G,{j_s}}(\Delta_{k_s})$ 
is independent of the product of the product of the 
remaining terms in this product because of property~(vi) 
of vector valued random spectral measures given in 
Section~3, and $EZ_{G,{j_s}}(\Delta_{k_s})=0$. Hence 
relation~(\ref{5.10}) holds in this case.

In the other case, there is an index $s'$ such that 
$l_{s'}=-k_s$. In this case the vector 
\begin{eqnarray*}
(Z_{G,{j_s}}(\Delta_{k_s}),\overline{Z_{G,{j_{s'}}}(\Delta_{l_{s'}})})
&=&(Z_{G,{j_s}}(\Delta_{k_s}),Z_{G,{j_{s'}}}(-\Delta_{l_{s'}})) \\
&=&(Z_{G,{j_s}}(\Delta_{k_s}),Z_{G,{j_{s'}}}(\Delta_{k_{s}}))
\end{eqnarray*}
is independent of the remaining terms, (because of 
property~(vi) of the vector valued random spectral measures). 
In the last relation we exploited that 
$-\Delta_{l_{s'}}=\Delta_{k_s}$). Hence
$$
EZ_{G,{j_s}}(\Delta_{k_s})\overline{Z_{G,{j_{s'}}}(\Delta_{l_{s'}})}=
EZ_{G,{j_s}}(\Delta_{k_s})\overline{Z_{G,{j_{s'}}}(-\Delta_{k_s})}=0,
$$
and relation (\ref{5.10}) holds in this case, too.

To prove (\ref{5.12}) let us observe that under its condition 
the investigated product can be written in the form
\begin{eqnarray*}
&&Z_{G,{j_1}}(\Delta_{k_1})\cdots Z_{G,{j_n}}(\Delta_{k_n})
\overline{Z_{G,{j'_1}}(\Delta_{l_1})\cdots Z_{G,{j'_n}}(\Delta_{l_n})} \\
&&\qquad =\prod_{p=1}^n Z_{G,{j_p}}(\Delta_{k_p}) 
\overline{Z_{G,{ j'_{\pi^{-1}(p)}}}(\Delta_{k_p})}.
\end{eqnarray*}
The terms in the product at the right-hand side are independent 
for different indices~$s$, and
$EZ_{G,{j_p}}(\Delta_{k_p})\overline{Z_{G,{ j'_{\pi^{-1}(p)}}}(\Delta_{k_p})}
=G_{j_p,j'_{\pi^{-1}(p)}}(\Delta_{k_p})$. Formula~(\ref{5.12}) follows from 
these relations and the independence between the terms in the
last product. (Here we use again property~(vi) of the random
spectral measures.)

To prove formula (\ref{5.4}) let us take a regular system 
$$
{\cal D}=\{\Delta_k,\;k=\pm1,\dots,\pm N\}
$$ 
such that both functions $f$ and $h$ are adapted to it. This can 
be done by means of a possible refinement of the original 
regular systems corresponding to the functions $f$ and 
$h$. Then we can write, by exploiting (\ref{5.2}) 
and (\ref{5.10}) that
\begin{eqnarray*}
&&EI_n(f|j_1,\dots,j_n)I_n(h|j'_1,\dots,j'_n) 
=EI_n(f|j_1,\dots,j_n)\overline{I_n(h|j'_1,\dots,j'_n)} \\
&&\qquad =\sum_{\pi\in\Pi(n)} 
\sum_{\substack{(k_1,\dots k_n),\; (l_1,\dots l_n)\\
k_p=\pm1,\dots,\pm N,\; p=1,\dots, n\\
l_p=k_{\pi(p)} \;p=1,\dots,n}}
f(u_{k_1},\dots u_{k_n})
\overline{h(u_{k_{\pi(1)}},\dots,u_{k_{\pi(n) }}) } \\
&&\qquad\qquad\qquad\quad 
\times E Z_{G,{j_1}}(\Delta_{k_1})\cdots Z_{G,{j_n}}(\Delta_{k_n})
\overline{Z_{G,{j'_1}}(\Delta_{l_1})\cdots Z_{G,{j'_n}}(\Delta_{l_n})},
\end{eqnarray*}
where $u_k\in\Delta_k$ for all $k=\pm1,\dots,\pm N$.

The expected value of the product at the right-hand side of 
this identity can be calculated with the help of~(\ref{5.12}), 
and this yields that
\begin{eqnarray*}
&&EI_n(f|j_1,\dots,j_n)I_n(h|j'_1,\dots,j'_n) \\
&&\qquad =\sum_{\pi\in\Pi(n)} 
\sum_{\substack{(k_1,\dots k_n),\; (l_1,\dots l_n)\\
k_p=\pm1,\dots,\pm N,\; p=1,\dots, n\\
l_p=k_{\pi(p)}, \;p=1,\dots,n}}
f(u_{k_1},\dots u_{k_n})
\overline{h(u_{l_1},\dots,u_{l_n}}) \\
&&\qquad\qquad\qquad\quad 
\times G_{j_1,j'_{\pi^{-1}(1)}}(\Delta_{k_1})\cdots G_{j_n,j'_{\pi^{-1}(n)}}(\Delta_{k_n}) \\
&&\qquad =\sum_{\pi\in\Pi(n)} \int f(x_1,\dots x_n)
\overline{h(x_{\pi(1)},\dots,x_{\pi(n)})} \\
&&\qquad\qquad\qquad\times G_{j_1,j'_{\pi^{-1}(1)}}(\,dx_1)
\dots G_{j_n,j'_{\pi^{-1}(n)}}(\,dx_n).
\end{eqnarray*}
Formula~(\ref{5.4}) is proved.

The proof of~(\ref{5.5}) is based on a similar idea, but it 
is considerably simpler. It can be proved similarly 
to relation~(\ref{5.10}) that for $n\neq n'$,
\begin{equation}
EZ_{G,{j_1}}(\Delta_{k_1})\cdots Z_{G,{j_n}}(\Delta_{k_n})
\overline{Z_{G,{j'_1}}(\Delta_{l_1})\cdots 
Z_{G,{j'_{n'}}}(\Delta_{l_{n'}})}=0 \label{5.13}
\end{equation}
if we define this expression by means a regular system
$$
{\cal D}=\{\Delta_k,\;k=\pm1,\pm2,\dots,\pm N\},
$$ 
some numbers
$j_1,\dots,j_n$ and $j_1'\dots,j'_{n'}$, all of them between 
1 and $d$, together with two sequences of numbers 
$k_1,\dots,k_n$ and $l_1,\dots,l_{n'}$ such that
$k_s,l_s\in\{\pm1,\dots,\pm N\}$ for all these numbers, and they
satisfy the relation $k_s\neq\pm k_{s'}$, and $l_s\neq\pm l_{s'}$ if 
$s\neq s'$. Then, if we express
$$
EI_n(f|j_1,\dots,j_n)I_{n'}(h|j'_1,\dots,j'_{n'}) 
=EI_n(f|j_1,\dots,j_n)\overline{I_{n'}(h|j'_1,\dots,j'_{n'})} 
$$
similarly as we have done in the proof of~(\ref{5.12}) we get 
such a sum where all terms equal zero because of~(\ref{5.13}). 
This implies relation~(\ref{5.5}).

To define the Wiener--It\^o integral for all functions 
$f\in{\cal K}_{n,j_1,\dots,j_n}$ we still need the following result.

\medskip\noindent
{\bf Lemma 5.1.} {\it The class of simple functions
$\hat{\cal K}_{n,j_1,\dots,j_n}$ is a dense linear subspace
of the (real) Hilbert space ${\cal K}_{n,j_1,\dots,j_n}$.}

\medskip
Lemma 5.1 is the multivariate version of Lemma~4.1 in~\cite{9}. 
(A more transparent proof of this result was given in the 
Appendix of~\cite{10}.) Actually, we do not have to prove Lemma~5.1, 
because it simply follows from Lemma~4.1 of~\cite{9}. By applying 
this result for $G=\sum_{j=1}^n G_{j,j}$ we get that all bounded 
functions of ${\cal K}_{n,j_1,\dots,j_n}$ are in the closure 
of $\hat{\cal K}_{n,j_1,\dots,j_n}$. But this implies that all 
functions of ${\cal K}_{n,j_1,\dots,j_n}$ are in this closure.

\medskip
Let us take the $L_2$ norm in the Hilbert space ${\cal H}$.
Then we have, for all $f\in\hat{\cal K}_{n,j_1,\dots,j_n}$,
$I_n(f|j_1,\dots,j_n)\in{\cal H}$, and by formula~(\ref{5.6}),
$$
\|I_n(f|j_1,\dots,j_n)\|=\left[E(I_n(f|j_1,\dots,j_n)^2)\right]^{1/2}
\le \sqrt{n!}\|f_{n,j_1,\dots,j_n}\|.
$$
Hence Lemma~5.1 enables us to extend the Wiener--It\^o 
integral $I_n(f|j_1,\dots,j_n)$ for all
$f\in{\cal K}_{n,j_1,\dots,j_n}$. Moreover, relations 
(\ref{5.2})---(\ref{5.6}) remain valid in the Hilbert space 
${\cal K}_{n,j_1,\dots,j_n}$ after this extension.

\medskip\noindent
{\it Remark.} In~(\ref{5.6}) we have given an upper bound for 
the second moment of a multiple Wiener--It\^o integral, but we
cannot write equality in this formula. In the scalar-valued 
case we had an identity in the corresponding relation. At 
least this was the case if we took the Wiener--It\^o integral 
of a symmetric function. On the other hand, working only 
with Wiener--It\^o  integrals of symmetric functions did 
not mean a serious restriction. This relative weakness of 
formula~(\ref{5.6}) (the lack of identity) is the reason 
why we cannot represent such a large class of random 
variables in the form of a sum of Wiener--It\^o integrals 
as in the scalar valued case. (This problem will be 
discussed in Section~2 of \cite{11}.)

\medskip
I would mention that there is a slightly stronger version of
Lemma 5.1 which is useful in the study in the second part of
this paper, in~\cite{11}, when we are interested in the 
question under what conditions we can state that a sequence 
of Wiener--It\^o integrals converges to a Wiener--It\^o 
integral. Here is this result.

\medskip\noindent
{\bf Lemma 5.2.} {\it For all functions 
$f\in{\cal K}_{n,j_1,\dots,j_n}$ and numbers $\varepsilon>0$ there 
is such a simple function $g\in\hat{\cal K}_{n,j_1,\dots,j_n}$ for
which $\|f-g\|\le\varepsilon$ in the norm of the Hilbert space 
${\cal K}_{n,j_1,\dots,j_n}$, and there is a regular system 
${\cal D}=\{\Delta_k,\;k=\pm1,\pm2,\dots,\pm N\}$ to which the
function $g$ is adapted, and the boundary of all sets
$\Delta_k\in{\cal D}$ has zero $\mu$-probability with the measure
$\mu$ we chose as the dominating measure for the complex measures
$G_{j,j'}$ in our considerations.}

\medskip
Lemma~5.2 also follows from the results of~\cite{9} or~\cite{10}.

\medskip
Finally, I make the following remark. If we define 
a new function by reindexing the variables of a function of 
$h\in{\cal K}_{n,j_1,\dots,j_n}$ by means of a permutation of the 
indices, and we change the indices of the spectral measure 
$Z_{G,{j_s}}$ in the Wiener-It\^o integral $I_n(h|j_1,\dots,j_n)$
in an appropriate way, then we get a new Wiener--It\^o 
integral whose value agrees with the original integral 
$I_n(h|j_1,\dots,j_n)$. More explicitly, the following result
holds.

\medskip\noindent
{\bf Lemma 5.3.} {\it Given a function $h\in{\cal K}_{n,j_1,\dots,j_n}$  
and a permutation $\pi\in\Pi(n)$ define the function
$h_\pi(x_1,\dots,x_n)=h(x_{\pi(1)},\dots,x_{\pi(n)})$. The following
identity holds. 
\begin{eqnarray}
&&\int h(x_1,\dots,x_n)Z_{G,{j_1}}(\,dx_1)\dots Z_{G,{j_n}}(\,dx_n) \nonumber \\
&&\qquad 
=\int h_\pi(x_1,\dots,x_n)
Z_{G,{j_{\pi(1) }}}(\,dx_1)\dots Z_{G,{j_{\pi(n)}}}(\,dx_n).
\label{5.14}
\end{eqnarray}
(In particular, $h_\pi\in{\cal K}_{n,j_{\pi(1)},\dots,j_{\pi(n)}}$, thus the 
integrals on both sides of the identity are meaningful.)}

\medskip\noindent
{\it Proof of Lemma 5.3.} This identity can be simply checked if 
$h$ is a simple function. It is enough to observe that if
$h(x_1,\dots,x_n)=h_1(x_1)\cdots h_n(x_n)$ with some 
$x_l\in\Delta_{k_l}$, $g(_l(\cdot)$ is some function on 
${\mathbb R}^{\nu}$, $1\le  l\le n$, then
$$
\int h(x_1,\dots,x_n)Z_{G,{j_1}}(\,dx_1)\dots Z_{G,{j_n}}(\,dx_n) 
=\prod_{l=1}^n h_l(x_l)Z_{G,j_l}(\Delta_{k_l}),
$$
$h_\pi(x_1,\dots,x_l)=h_1(x_{\pi_1})\cdots h_n(x_{\pi_n})$, 
$$
\int h_\pi(x_1,\dots,x_n)
Z_{G,{j_{\pi(1) }}}(\,dx_1)\dots Z_{G,{j_{\pi(n)}}}(\,dx_n)=
\prod_{l=1}^n h(x_{\pi_l})Z_{G,j_{\pi_l}}(\Delta_{k_{\pi(l)}}),
$$
and the last two Wiener--It\^o integrals equal. Then a 
simple limiting procedure implies it in the general case. 
Lemma~5.3 is proved.

\medskip
We saw in~\cite{9} that in the scalar valued case the value 
of a Wiener--It\^o integral 
$\int f(x_1,\dots,x_n)Z_G(\,dx_1)\dots Z_G(\,dx_n)$ 
does not change if we replace the kernel function~$f$ by 
the function we get by permuting its variables $x_1,\dots,x_n$ 
in an arbitrary way. Lemma~5.3 is the generalization of this 
result to the case when we integrate with respect to the 
coordinates of a vector valued random spectral measure. 

\medskip\noindent
{\it Remark.} A consequence of the result of Lemma~5.3 shows
an essential difference between the behaviour of multiple
Wiener--It\^o integrals with respect to scalar and vector 
valued random spectral measures. It follows from the scalar
valued version of Lemma~5.3 that in the scalar valued case
the Wiener--It\^o integral of a kernel function agrees with
the Wiener--it\^o integral of the symmetrization of this
kernel function. This has the consequence that in the scalar
valued case we can restrict our attention to the Wiener--It\^o
integrals of symmetrical functions which do not change their
values by any permutation of their variables. It can be seen
that any random variable which can be written as the sum 
of Wiener--It\^o integrals can be written in a unique form
as a sum of Wiener--It\^o integrals of different multiplicity
with symmetric kernel functions. The analogous result does 
not hold in the vector valued case. Indeed, if there is 
some linear dependence among the coordinates of the
underlying vectors in a vector valued stationary random
field, then such functions $f_j$ can be found for
which $\sum_{j=1}^d \int f_j(x)Z_{G,j}(\,dx)\equiv0$,
and not all kernel functions $f_j$ disappear in the above
sum. This shows that the unique representation of the
random variables by means of a sum of Wiener--It\^o 
integrals may not hold in vector valued models.

\section{The diagram formula for the product of multiple 
Wiener--It\^o integrals}

Let us consider a vector valued random spectral measure 
$(Z_{G,1,}\dots,Z_{G,d})$ corresponding to the matrix valued 
spectral measure $(G_{j,j'})$, $1\le j,j'\le d$, of a vector 
valued stationary Gaussian random field with expectation zero 
(either to a discrete random field 
$X(p)=(X_1(p),\dots,X_d(p))$, $p\in{\mathbb Z}^\nu$, or to a 
generalized one $X(\varphi)=(X_1(\varphi),\dots,X_d(\varphi))$, 
$\varphi\in{\cal S}^\nu$). Let us assume that the spectral 
measure $G_{j,j'}$, $1\le j,j'\le d$, is non-atomic, and take 
two Wiener--It\^o integrals 
\begin{equation}
I_n(h_1|j_1,\dots,j_n)=\int h_1(x_1,\dots,x_n)Z_{G,{j_1}}(\,dx_1)\dots 
Z_{G,{j_n}}(dx_n) \label{6.1}
\end{equation}
and
\begin{equation}
I_m(h_2|j'_1,\dots,j'_m)=\int h_2(x_1,\dots,x_m)Z_{G,{j'_1}}(\,dx_1)\dots 
Z_{G,{j'_m}}(dx_m) \label{6.2}
\end{equation}
with some kernel functions $h_1\in{\cal K}_{n,j_1,\dots,j_n}$ and 
$h_2\in{\cal K}_{m,j'_1,\dots,j'_m}$, where $j_s, j'_t\in\{1,\dots,d\}$
for all $1\le s\le n$ and $1\le t\le m$.

Actually we state our problems a bit differently, 
which is more appropriate for our discussion.\!
 We take two 
functions $h_1(x_1,\dots,x_n)$ and $h_2(x_{n+1},\dots,x_{n+m})$ 
in the space ${\mathbb R}^{(n+m)\nu}$, and define the function 
\hfill\break 
$h^{(0)}_2(x_1,\dots,x_m)$ by the identity
$$
h^{(0)}_2(x_1,\dots,x_m)=h_2(x'_{n+1},\dots,x'_{n+m}))
\textrm{ if } (x_1,\dots,x_m)=(x'_{n+1},\dots,x'_{n+m}).
$$
We assume that $h_1\in{\cal K}_{n,j_1,\dots,j_n}$,
$h^{(0)}_2\in{\cal K}_{m,j'_1,\dots,j'_m}$. Then we define the 
Wiener--It\^o integrals~(\ref{6.1}) and~(\ref{6.2}) with 
the kernel functions $h_1$ and $h_2^{(0)}$. In 
formula~(\ref{6.2}) we  should have written the function 
$h_2^{(0)}$, but we omitted the superscript~$^{(0)}$.

I shall present a result in which we express the product of these 
two Wiener--It\^o integrals as a sum of Wiener--It\^o integrals.
This result is called the diagram formula, since the kernel 
functions of the Wiener--It\^o integrals appearing in this sum 
are expressed by means of some diagrams. This result is a 
multivariate version of the diagram formula proved in Chapter~5 
of \cite{9}. In that work also the product of more than two 
Wiener--It\^o integrals is expressed in the form of a sum of 
Wiener--It\^o integrals. But actually the main point of the proof 
is to show the validity of the diagram formula for the product of 
two Wiener--It\^o integrals, and we shall need only this result. 
So I restrict my attention to this case. Actually we need 
the diagram formula only in a special case. The result in this 
special case will be given in a corollary.

To express the product of the two Wiener--It\^o integrals 
in formulas~(\ref{6.1}) and~(\ref{6.2}) as a sum of 
Wiener--It\^o integrals first I introduce a class of coloured 
diagrams $\Gamma=\Gamma(n,m)$ that will be used in the 
definition of the Wiener--It\^o integrals we shall be working 
with. A coloured diagram $\gamma\in\Gamma$ is a graph whose 
vertices are the pairs of integers $(1,s)$, $1\le s\le n$, and 
$(2,t)$, $1\le t\le m$. Each vertex is coloured with one of 
the numbers $1,\dots,d$. The colour of the vertex $(1,s)$ is 
$j_s$, $1\le s\le n$, and the colour of the vertex $(2,t)$ is 
$j'_t$, $1\le t\le m$. The set of vertices of the form $(1,s)$ 
will be called the first row and the set of vertices of the 
form $(2,t)$ will be called the second row  of a diagram 
$\gamma\in\Gamma$. The coloured diagrams $\gamma\in\Gamma$ are 
those undirected graphs with the above coloured vertices for 
which edges can go only between vertices of the first and 
second row, and from each vertex there starts zero or one edge. 
Given a coloured diagram $\gamma\in\Gamma$ we shall denote 
the number of its edges by $|\gamma|$.

I shall define for all coloured diagrams $\gamma\in\Gamma$ a
multiple Wiener--It\^o integral depending on~$\gamma$. The
diagram formula states that the product of the Wiener--It\^o 
integrals in~(\ref{6.1}) and~(\ref{6.2}) equals the sum of 
these Wiener--It\^o integrals. 

When stating the diagram formula I shall work
with the functions \hfill\break 
$h_1(x_1,\dots,x_n)$ and $h_2(x_{n+1},\dots,x_{n+m})$ in 
${\mathbb R}^{n+m}$. The function  \hfill\break
$h_2(x_{n+1},\dots,x_{n+m})$ 
is the function which corresponds to the kernel function 
$h^{(0)}_2(x_1,\dots,x_m)$ in the definition of the
Wiener--It\^o integral in~(\ref{6.2}). We define
with their help the function
\begin{equation}
H(x_1,\dots,x_{n+m})=h_1(x_1,\dots,x_n)h_2(x_{n+1},\dots,x_{n+m}). 
\label{6.3}
\end{equation}
We shall define the kernel functions appearing in the 
Wiener--it\^o integrals in the diagram formula with the
help of the functions $H(x_1,\dots,x_{n+m})$. In the definition
of these kernel functions I shall apply the following 
natural bijection~$S$ between the coordinates of the vectors 
in~${\mathbb R}^{n+m}$, i.e., the set $\{1,\dots,n+m\}$ and
the vertices of the diagrams of $\gamma\in\Gamma$. 
\begin{equation}
S((1,k))=k \textrm{ for } 1\le k\le n, \quad\textrm{and} \quad
S((2,k))=n+k \textrm{ for } 1\le k\le m. \label{6.4}
\end{equation} 
To simplify the formulation of the diagram formula 
I shall introduce the following notation with the
help of the colours of the diagrams.
\begin{equation}
J(1,k)=j_k,\;\; 1\le k\le n \quad\textrm{and}\quad
J(2,l)=j'_l, \;\; 1\le l\le m. \label{6.5}
\end{equation}

First I give the formal definition of the Wiener--It\^o 
integrals that appear in the diagram formula. These 
Wiener-It\^o integrals correspond to the diagrams
$\gamma\in\Gamma$ introduced before. Then I describe 
the diagram formula with the help of these Wiener--It\^o 
integrals. The definition of the Wiener--It\^o integrals 
we need in the diagram formula applies a rather 
complicated notation, but its informal explanation given 
after formula~(\ref{6.16}) may help to understand it. For 
the sake of a better comprehension of the calculations in 
the diagram formula I shall present an example after the 
formulation of this result, where the product of two 
Wiener--It\^o integrals is considered, and I show how to 
calculate a typical term in the sum of Wiener--It\^o 
integrals which appears in the diagram formula for this 
product.

Fix some diagram $\gamma\in\Gamma$. I explain how 
to define the Wiener--It\^o integral corresponding to 
$\gamma$ in the diagram formula. First I define a function 
$H_\gamma(x_1,\dots,x_{n+m})$ which we get by means of an
appropriate permutation of the indices of the function~$H$
defined in (\ref{6.3}). This permutation of the indices
depends on the diagram~$\gamma$.

To define this permutation of the indices first
I define a map $T_\gamma$ which maps the set $\{1,\dots,n+m\}$ 
to the elements in the rows of the diagrams. This map depends 
on the diagram~$\gamma$. 

To define this map first I introduce the 
following sets depending on the diagram $\gamma$:
\begin{eqnarray}
A_1=A_1(\gamma)&=&\{r_1,\dots,r_{n-|\gamma|}\colon\; 1\le r_1<r_2<\cdots
<r_{n-|\gamma|}\le n   \label{6.6} \\ 
&&\qquad\qquad\textrm{no edge of } \gamma \textrm{ starts from }
(1, r_k),\quad 1\le k\le n-|\gamma|\}, \nonumber 
\end{eqnarray}
\begin{eqnarray}
A_2=A_2(\gamma)&=&\{t_1,\dots,t_{m-|\gamma|}\colon\; 1\le t_1<t_2<\cdots
<t_{m-|\gamma|}\le m,  \label{6.7} \\ 
&&\qquad\qquad\textrm{no edge of } \gamma \textrm{  starts from }(2,t_k),
\quad 1\le k\le m-|\gamma|\},  \nonumber
\end{eqnarray}
and
\begin{eqnarray}
B=B(\gamma)&=&\{(v_1,w_1),\dots, (v_{|\gamma|},w_{|\gamma|}))\colon\; 
1\le v_1<v_2<\cdots v_{|\gamma|}\le n,  \nonumber \\
&&\qquad ((1,v_k),(2,w_k)) \textrm{ is an edge of } |\gamma|, 
\;\; 1\le k\le|\gamma|\}.
\label{6.8}
\end{eqnarray}
Let us also define with the help of the set $B$ the sets
\begin{equation}
B_1=B_1(\gamma)=\{v_1,\dots,v_{|\gamma|}\}, \quad B_2=B_2(\gamma)
=\{w_1,\dots,w_{|\gamma|}\} \label{6.9}
\end{equation}
with the numbers $v_k$ and $w_l$ appearing in the set
$$
B=B(\gamma)=\{(v_1,w_1)),\dots, (v_{|\gamma|},w_{|\gamma|}))\}.
$$
Now, I define the map  $T_\gamma$ in the following way:
\begin{eqnarray}
&&T_\gamma(k)=(1,r_k) \;\;\textrm{for} \;\;1\le k\le n-|\gamma|, 
\label{6.10} \\
&&T_\gamma(n-|\gamma|+k)=(2,t_k) \;\;\textrm{for} \;\;
1\le k\le m-|\gamma|, \nonumber  \\ 
&&T_\gamma(n+m-2|\gamma|+k)=(1,v_k) \;\;\textrm{for}\;\; 
1\le k \le |\gamma|,  \nonumber \\ 
&&T_\gamma(n+m-|\gamma|+k)=(2,w_k) \;\;\textrm{for}\;\; 
1\le k \le|\gamma|. \nonumber
\end{eqnarray}
In formula~(\ref{6.10}) we worked with the numbers $r_k$, $t_k$,
$v_k$ and $w_k$ defined in (\ref{6.6})---(\ref{6.9}). It 
has the following meaning. We listed the vertices of the diagram 
$\gamma$ in the form $T_\gamma(s)$, $1\le s\le n+m$. If the vertex
$T_\gamma(s)$ gets the index $s$, then the first  $n-|\gamma|$ 
indices are given in increasing order to the vertices from the 
first row from which no edge starts. The vertices of the
second row from which no edge starts get the next $m-|\gamma|$ 
indices also in increasing order. Then the $|\gamma|$ vertices 
from the first row from which an edge starts get the subsequent 
$|\gamma|$ indices in increasing order. The remaining $|\gamma|$ 
vertices from the second row from which an edge starts get the 
indices between $n+m-|\gamma|+1$ and $n+m$. They are indexed 
in such a way that if two vertices $(1,v_k)$ and $(2,w_k)$ are 
connected by en edge then the index of $(2,w_k)$ is obtained if 
we add $|\gamma|$ to the index of $(1,v_k)$. 

I define with the help of the function $T_\gamma$ and the
map  $S(\cdot)$ defined in (\ref{6.4}) the permutation
\begin{equation}
\pi_\gamma(k)=S(T_\gamma(k)), \quad 1\le k\le n+m \label{6.11}
\end{equation}
of the set $\{1,\dots,n+m\}$. Next I introduce the Euclidean
space ${\mathbb R}^{n+m}_\gamma$ with elements 
$x(\gamma)=(x(\gamma)_1,\dots,x(\gamma)_{n+m})$ by reindexing 
the arguments of the  Euclidean space  ${\mathbb R}^{n+m}$, 
where the functions $h_1(x_1,\dots,x_n)$ and 
$h_2(x_{n+1},\dots,x_{n+m})$ are defined in the following way.
$$
(x(\gamma)_1,\dots,x(\gamma)_{n+m})
=(x_{\pi_\gamma(1)},\dots,x_{\pi_\gamma(n+m)}) 
$$
with
$(x(\gamma)_1,\dots,x(\gamma)_{n+m})\in{\mathbb R}^{n+m}_\gamma$ 
and $(x_1,\dots,x_{n+m})\in{\mathbb R}^{n+m}$. It will be simpler 
to define the quantities needed in the definition of the 
Wiener--It\^o integral corresponding to the diagram~$\gamma$ 
as functions defined in the space $R^{n+n}_\gamma$. First we 
define the function $H_\gamma$ as
\begin{eqnarray}
&&
\!\!\!\!\!\!\!\!\!\!\!
H_\gamma(x(\gamma)_1,\dots,x(\gamma)_{n+m}) \label{6.12} \\ 
&&
\!\!\!\!\!\!\!\!\!\!\!
\qquad=H(x(\gamma)_1,\dots,x(\gamma)_{n-|\gamma|},
x(\gamma)_{n+m-2|\gamma|+1},\dots,x(\gamma)_{n+m-|\gamma|}, 
\nonumber\\
&&
\!\!\!\!\!\!\!\!\!\!\!
\qquad\qquad\qquad x(\gamma)_{n-|\gamma|+1},\dots,x(\gamma)_{n+m-2|\gamma|+1},
x(\gamma)_{(n+m-|\gamma|+1},\dots, x(\gamma)_{n+m}) \nonumber \\
&&
\!\!\!\!\!\!\!\!\!\!\!
\qquad=h_1(x(\gamma)_1,\dots,x(\gamma)_{n-|\gamma|},
x(\gamma)_{\pi_\gamma(n+m-2|\gamma|+1)},\dots,
x(\gamma)_{n+m-|\gamma|}) \nonumber \\
&&
\!\!\!\!\!\!\!\!\!\!\!
\qquad\qquad \times h_2(x(\gamma)_{n-|\gamma|+1},\dots,x(\gamma)_{n+m-2|\gamma|+1},
x(\gamma)_{n+m-|\gamma|+1},\dots, x(\gamma)_{n+m}).
\nonumber 
\end{eqnarray}

Next I define the function 
$\bar h_\gamma(x(\gamma)_1,\dots,x(\gamma)_{n+m-|\gamma|}))$ (with 
$n+m-|\gamma|$ arguments) which we get by replacing 
$x(\gamma)_{n+m-|\gamma|+k}$ by               
$-x(\gamma)_{n+m-2|\gamma|+k})$ in the 
function $H_\gamma$ defined in formula~(\ref{6.12}) for all 
$1\le k\le \gamma$, i.e., I define
\begin{eqnarray}
&&
\!\!\!\!\!\!\!\!\!\!\!
\bar h_\gamma(x(\gamma)_1,\dots,x(\gamma)_{n+m-|\gamma|})  \label{6.13} \\
&&
\!\!\!\!\!\!\!\!\!\!\!
\qquad=H_\gamma(x(\gamma)_1,\dots,x(\gamma)_{n+m-|\gamma|},
-x(\gamma)_{n+m-2|\gamma|+1},\dots,-x(\gamma)_{n+m-|\gamma|})  \nonumber \\ 
&&
\!\!\!\!\!\!\!\!\!\!\!
\qquad=H(x(\gamma)_1,\dots,x(\gamma)_{n-|\gamma|},
x(\gamma)_{n+m-2|\gamma|+1},\dots,x(\gamma)_{n+m-|\gamma|}, 
\nonumber\\
&&
\!\!\!\!\!\!\!\!\!\!\!
\qquad\qquad\qquad x(\gamma)_{n-|\gamma|+1},\dots,x(\gamma)_{n+m-2|\gamma|+1},
\nonumber \\
&&
\!\!\!\!\!\!\!\!\!\!\!
\qquad\qquad\qquad\qquad 
-x(\gamma)_{n+m-2|\gamma|+1},\dots, -x(\gamma)_{n+m-|\gamma|}) 
\nonumber \\
&&
\!\!\!\!\!\!\!\!\!\!\!
\qquad=h_1(x(\gamma)_1,\dots,x(\gamma)_{n-|\gamma|},
x(\gamma)_{n+m-2|\gamma|+1},\dots,
x(\gamma)_{n+m-|\gamma|}) \nonumber \\
&&
\!\!\!\!\!\!\!\!\!\!\!
\qquad\qquad \times h_2(x(\gamma)_{n-|\gamma|+1},\dots,x(\gamma)_{n+m-2|\gamma|+1},
\nonumber \\
&&
\!\!\!\!\!\!\!\!\!\!\!
\qquad\qquad\qquad\qquad
-x(\gamma)_{n+m-2|\gamma|+1},\dots, -x(\gamma)_{n+m)-|\gamma|}).
\nonumber 
\end{eqnarray}
 In the next step I define the function
$\bar{\bar h}_\gamma(x(\gamma)_1,\dots,x(\gamma)_{n+m-2|\gamma|})$. 
This will be the kernel
function of the Wiener--It\^o integral which corresponds to
the diagram~$\gamma$ in the diagram formula if we
express it as a Wiener--It\^o integral with respect to the
variables $x(\gamma)_1,\dots,x(\gamma)_{n+m-2|\gamma|}$, 
\begin{eqnarray}
&&\bar{\bar h}_\gamma(x_\gamma)_1,\dots,x(\gamma)_{n+m-2|\gamma|}) 
=\int \bar h_\gamma(x(\gamma)_1,\dots,x(\gamma)_{n+m-|\gamma|}) \label{6.14} \\
&&\qquad\qquad \times \prod_{k=1}^{|\gamma|} 
G_{J(S^{-1}(n+m-2|\gamma|+k)),J(S^{-1}(n+m-|\gamma|+k))}
(\,dx(\gamma)_{n+m-2|\gamma|+k}) \nonumber \\
&&\qquad=\int \bar h_\gamma(x(\gamma)_1,\dots,x(\gamma)_{n+m-|\gamma|}) 
\prod_{k=1}^{|\gamma|} G_{j_{v_k},j'_{w_k}}(\,dx(\gamma)_{n+m-2|\gamma|+k}) \nonumber 
\end{eqnarray}
with the function $J(\cdot)$ defined in~(\ref{6.5}), the indices
$v_k$ and $w_k$ defined in~(\ref{6.8}) and the
function $T_\gamma$ defined in~(\ref{6.10}). 

I shall show that the Wiener--It\^o integrals
\begin{eqnarray}
&&I_{n+m-2|\gamma|}(\bar{\bar h}_\gamma|j_{r_1},\dots,j_{r_{n-|\gamma|}},
j'_{t_1},\dots,j'_{t_{m-|\gamma|}}) \label{6.15} \\
&&\qquad=\int \bar{\bar h}_\gamma(x(\gamma)_1,\dots,x(\gamma)_{n+m-2|\gamma|})  
\prod_{k=1}^{n+m-2|\gamma|}Z_{G,{J(S^{-1}(k))}}(dx(\gamma)_k)  \nonumber \\
&&\qquad=\int\bar{\bar h}_\gamma(x(\gamma)_1,\dots,x(\gamma)_{n+m-2|\gamma|}) 
\nonumber\\ 
&& \qquad\qquad\qquad \prod_{k=1}^{n-|\gamma|}Z_{G,{j_{r_k}}}(\,dx(\gamma)_k)
\prod_{l=1}^{m-|\gamma|}Z_{G,{j'_{t_l}}}(\,dx(\gamma)_{l+n-|\gamma|})
\nonumber 
\end{eqnarray}
exist for all $\gamma\in\Gamma$, and these Wiener--It\^o integrals
appear in the diagram formula. The numbers $r_k$ and $t_l$ in
this formula were defined in~(\ref{6.6}) and~(\ref{6.7}).

In formula~(\ref{6.15}) we integrated with respect to the coordinates
$x(\gamma)_s$, $1\le s\le n+m$, of the vectors in the Euclidean space 
${\mathbb R}^{n+m}_\gamma$. If we replace the variables $x(\gamma)_s$
by $x_s$ in~(\ref{6.15}), then we get a Wiener--it\^o integral in
the space ${\mathbb R}^{n+m}$  with the same value. This means
that the following relation holds:
\begin{eqnarray}
&&I_{n+m-2|\gamma|}(\bar{\bar h}_\gamma|j_{r_1},\dots,j_{r_{n-|\gamma|}},
j'_{t_1},\dots,j'_{t_{m-|\gamma|}}) \label{6.16} \\
&&\qquad=I_{n+m-2|\gamma|}(h_\gamma|j_{r_1},\dots,j_{r_{n-|\gamma|}},
j'_{t_1},\dots,j'_{t_{m-|\gamma|}}) \nonumber \\
&&\qquad=\int h_\gamma(x_1,\dots,x_{n+m-2|\gamma|}) 
\nonumber\\ 
&& \qquad\qquad\qquad \prod_{k=1}^{n-|\gamma|}Z_{G,{j_{r_k}}}(\,dx_k)
\prod_{l=1}^{m-|\gamma|}Z_{G,{j'_{t_l}}}(\,dx_{l+n-|\gamma|})
\nonumber 
\end{eqnarray}
with
\begin{eqnarray*}
h_\gamma(x_1,\dots,x_{n+m-2|\gamma|})
&=&\bar{\bar h}_\gamma(x(\gamma)_1,\dots,x(\gamma)_{n+m-2|\gamma|}) \\
&=&\bar{\bar h}_\gamma(x_{\pi_\gamma(1)},\dots,x_{\pi_\gamma(n+m-2|\gamma|)}).
\end{eqnarray*}

Before describing the diagram formula I explain the content
of the above defined formulas.

Let us fix a diagram $\gamma\in\Gamma$, and let us call a 
vertex of it  from which no edge starts open, and a vertex 
from which an edge starts closed. We listed the open vertices 
from the first row in increasing order as 
$(1,r_1),\dots,(1,r_{n-|\gamma|})$, and the open vertices from 
the second row as $(2,t_1),\dots,(2,t_{m-|\gamma|})$. We 
listed the closed vertices from the first row in increasing 
order as $(1,v_1),\dots,(1,v_\gamma)$. Finally we listed the 
closed vertices from the second row as 
$(2,w_1),\dots,(2,w_\gamma)$, and we indexed them in such a way
that the vertices $(1,v_k)$ and $(2,w_k)$ are connected by
an edge for all $1\le k\le\gamma$.  

In formula~(\ref{6.10}) we defined the map $T_\gamma$ from 
the set $\{1,\dots,n+m\}$ to the set of vertices of the 
diagram~$\gamma$ with the help of the above listing of 
the vertices. First we considered the open vertices from 
the first row, then the open vertices from the second row, 
and then we finished with the closed vertices first from 
the first and then from the second row. We defined 
in~(\ref{6.11}) the permutation $\pi_\gamma$ of the set 
$\{1,\dots,n+m\}$ by applying first the map 
the map $T_\gamma$ and then the map $S$ defined~(\ref{6.4}).
We defined the function $H_\gamma$ in (\ref{6.13}) with 
the help of this permutation. We have introduced a 
Euclidean space ${\mathbb R}^{n+m}_\gamma$ whose elements 
we get by rearranging the indices of the coordinates of 
the Euclidean space ${\mathbb R}^{n+m}$ where we are 
working with the help of the permutation~$\pi_\gamma$,
and we have defined our functions in this space.

We defined the function $H_\gamma$ on the space 
${\mathbb R}^{n+m}_\gamma$ as the product of the functions $h_1$ 
and $h_2$ with reindexed variables. In the function $h_1$ 
first we took the variables $x(\gamma)_s=x_{\pi_\gamma(s)}$
with those indices $\pi_\gamma(s)$ which correspond to the
open vertices of the first row, and then the variables with 
indices corresponding to the closed vertices of the first row.
We defined the reindexation of the variables in the second
row similarly. First we took those variables whose indices
correspond to the open vertices and then the variables whose
indices correspond to the closed vertices of the second row. 

The variables 
$$
x(\gamma)_{n+m-2|\gamma|+k}=x_{\pi_\gamma(n+m-2|\gamma|+k)} \textrm{ and }
x(\gamma)_{n+m-|\gamma|+k}=x_{\pi_\gamma(n+m-|\gamma|+k)}
$$ 
in the function $H_\gamma$ are variables with indices 
corresponding to vertices connected by an edge. So in 
the definition of the function
$\bar h_\gamma$ in~(\ref{6.14}) I replaced in $H_\gamma$ the 
variable corresponding to the endpoint of an edge from 
the second row of the diagram~$\gamma$ by the variable 
corresponding to the other endpoint of this edge, and 
multiplied this variable by~$-1$. Thus the variables
$x(\gamma)_{n+m-2|\gamma|+k}=x_{\pi_\gamma(n+m-2|\gamma|+k)}$, 
$1\le k\le|\gamma|$, of the function $\bar h_\gamma$ 
correspond to the edges of the diagram~$\gamma$. I 
defined the function $\bar{\bar h}_\gamma$ by integrating the 
function $\bar h_\gamma$ by these variables. The variable 
$x(\gamma)_{n+m-2|\gamma|+k}=x_{\pi_\gamma(n+m-2|\gamma|+k)}$ corresponds
to the $k$-th edge of the diagram, and we integrate this
variable with respect to the measure $G_{j_{v_k},j'_{w_k}}$,
that is with respect to the measure $G_{u,v}$ whose
coordinates are the colours of the endpoints of the
$k$-th edge. 

Finally we define the Wiener--It\^o integral corresponding
to the diagram $\gamma$ with kernel 
function~$\bar{\bar h}_\gamma$. We integrate the argument 
$x(\gamma)_k$ with respect to that random
spectral measure~$Z_{G,j}$ whose parameter agrees with
the colour of the vertex corresponding to this variable.
Thus we choose $Z_{G,j_{r_k}}(\,dx(\gamma)_k)$ for 
$1\le k\le n-|\gamma|$ and 
$Z_{G_{j'_{t_{k-n+|\gamma}}}}(\,dx(\gamma)_k)$ if 
$n-|\gamma|+1\le k\le n+m-2|\gamma|$. We can replace
this Wiener--It\^o integral defined in~(\ref{6.15}) 
with kernel function~$\bar{\bar h}_\gamma$ by the 
Wiener--It\^o integral defined in~(\ref{6.16}) with
kernel function~$h_\gamma$.

Next I formulate the diagram formula.

\medskip\noindent
{\bf Theorem 6.1. The diagram formula.} {\it Let us consider 
the Wiener--It\^o integrals $I_n(h_1|j_1,\dots,j_n)$ and 
$I_m(h_2|j_1',\dots,j_m')$ introduced in formulas (\ref{6.1}) and 
(\ref{6.2}). The following results hold.

\medskip
\begin{description}
\item[{\rm(A)}] The function $\bar{\bar h}_\gamma$ defined 
in~(\ref{6.14}) satisfies the relations
$$
\bar{\bar h}_\gamma\in {\cal K}_{n+m-2|\gamma|,
j_{r_1},\dots,j_{r_{n-|\gamma|}},j'_{t_1},\dots,j'_{t_{m-|\gamma|}}},
$$ 
and $\|\bar{\bar h}_\gamma\|\le \|h_1\| \|h_2\|$ for all $\gamma\in\Gamma$. 
Here the norm of  the function $h_1$ in ${\cal K}_{n,j_1,\dots,j_n}$, 
the norm of $\bar{\bar h}_2$ in ${\cal K}_{m,j'_1,\dots,j'_m}$, and 
the norm of $\bar{\bar h}_\gamma$ in 
${\cal K}_{n+m-2|\gamma|,
j_{r_1},\dots,j_{r_{n-|\gamma|}},j'_{t_1},\dots,j'_{t_{m-|\gamma|}}}$ is taken.

\item[{\rm(B)}] One has
\begin{eqnarray}
&&I_n(h_1|j_1,\dots,j_n)I_m(h_2|j'_1,\dots,j'_m) \label{6.17} \\
&&\qquad =\sum_{\gamma\in\Gamma}  
I_{n+m-2|\gamma|}(\bar{\bar h}_\gamma|j_{r_1},\dots,j_{r_{n-|\gamma|}},j'_{t_1},
\dots,j'_{t_{m-|\gamma|}}). \nonumber
\end{eqnarray}
The terms in the sum at the right-hand side of formula~(\ref{6.17}) 
were defined in formulas~(\ref{6.12})---(\ref{6.15}). The 
Wiener--It\^o integral
$$
I_{n+m-2|\gamma|}(\bar{\bar h}_\gamma|j_{r_1},\dots,j_{r_{n-|\gamma|}},j'_{t_1},
\dots,j'_{t_{m-|\gamma|}})
$$ 
in formula~(\ref{6.17}) can be replaced by the Wiener--It\^o integral
$$I_{n+m-2|\gamma|}(h_\gamma|j_{r_1},\dots,j_{r_{n-|\gamma|}},j'_{t_1},
\dots,j'_{t_{m-|\gamma|}})
$$ 
defined in~(\ref{6.16}).
\end{description}
}

\medskip
To understand the formulation of the diagram formula better 
let us consider the following example. We take a five 
dimensional stationary Gaussian random field with some 
spectral measure $(G_{j,j'}(x))$, $1\le j,j'\le 5$, and random 
spectral measure $Z_{G,j}(\,dx)$, $1\le j\le 5$, corresponding 
to it. Let us understand how we define the Wiener--It\^o 
integral corresponding to a typical diagram when we apply 
the diagram formula in the following example. Take the 
product of two Wiener--It\^o integrals of the following form:
$$
I_3(h_1|2,3,5)=\int h_1(x_1,x_2,x_3)Z_{G,2}(\,dx_1)Z_{G,3}(\,dx_2)Z_{G,5}(\,dx_3)
$$
and
\begin{eqnarray*}
I_4(h_2|1,5,4,1)&=&\int h_2(x_1,x_2,x_3,x_4) \\
&& \qquad Z_{G,1}(\,dx_1)Z_{G,5}(\,dx_2) Z_{G,4}(\,dx_3)Z_{G,2}(\,dx_4),
\end{eqnarray*}
and let us write it in the form of a sum of Wiener--It\^o 
integrals with the help of the diagram formula.

First I give the vertices of the coloured diagrams we shall be 
working with together with their colours.

\medskip

\begin{figure}[h]
\begin{center}
\epsfig{file=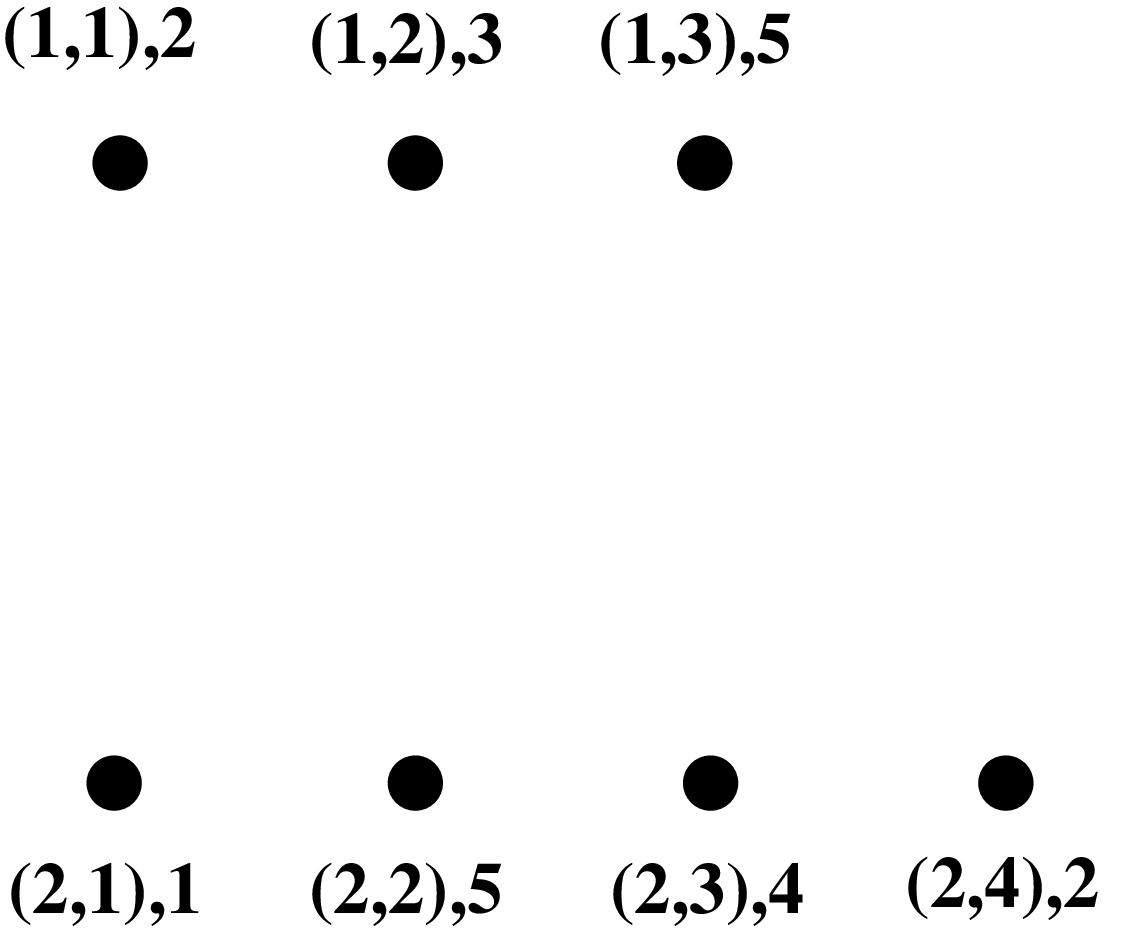, width=5cm}
\end{center}
\caption{the vertices of the diagrams together with their colours}
\end{figure}

\medskip

Next I consider a diagram~$\gamma$ which yields one of the terms 
in the sum  expressing the product of these two Wiener--It\^o 
integrals. I take the diagram which has two edges, one edge 
connecting the vertices $(1,2)$ and $(2,4)$, and another edge 
connecting the vertices $(1,3)$ and $(2,1)$. Let us calculate 
which Wiener--It\^o integral corresponds to this diagram~$\gamma$.

\medskip

\begin{figure}[h]
\begin{center}
\epsfig{file=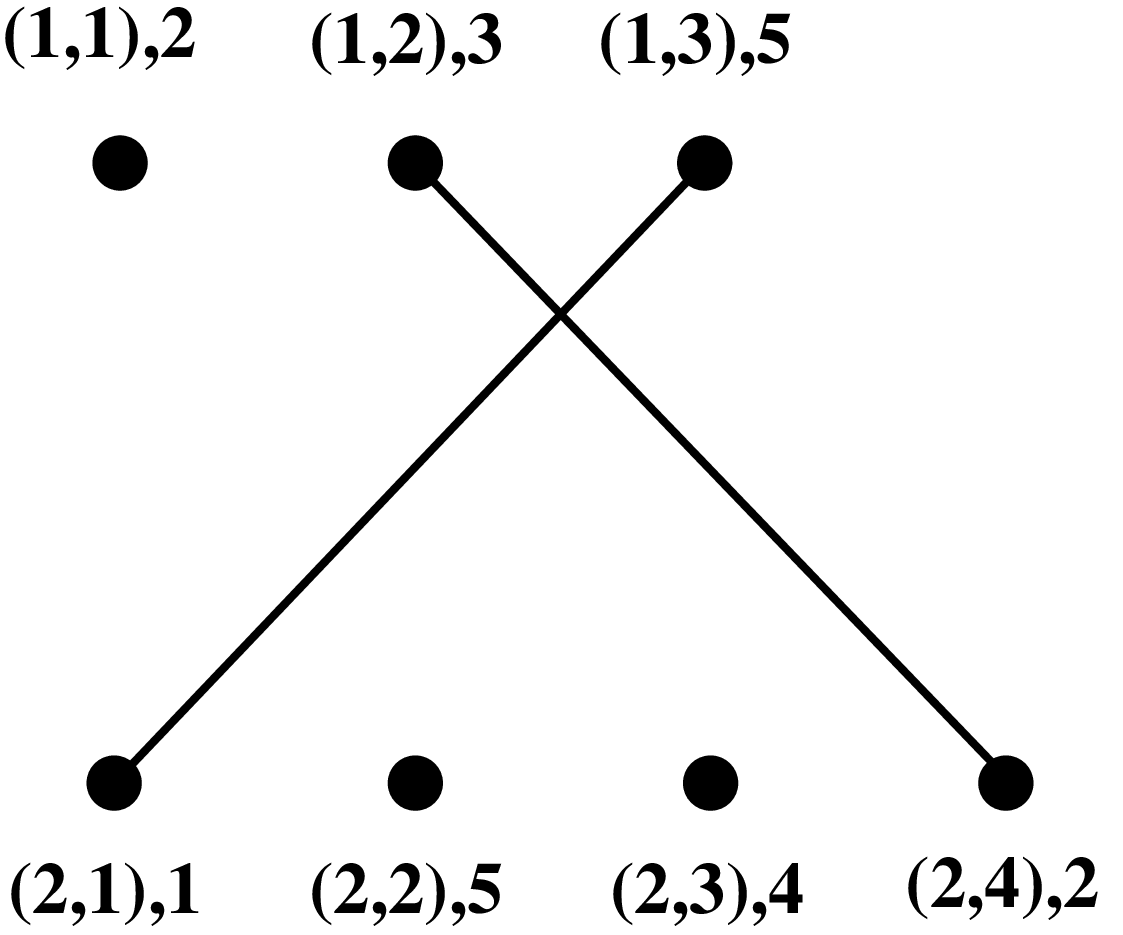, width=5cm}
\end{center}
\caption{a typical diagram}
\end{figure}

\medskip

In the next step I take this diagram~$\gamma$, and I show not 
only the indices and colours of its vertices, but for each 
vertex I also tell which value $T_\gamma(k)$ it equals. Here
$T_\gamma(k)$ is the function defined in formula~(\ref{6.10}). 

\medskip

\psfrag{(1,1),2}{$(1,1)=T_\gamma(1),2$}
\psfrag{(1,2),3}{$\;\; (1,2)=T_\gamma(4),3$}
\psfrag{(1,3),5}{$\quad \;\;(1,3)=T_\gamma(5),5$}
\psfrag{(2,1),1}{$(2,1)=T_\gamma(7),1$}
\psfrag{(2,2),5}{$(2,5)=T_\gamma(2),5$}
\psfrag{(2,3),4}{$(2,3)=T_\gamma(3),4$}
\psfrag{(2,4),2}{$(2,4)=T_\gamma(6),2$}

\begin{figure}[h]
\epsfig{file=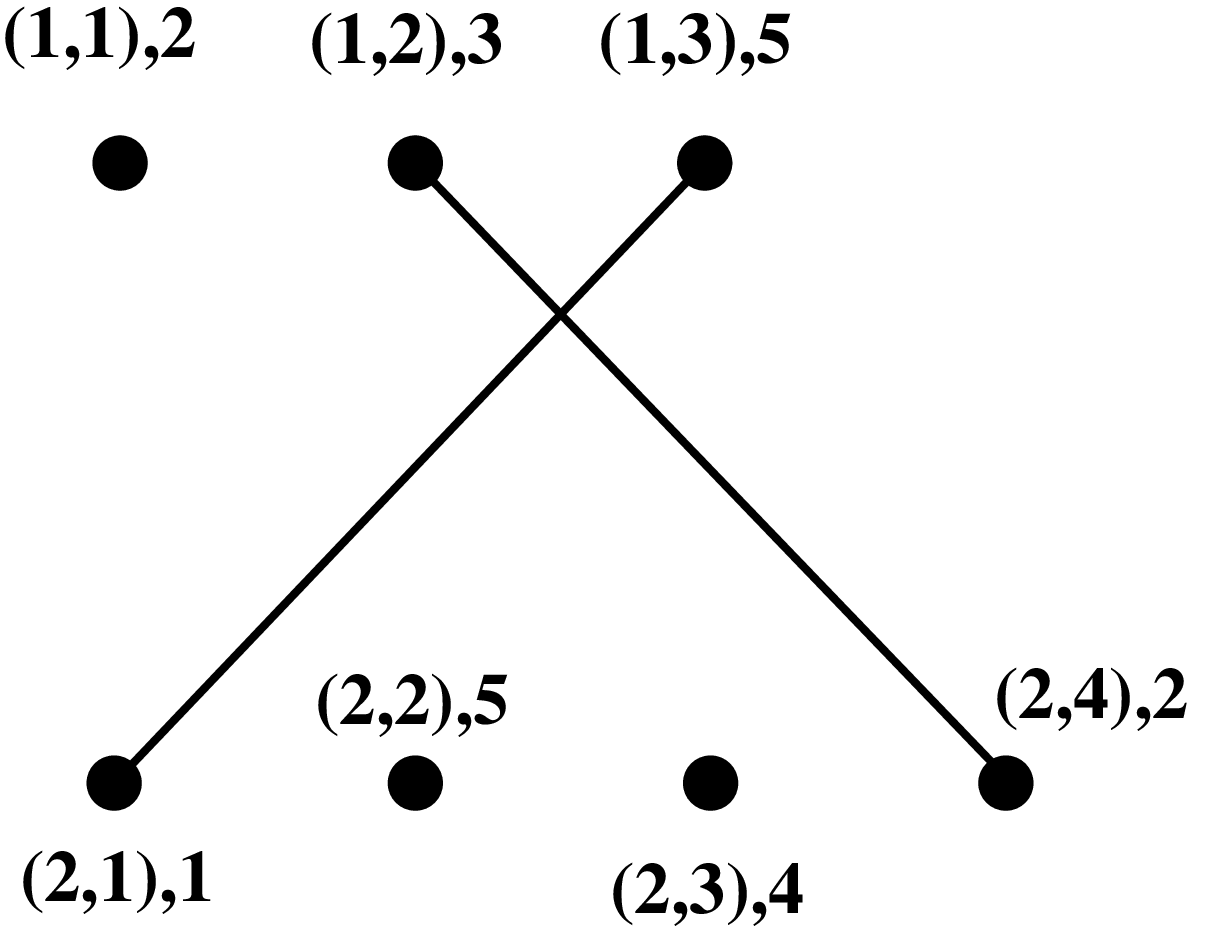, width=10cm}
\caption{the previous diagram and the enumeration of their vertices with
the help of the function $T_\gamma$}
\end{figure}

\medskip

To define the Wiener--It\^o integral corresponding to this
diagram let us first consider the function
$$
H(x_1,\dots,x_7)=h_1(x_1,x_2,x_3)h_2(x_4,x_5,x_6,x_7)
$$
defined  in~(\ref{6.3}). Simple calculation shows that
the function $\pi_\gamma(\cdot)= S(T_\gamma(\cdot))$ has the 
following form in this example. $\pi_\gamma(1)=1$, 
$\pi_\gamma(2)=5$, $\pi_\gamma(3)=6$, $\pi_\gamma(4)=2$, 
$\pi_\gamma(5)=3$, $\pi_\gamma(6)=7$, $\pi_\gamma(7)=4$.
This also means that the coordinates of the vectors 
in the Euclidean space ${\mathbb R}^7_\gamma$ which we get 
by reindexing the coordinates of the vectors in
${\mathbb R}^7$ have the form
$$
(x(\gamma)_1,x(\gamma)_2,x(\gamma)_3,x(\gamma)_4,x(\gamma)_5,
x(\gamma)_6,x(\gamma)_7)=(x_1,x_5,x_6,x_2,x_3,x_7,x_4).
$$
Then we can write the function $\bar H_\gamma$ and $\bar h_\gamma$ 
defined in~(\ref{6.12}) and~(\ref{6.13}) as
$$
H_\gamma(x(\gamma)_1,\dots,x(\gamma)_7)
=h_1(x(\gamma)_1,x(\gamma)_4,x(\gamma)_5)
h_2(x(\gamma)_2,x(\gamma)_3,x(\gamma)_6,x(\gamma)_7),
$$
and
$$
\bar h_\gamma(x(\gamma)_1,\dots,x(\gamma)_5)
=h_1(x(\gamma)_1,x(\gamma)_4,x(\gamma)_5)
h_2(x(\gamma)_2,x(\gamma)_3,-x(\gamma)_4,-x(\gamma)_5).
$$
Then we have
$$
\bar{\bar h}_\gamma(x(\gamma)_1,x(\gamma)_2,x(\gamma)_3)
=\int\bar h_\gamma(x(\gamma)_1,\dots,x(\gamma)_5)G_{3,2}(\,dx(\gamma)_4)
G_{5,1}(\,dx(\gamma)_5),
$$
and 
\begin{eqnarray*}
&&I_3(\bar{\bar h}_\gamma|2,5,4) \\
&&\qquad
=\int \bar{\bar h}_\gamma(x(\gamma)_1,x(\gamma)_2,x(\gamma)_3)
Z_{G,2}(\,dx(\gamma)_1)Z_{G,5}(\,dx(\gamma)_2)Z_{G,4}(\,dx(\gamma)_3)
\end{eqnarray*}
is the multiple Wiener--It\^o integral corresponding to
the diagram~$\gamma$ in the diagram formula. To understand
the definition of the function $\bar{\bar h}_\gamma$ and of the
Wiener--It\^o integral $I_3(\bar{\bar h}_\gamma)$ let us
observe that the first edge of the diagram connects the
vertices $(1,2)$ and $(2,4)$ with colours 3 and 2, hence
in the definition of $\bar{\bar h}_\gamma$ we integrate the 
argument $x(\gamma)_4$ by $G_{3,2}(\,dx(\gamma)_4)$, the 
second edge connects the vertices $(1,3)$ and $(2,1)$ 
with colours 5 and 1, hence we integrate the variable
$x(\gamma)_5$ by $G_{5,1}(\,dx(\gamma)_5)$. In the definition
of the Wiener integral the variable $x(\gamma)_1$ corresponds
to the vertex $S^{-1}(\pi_\gamma(1))=(1,1)$ which has colour 2,
hence we integrate the variable $x(\gamma)_1)$ by 
$Z_{G,2}(\,dx(\gamma)_1)$. Similarly, we define the variable
$x(\gamma)_2$ by the measure determined by the colour of 
$S^{-1}(\pi_\gamma(2))=(2,2)$, which is 5, i.e., we integrate 
by $Z_{G,5}(\,dx(\gamma)_2)$. Finally $S^{-1}(\pi_\gamma(3))=(2,3)$ 
has colour 4, and we integrate the variable $x(\gamma)_3$ 
by $Z_{G,4}(\,dx(\gamma)_3)$.

The Wiener--It\^o integral
$I_3(\bar{\bar h}_\gamma|3,1,3)$ can be rewritten with the
help of formula~(\ref{6.16}) in the following form:
$$
I_3(\bar{\bar h}_\gamma|2,5,4)
=I_3(h_\gamma|2,5,4)
=\int h_\gamma(x_1,x_2,x_3)
Z_{G,2}(\,dx_1)Z_{G,5}(\,dx_2)Z_{G,4}(\,dx_3)
$$
with
$$
h_\gamma(x_1,x_2,x_3)=\int h_1(x_1,x_4,x_5)h_2(x_2,x_3,-x_4,-x_5)
G_{3,2}(\,dx_4) G_{5,1}(\,dx_5).
$$
This expression can be calculated similarly to
$I_3(\bar{\bar h}_\gamma|2,5,4)$, only we have to replace
$x(\gamma)_s$ everywhere by $x_s$ in the calculation.

\medskip\medskip
I formulate a corollary of the diagram formula in which 
I consider the special case of this result when the second 
Wiener--It\^o integral defined in formula~(\ref{6.2}) is a 
one-fold integral. In this case it has the simpler form
\begin{equation}
I_1(h_2|j'_1)=\int h_2(x_1)Z_{G,{j'_1}}(\,dx_1) 
\quad\textrm{with }h_2\in{\cal K}_{1,j'_1}.  \label{6.18}
\end{equation}
Here again we formulate the problem in the following way. We 
take a pair of functions $h_1(x_1,\dots,x_n)$ and $h_2(x_{n+1})$ 
on ${\mathbb R}^{(n+1)\nu}$. Then we define a function 
$h_2^{(0)}(x_1)$ on ${\mathbb R}^1$ by the formula 
$h^{(0)}_2(x_1)=h_2(x_{n+1})$ if $x_1=x_{n+1}$. We integrate
the function $h_2^{(0)}(x)$ in formula~(\ref{6.18}), but we 
omit the superscript~${}^{(0)}$ in our notation. We assume 
that $h_1\in{\cal K}_{n,j_1,\dots,j_n}$, and $h_2\in{\cal K}_{1,j'_1}$. 

In the next Corollary I express the product of the 
Wiener--It\^o integrals given in~(\ref{6.1}) and (\ref{6.18}) 
as a sum of Wiener--It\^o integrals. This formula will be 
needed in the proof of the multivariate version of It\^o's 
formula in paper~\cite{11}. 

The diagram formula in this case has a simpler form, since the 
second row of the diagrams we are working with consists only of one 
point $(2,1)$. Hence there are only the diagram $\gamma_0\in\Gamma$ 
that contains no edges and the diagrams $\gamma_p\in\Gamma$, 
$1\le p\le n$, which contain one edge that connects the vertices 
$(1,p)$ and $(2,1)$. 

\medskip\noindent
{\bf Corollary of Theorem 6.1.} {\it The product of the Wiener--It\^o 
integrals 
$$
I_n(h_1|j_1,\dots,j_n) \textrm{ \ and \ } I_1(h_2|j_1')
$$ 
introduced in formulas (\ref{6.1}) and (\ref{6.18}) satisfy the identity
\begin{eqnarray}
&&I_n(h_1|j_1,\dots,j_n)I_1(h_2|j_1') \label{6.19} \\
&&\qquad =\int h_{\gamma_0}(x_1,\dots,x_{n+1}) 
Z_{G,{j_1}}(\,dx_1)\cdots  Z_{G,{j_n}}(\,dx_n) 
Z_{G,{j'_1}}(\,dx_{n+1}) \nonumber \\
&&\qquad\qquad +\sum_{p=1}^n \int h_{\gamma_p}(x_1,\dots,x_{n-1})
\prod_{s=1}^{p-1}Z_{G,{j_s}}(\,dx_s) \prod_{s=p}^{n-1} 
Z_{G,{j_{s+1}}}(\,dx_s) \nonumber \\
&&\qquad=I_{n+1}(h_{\gamma_0}|j_1,\dots,j_n, j_1')
+\sum_{p=1}^n 
I_{n-1}(h_{\gamma_p}|j_1,\dots,j_{p-1},j_{p+1},\dots,j_n), \nonumber 
\end{eqnarray}
where $h_{\gamma_0}(x_1,\dots,x_{n+1})=h_1(x_1,\dots,x_n)h_2(x_{n+1})$,
and for $1\le p\le n$  
$$
h_{\gamma_p}(x_1,\dots,x_{n-1})=
\int h_{1,\gamma_p}(x_1,\dots,x_n)\overline{h_2(x_n)}G_{j_p,j'_1}(\,dx_n)
$$
with
$h_{1,\gamma_p}(x_1,\dots,x_n)
=h_1(x_{\pi_p(1)},\dots,x_{\pi_p(n)})$,
where $\pi_p(k)=k$ if $1\le k\le p-1$, $\pi_p(p)=n$,
and $\pi_p(k)=k-1$ if $p+1\le k\le n$.

To make the definition of formula (\ref{6.19}) complete I remark that
for $p=1$ we put $\prod\limits_{s=1}^0Z_{G,{j_s}}(\,dx_s)\equiv1$ and for
$p=n$  $\prod\limits_{s=n}^{n-1} Z_{G,{j_s}}(\,dx_s)\equiv1$.}

\medskip\noindent
{\it Proof of the Corollary.}\/ We get the result of the corollary 
by applying Theorem~6.1 in the special case when the second  
Wiener--It\^o integral is defined by formula~(\ref{6.18}) 
instead of~(\ref{6.2}). We have to check that in this case 
the function $h_{\gamma_0}$ corresponding to the diagram 
$\gamma_0$ agrees with the function $h_{\gamma_0}$ defined in 
the corollary, and to calculate the functions $h_{\gamma_p}$ 
defined in~(\ref{6.14}) for the remaining 
diagrams~$\gamma_p$, $1\le p\le n$. In this case 
$\pi_{\gamma_p}(k)=k$ for $1\le k\le p-1$, 
$\pi_{\gamma_p}(k)=k+1$ for $p\le k\le n-1$, 
$\pi_{\gamma_p}(n)=p$, $\pi_{\gamma_p}(n+1)=n+1$, hence
$$
(x(\gamma_p)_1,\dots,x(\gamma_p)_{n+1})
=(x_1,\dots,x_{p-1},x_{p+1},\dots,x_n,x_p,x_{n+1}),
$$
and 
$$
\bar h_{\gamma_p}(x(\gamma_p)_1,\dots,x(\gamma_p)_{n+1})=
h_1(x(\gamma_p)_1,\dots,x(\gamma_p)_n) h_2(-x(\gamma_p)_n)
$$
for $1\le p\le n$. On the other hand, 
$h_2(-x)=\overline{h_2(x)}$, since $h_2\in{\cal K}_{1,j_1'}$.
Thus
\begin{eqnarray*}
&&\bar{\bar h}_{\gamma_p}(x(\gamma_p)_1,\dots,x(\gamma_p)_{n-1}) \\
&&\qquad=
\int h_1(x(\gamma_p)_1,\dots,x(\gamma_p)_{n-1},x(\gamma_p)_n) 
\overline{h_2(x(\gamma_p)_n)}G_{j_p,j'_1}(\,dx(\gamma_p)_n).
\end{eqnarray*}
Then simple calculation shows that for $\gamma=\gamma_p$
the kernel function $h_\gamma=h_{\gamma_p}$ in formula~(\ref{6.16})
agrees with the function~$h_{\gamma_p}$ defined in the corollary
for all $1\le p\le n$, and Theorem~6.1 yields 
identity~(\ref{6.19}) under the conditions of the~corollary.
The corollary is proved.

\medskip
The proof of Theorem~6.1 is similar to the proof of the diagram 
formula (Theorem~5.3 in~\cite{9}). It applies the same method, 
only the notation becomes more complicated than the also rather 
complicated notation of the original proof, since we have to 
work with spectral measures of the form $G_{j_s,j'_t}$ and random 
spectral measures of the form $Z_{G,{j_s}}$ or $Z_{G,{j'_t}}$ instead 
of the  spectral measure $G$ and random spectral measure $Z_G$. 
Hence I decided not to describe the complete proof, I only 
concentrate on its main ideas and the formulas that explain 
why such a result appears in the diagram formula. The
interested reader can reconstruct the proof by means of a 
careful study of the proof of Theorem~5.3 in~\cite{9}.

\medskip\noindent
{\it A sketch of proof for Theorem 6.1.} The proof of Part~A 
is relatively simple. One can check that the function $h_\gamma$ 
satisfies relation~(a) in the definition of the functions in
${\cal K}_{n+m-2|\gamma|,j_{r_1},\dots,j_{r_{n-|\gamma|}},
j'_{t_1},\dots,j'_{t_{m-|\gamma|}}}$
given in Section~5 by exploiting formula~(\ref{6.14}), the similar 
property of the functions $h_1$ and $h_2$ together with the 
symmetry property $G_{j,j'}(-A)=\overline{G_{j,j'}(A)}$ for all 
$1\le j,j'\le d$ and sets~$A$ of the spectral measure~$G$.

To prove the inequality formulated in Part~A let us first 
rewrite the definition of $h_\gamma$ in~(\ref{6.14}) by replacing all 
measures of the form $G_{j.j'}(dx)$ by 
$g_{j,j'}(x)\mu(\,dx)=G_{j,j'}(\,dx)$, where $\mu$ is a 
dominating measure for all complex measures 
$G_{j,j'}$, $g_{j,j'}$ is the Radon--Nikodym derivative of 
$G_{j,j'}$ with respect to~$\mu$, and observe that the 
inequality~(\ref{3.2}) and formula~(\ref{6.13}) 
and~(\ref{6.14}) imply that
\begin{eqnarray*}
&&|\bar{\bar h}_\gamma(x(\gamma)_1,\dots,x(\gamma)_{n+m-2|\gamma|})| \\
&&
\qquad\le\int 
h_1(x_{\pi_\gamma(1)},\dots,x_{\pi_\gamma(n-|\gamma|)},
x_{\pi_\gamma(n+m-2|\gamma|+1)},\dots,x_{\pi_\gamma(n+m-|\gamma|+1)})  \\
&& 
\qquad\qquad \times h_2(x_{\pi_\gamma(n-|\gamma|+1)},\dots, 
x_{\pi_\gamma(n+m-2|\gamma|)},  \\ 
&&\qquad\qquad\qquad\qquad -x_{\pi_\gamma(n+m-|2\gamma|+1)},\dots, 
-x_{\pi_\gamma(n+m-|\gamma|)})   \\
&&\qquad\qquad\quad \times \prod_{k=1}^{|\gamma|} 
\sqrt{g_{j_{v_k},j_{v_k}}(x_{\pi_\gamma(n+m-2|\gamma|+k)})} 
\sqrt{g_{j'_{w_k},j'_{w_k}}(x_{\pi_\gamma(n+m-2|\gamma|+k)})} \\ 
&&\qquad\qquad\qquad\qquad
\times \mu(\,dx_{\pi_\gamma(n+m-2|\gamma|+k)}).
\end{eqnarray*}
We get, by applying the Schwarz inequality the evenness of the
measures $G_{j,j}$ and by 
replacing the measures of the form $g_{j,j}(x)\mu(\,dx)$ 
or $g_{j',j'}(x)\mu(\,dx)$ by the measures of the form 
$G_{j,j}(\,dx)$ and $G_{j',j'}(\,dx)$ that 
\begin{eqnarray*}
&&|\bar{\bar h}_\gamma(x(\gamma)_1,\dots,x(\gamma)_{n+m-2|\gamma|})|^2 \\
&&\qquad\le\int
|h_1(x_{\pi_\gamma(1)},\dots,x_{\pi_\gamma(n-|\gamma|)},
x_{\pi_\gamma(n+m-2|\gamma|+1)},\dots,x_{\pi_\gamma(n+m-|\gamma|+1)})|^2  \\
&&\qquad\qquad\qquad \times \prod_{k=1}^{|\gamma|}
G_{j_{v_k},j_{v_k}}(\,dx_{\pi_\gamma(n+m-2|\gamma|+k)})\\ 
&&\qquad\qquad \times \int 
|h_2(x_{\pi_\gamma(n-|\gamma|+1)},\dots, 
x_{\pi_\gamma(n+m-2|\gamma|)},  \\ 
&&\qquad\qquad\qquad\qquad\qquad -x_{\pi_\gamma(n+m-|2\gamma|+1)},\dots, 
-x_{\pi_\gamma(n+m-|\gamma|)})|^2\\
&&\qquad\qquad\qquad\qquad \times \prod_{k=1}^{|\gamma|}
G_{j'_{w_k,w_k}} (\,dx_{\pi_\gamma(n+m-2|\gamma|+k)}).
\end{eqnarray*}
Let us integrate the last inequality with respect to the
product measure
\begin{eqnarray*} 
&&\prod_{k=1}^{n-|\gamma|}G_{j_{r_k},j_{r_k}}(\,dx(\gamma)_k)
\prod_{l=1}^{m-|\gamma|} G_{j'_{t_l},j'_{t_l}}(\,dx(\gamma)_{n-|\gamma|+l}) \\
&&\qquad=\prod_{k=1}^{n-|\gamma|}G_{j_{r_k},j_{r_k}}(\,dx_{\pi_\gamma(k)})
\prod_{l=1}^{m-|\gamma|}G_{j'_{t_l},j'_{t_l}}(\,dx_{\pi_\gamma(n-|\gamma|+l)}).
\end{eqnarray*}
A careful analysis shows that the inequality we get in such a
way agrees with the inequality formulated in Part~A of Theorem~6.1.
Indeed, we get at the left-hand side of this inequality 
$\|\bar{\bar h}_\gamma\|$ with the norm formulated in Part~A of 
Theorem~6.1, and the right-hand side equals the product 
$\|h_1\|\|h_2\|$. We got the same integrals as the integrals 
defining these norms, only we integrate by the variables of 
the functions $h_1$ and $h_2$ in a different order. We also 
have to exploit that the measures $G_{j,j}$ are symmetric, 
hence the value of the integrals we are investigating does 
not change if we replace the coordinate $x_k$ by $-x_k$ in 
the kernel function for certain coordinates $k$.

Next I turn to the proof of Part~B of Theorem~6.1. First we 
prove this result, i.e., identity~(\ref{6.17}) in the special 
case when both $h_1$ and $h_2$ are simple functions. We may 
also assume that they are adapted to the same regular system
$$
{\cal D}=\{\Delta_p,\;p=\pm1,\pm2,\dots,\pm N\},
$$ 
and by a possible further division of the sets $\Delta_p$ we may 
also assume that the elements of ${\cal D}$ are very small. More 
explicitly, first we choose such a measure $\mu$ on ${\mathbb R}^\nu$ 
which has finite value on all compact sets, all complex  
measures $G_{k,l}$, $1\le k,l\le d$, are absolutely continuous with 
respect to $\mu$, and their Radon--Nikodym derivatives satisfy 
the inequality $|\frac{dG_{k,l}}{d\mu}(x)|\le1$ for all 
$x\in {\mathbb R}^\nu$. Fix a small number~$\varepsilon>0$. We may achieve, 
by splitting up the sets $\Delta_p$ into smaller sets if it is 
necessary, that $\mu(\Delta_p)\le\varepsilon$ for all 
$\Delta_p\in{\cal D}$.
Let us fix a number $u_p\in\Delta_p$ in all sets 
$\Delta_p\in{\cal D}$. We can express the product 
$I_n(h_1|j_1,\dots,j_n)I_m(h_2|j_1',\dots,j'_m)$ as
\begin{eqnarray*}
&& \!\!\!\!\!\!\!\!\!\!
I=I_n(h_1|j_1,\dots,j_n)I_m(h_2|j_1',\dots,j'_m)
={\sum}' h_1(u_{p_1},\dots,u_{p_n})h_2(u_{q_1},\dots,u_{q_m})\\
&&\qquad\qquad \times Z_{G,{j_1}}(\Delta_{p_1})\cdots Z_{G,{j_n}}(\Delta_{p_n})
Z_{G,{j'_1}}(\Delta_{q_1})\cdots Z_{G,{j'_m}}(\Delta_{q_m}).
\end{eqnarray*}
The summation in the sum $\sum'$ goes through all pairs 
$((p_1,\dots,p_n),(q_1,\dots,q_m))$ such that 
$p_k,\,q_l\in\{\pm1,\dots,\pm N\}$, \ $k=1,\dots,n$, $l=1,\dots,m$,
and $p_k\neq\pm p_{\bar k}$, if $k\neq \bar k$, and 
$q_l\neq\pm q_{\bar l}$ if $l\neq\bar l$.

Write
\begin{eqnarray*}
I&=&\sum_{\gamma\in\Gamma}
{\sum}^\gamma \,
h_1(u_{p_1},\dots,u_{p_n}) h_2(u_{q_1},\dots,u_{q_m})  \\
&&\qquad\qquad \times Z_{G,{j_1}}(\Delta_{p_1})\cdots Z_{G,{j_n}}(\Delta_{p_n})
Z_{G,{j'_1}}(\Delta_{q_1})\cdots Z_{G,{j'_n}}(\Delta_{q_m}).
\end{eqnarray*}
where $\sum^\gamma$ contains those terms of $\sum'$ for which
$p_k=q_l$ or $p_k=-q_l$ if the vertices $(1,k)$ and $(2,l)$ are
connected in $\gamma$, and $p_k\neq \pm q_l$ if $(1,k)$ and $(2,l)$
are not connected in~$\gamma$.

Let us introduce the notation
\begin{eqnarray*}
\Sigma^\gamma&=&
{\sum}^\gamma \,
h_1(u_{p_1},\dots,u_{p_n}) h_2(u_{q_1},\dots,u_{q_m})  \\
&&\qquad\qquad\times Z_{G,{j_1}}(\Delta_{p_1})\cdots Z_{G,{j_n}}(\Delta_{p_n})
Z_{G,{j'_1}}(\Delta_{q_1})\cdots Z_{G,{j'_n}}(\Delta_{q_m}).
\end{eqnarray*}
for all $\gamma\in\Gamma$.

We want to show that for small $\varepsilon>0$ (where $\varepsilon$
is an upper bound for the measure $\mu$ of the sets 
$D_p\in{\cal D}$) the expression $\Sigma^\gamma$ is very close to
\begin{equation}
I_\gamma=I_{n+m-2|\gamma|}(\bar{\bar h}_\gamma|j_{v_1},\dots,j_{v_{(n-|\gamma|}},
j'_{w_1},\dots,j'_{w_{m-|\gamma|}})  \label{6.20}
\end{equation}
for all $\gamma\in\Gamma$. For this goal we make the decomposition
$\Sigma^\gamma=\Sigma^\gamma_1+\Sigma^\gamma_2$ of $\Sigma^\gamma$ with
\begin{eqnarray*}
\Sigma_1^\gamma&=&{\sum}^\gamma\, h_1(u_{p_1},\dots,u_{p_n})
h_2(u_{q_1},\dots,u_{q_m})
\prod_{k\in A_1}Z_{G,{j_k}}(\Delta_{p_k})
\prod_{l\in A_2} Z_{G,{j'_l}}(\Delta_{q_l})\\
&&\qquad\qquad\qquad
\times \prod_{(k,l)\in B} 
E\left(Z_{G,{j_k}}(\Delta_{p_k})Z_{G,{j'_l}}(\Delta_{q_l})\right)
\end{eqnarray*}
and
$$
\Sigma^\gamma_2=\Sigma^\gamma-\Sigma^\gamma_1,
$$
where the sets $A_1$, $A_2$ and $B$ were defined in 
formulas~(\ref{6.6}), (\ref{6.7}) and~(\ref{6.8}).

It is not difficult to check that both $\Sigma^\gamma_1$ and 
$\Sigma^\gamma_2$ are real valued random variables. We want to
show that $\Sigma^\gamma_1$ is close to the random variable 
$I_\gamma$ introduced in~(\ref{6.20}), while $\Sigma^\gamma_2$ is a 
small error term. To understand the behaviour of 
$\Sigma^\gamma_1$ observe that
$$
E(Z_{G,{j_k}}(\Delta_{p_k})Z_{G,{j'_l}}(\Delta_{q_l})
=E(Z_{G,{j_k}}(\Delta_{p_k})\overline{Z_{G,{j'_l}}(-\Delta_{q_l})}=0
$$
if $\Delta_{p_k}=\Delta_{q_l}$
(and as a consequence
if $\Delta_{p_k}\cap(-\Delta_{q_l})=\emptyset$), and
$$
E(Z_{G,{j_k}}(\Delta_{p_k})Z_{G,{j'_l}}(\Delta_{q_l})
=E(Z_{G,{j_k}}(\Delta_{p_k})\overline{Z_{G,{j'_l}}(-\Delta_{q_l})}
=G_{j_k,j'_l}(\Delta_{p_k})
$$
if $\Delta_{p_k}=-\Delta_{q_l}$. In the case $(k,l)\in B$ one 
of these possibilities happens.

These relations make possible to rewrite $\Sigma^\gamma_1$ 
in a simpler form. It can be rewritten in the form of a 
Wiener--It\^o integral of order $n+m-2|\gamma|$ with 
integration with respect to the  random measure 
$\prod\limits_{k\in A_1}Z_{G,{j_k}}(\,dx_k)
\prod\limits_{l\in A_2} Z_{G,{j'_l}}(\,dx_l)$, (where the sets
$A_1$ and $A_2$ were defined in~(\ref{6.6}) and~(\ref{6.7})).
Then we can rewrite this integral, by reindexing its 
variables in a right way to an integral very similar 
to the Wiener--It\^o integral~(\ref{6.15}) (with the same 
parameter $\gamma$). The difference between these 
two expressions is that the kernel function $h'_\gamma$ 
of the Wiener--It\^o integral expressing 
$\Sigma^\gamma_1$ is slightly different from the kernel 
function $\bar{\bar h}_\gamma$ appearing in the other integral. 
The main difference between these two kernel functions
is that there is a small set in the domain of 
integration where $h'_\gamma$ disappears, while 
$\bar{\bar h}_\gamma$ may not disappear. But the two 
Wiener--It\^o integrals are very close to each other. An 
adaptation of the argument in the proof of Theorem~5.3 
in~\cite{9} shows that
$$
E(\Sigma_1^\gamma-I_\gamma)^2\le C\varepsilon
$$
with an appropriate constant $C>0$.

We also want to show that $\Sigma^\gamma_2$ is a 
negligibly small error term. To get a good upper 
bound on $E(\Sigma^\gamma_2)^2$ we write it in the form
\begin{eqnarray*}
E(\Sigma_2^\gamma)^2&=&{\sum}^{\gamma}_2\, 
h_1(u_{p_1},\dots,u_{p_n}) h_2(u_{q_1},\dots,u_{q_m}) \\
&&\qquad \times h_1(u_{\bar p_1},\dots,u_{\bar p_n}) 
h_2(u_{\bar q_1},\dots,u_{\bar q_m}) \\
&&\qquad\qquad \times \Sigma^\gamma_3(p_k,q_l,p_{\bar k},q_{\bar l},
\;k,\bar k\in\{1,\dots,n\},
\; l,\bar l\in\{1,\dots,m\}) 
\end{eqnarray*}
with
\begin{eqnarray*}
&& \!\!\!\!\!\!\!\!\!
\Sigma^\gamma_3(p_k,q_l,p_{\bar k},q_{\bar l},
\;k,\bar k\in\{1,\dots,n\},
\; l,\bar l\in\{1,\dots,m\}) \\
&&=E\Biggl(\left(\prod_{k\in A_1}Z_{G,{j_k}}(\Delta_{p_k})
\prod_{l\in A_2}Z_{G,{j'_l}}(\Delta_{q_l})
\prod_{\bar k\in A_1}Z_{G,{j_{\bar k}}}(\Delta_{p_{\bar k}})
\prod_{\bar l\in A_2}Z_{G,{j'_{\bar l}}}(\Delta_{q_{\bar l}})\right)\  \\
&&\qquad\times\left[\prod_{(k,l)\in B} 
Z_{G,{j_k}}(\Delta_{p_k})Z_{G,{j'_l}}(\Delta_{q_l})
-E\prod_{(k,l)\in B} Z_{G,{j_k}}(\Delta_{p_k})Z_{G,{j'_l}}(\Delta_{q_l})\right] \\
&&\qquad\times \left[\prod_{(\bar k,\bar l)\in B} 
Z_{G,{j_{\bar k}}}(\Delta_{p_{\bar k}}) Z_{G,{j'_{\bar l}}}(\Delta_{q_{\bar l}})
-E\prod_{(\bar k,\bar l)\in B} Z_{G,{j_{\bar k}}}(\Delta_{p_{\bar k}})
Z_{G,{j'_{\bar l}}}(\Delta_{q_{\bar l}})\right]\Biggr), 
\end{eqnarray*}
where we sum in ${\sum}^\gamma_2$ for such sequences of indices 
$p_k,\,q_l,\,p_{\bar k},\,q_{\bar l}$,  
$k,\bar k\in\{1,\dots,n\}$, $l,\bar l\in\{1,\dots,m\}$,
$p_k,p_{\bar k},q_l,q_{\bar l}\in\{\pm1,\dots,\pm N\}$
which satisfy the following properties. For all indices 
$k,l,\bar k$ and $\bar l$,  $p_k=q_l$ or $p_k=-q_l$ if 
$(k,l)\in B$, and similarly $p_{\bar k}=q_{\bar l}$ or 
$p_{\bar k}=-q_{\bar l}$ if $(\bar k,\bar l)\in B$. 
Otherwise all numbers $\pm p_k$ and $\pm q_l$ are
different, and similarly otherwise all $\pm p_{\bar k}$
and $\pm q_{\bar l}$ are different. 

We get a good estimate on $E(\Sigma_2^\gamma)^2$ by giving a good
bound on all terms
\begin{equation}
\Sigma^\gamma_3(p_k,q_l,p_{\bar k},q_{\bar l},
\;k,\bar k\in\{1,\dots,n\},\; l,\bar l\in\{1,\dots,m\}) \label{6.21}
\end{equation}
in the formula expressing it.
This can be done by adapting the corresponding argument in
Theorem~5.3 of~\cite{9}. This argument shows that for most
sets of parameters $p_k,q_k,p_{\bar k},q_{\bar l}$ the term
in~(\ref{6.21}) equals zero. More explicitly, it is
equal to zero if ${\cal A}\neq-\bar{\cal A}$ with
$$
{\cal A}=\{p_k\colon\; k\in A_1\}\cup \{q_l\colon\;l\in A_2\} \quad
\textrm{and} \quad \bar{\cal A}=\{p_{\bar k}\colon\; 
\bar k\in A_1\}\cup \{q_{\bar l}\colon\;\bar l\in A_2\},
$$
and it also equals zero if ${\cal F}\cup(-{\cal F})$ and 
$\bar {\cal F}\cup(-\bar {\cal F})$ are disjoint, where
$$
{\cal F}=\bigcup_{(k,l)\in B}\{p_k,q_l\} \quad\textrm{and}\quad
\bar{\cal F}=\bigcup_{(\bar k,\bar l)\in B}\{p_{\bar k},q_{\bar l}\}.
$$

These statements can be proved by adapting the corresponding 
argument in Theorem~5.3 of~\cite{9}. More precisely, in the proof 
of the first statement we still need the following additional 
observation. If $(X,Y,Z)$ is a three-dimensional Gaussian 
vector with $EX=EY=EZ=0$, then $EXYZ=0$. (In the proof of 
Theorem~5.3 in~\cite{9} we needed this statement only in a 
special case when it trivially holds.)

To prove this statement let us apply the following orthogonalization
for the random variables $X$, $Y$ and $Z$. Write $Y=\alpha X+\eta$,
$Z=\beta_1 X+\beta_2\eta+\zeta$, where $X,\eta,\zeta$ are orthogonal,
(jointly) Gaussian random variables with expectation zero. Then they
are also independent, hence 
$EXYZ=EX(\alpha X+\eta)(\beta_1X+\beta_2\eta+\zeta)=0$.

In the remaining cases the expression in~(\ref{6.21}) can be estimated 
(again by adapting the argument of Theorem~5.3 in~\cite{9}) in the 
following way.
\begin{eqnarray*}
&&\Sigma^\gamma_3(p_k,q_l,p_{\bar k},q_{\bar l},
\;k,\bar k\in\{1,\dots,n\},\; l,\bar l\in\{1,\dots,m\}) \\
&&\qquad \le C\varepsilon {\prod}'
\mu(\Delta_{p_k})\mu(\Delta_{l_q})
\mu(\Delta_{p_{\bar k}})\mu(\Delta_{q_{\bar l}})
\end{eqnarray*}
with some constant $C$ (not depending on~$\varepsilon$) and the
measure $\mu$ dominating the complex measures 
$G_{j,k}$ with the properties we demanded at the start of 
the proof. The sign~$'$ in the product ${\prod}'$ means that 
first we take the sets
$\Delta_{p_k}$, $\Delta_{q_l}$, $\Delta_{p_{\bar k}}$, $\Delta_{q_{\bar l}}$
for all parameters $k,\bar k\in\{1,\dots,n\}$ and 
$l,\bar l\in\{1,\dots,m\}$, then if a set $\Delta$ appears twice 
in the sequence of these sets we omit one of them. Then if both 
the sets $\Delta$ and $-\Delta$ appear for some set $\Delta$, 
then we omit one of them from this sequence. Then we take 
in~${\prod}'$ the product of the terms~$\mu(\Delta)$ with
the sets~$\Delta$ in the remaining sequence.

It can be proved with the help of the estimates on the terms 
in~(\ref{6.21}) (see again Theorem~5.3 in~\cite{9}) that
$$
E(\Sigma_2^\gamma)^2\le C\varepsilon.
$$
It is not difficult to prove part~B of Theorem~6.1 with the help of 
the estimates on $E(\Sigma_1^\gamma-I_\gamma)^2\le C\varepsilon$ and 
$E(\Sigma_2^\gamma)^2\le C\varepsilon$ if $h_1$ and $h_2$ are simple
functions. One only has to make an appropriate limiting procedure
with $\varepsilon\to0$. Then we can complete the proof of Theorem~6.1
similarly to the proof of Theorem~5.3 in~\cite{9} by means of an 
appropriate approximation of Wiener--It\^o integrals with 
Wiener--It\^o integrals of simple functions. In this approximation 
we have to apply Lemma~5.1 and the properties of the Wiener--It\^o
integrals, in particular the already proved Part~A of Theorem~6.1.

\end{document}